\documentclass{article}

\usepackage[english]{babel}

\usepackage[letterpaper,top=2cm,bottom=2cm,left=3cm,right=3cm,marginparwidth=1.75cm]{geometry}

\usepackage{amsmath} 
\usepackage{bm}
\usepackage{amsfonts}
\usepackage{graphicx}
\usepackage[colorlinks=true, allcolors=blue]{hyperref}
\usepackage{multirow}
\usepackage{verbatim}
\usepackage{float}
\usepackage{booktabs}
\usepackage{hyperref}
\usepackage{authblk}
\usepackage{subcaption}
\usepackage{physics}

\usepackage{tikz}
\usepackage{bm}       
\usetikzlibrary{shapes.geometric, arrows.meta, positioning, fit, calc, decorations.pathreplacing, fit, backgrounds, shapes.misc}
\definecolor{fomBlue}{HTML}{00AEEF}
\definecolor{offlineOrange}{HTML}{F7941D}
\definecolor{nnGreen}{HTML}{8DC63F}
\definecolor{onlinePurple}{HTML}{9E005D}
\definecolor{darkGray}{HTML}{555555}

\date{}

\title{Machine-Learning-Enhanced Discretize-then-Project Reduced-Order Modeling of Turbulent Flows on Collocated Grids}

\author[1]{Nadim Rooholamin\footnote{nadim.r99@gmail.com}}
\author[1,2]{Kabir Bakhshaei\footnote{kabir.bakhshaei@santannapisa.it}}
\author[1]{Giovanni Stabile\footnote{giovanni.stabile@santannapisa.it}}

\affil[1]{Biorobotics Institute, SMART Lab, Sant'Anna School of Advanced Studies, Pisa, Italy}
\affil[2]{Department of Computer Science, University of Pisa, Pisa, Italy}

\begin{document}
\date{}
\maketitle

\begin{abstract}

This study presents a hybrid reduced-order modeling (ROM) framework for incompressible flows on collocated finite-volume grids. The method combines a discretize-then-project consistent-flux formulation for velocity and pressure with a non-intrusive neural-network closure for the turbulent-viscosity coefficients. The intrusive formulation preserves the discrete mass-conservation constraint and pressure--velocity coupling, while a reduced pressure reference-cell constraint fixes the pressure gauge ambiguity. Three neural-network architectures are evaluated for the closure: a Multilayer Perceptron (MLP), a Transformer, and a Long Short-Term Memory (LSTM) network. The framework is first assessed on a three-dimensional lid-driven cavity simulation at $Re=100$, where the LSTM-based ROM achieves relative errors of 0.7\% for velocity and 4\% for turbulent viscosity. A more demanding case at $Re=3200$ is then considered. A systematic mode-sensitivity study identifies $N=15$ POD modes as the best overall configuration, balancing velocity and pressure accuracy, reduced dimension, closure robustness, and computational cost. For this configuration, the LSTM-based ROM achieves a final relative velocity error of approximately 12.3\% and an online-phase wall-clock speedup of approximately $50\times$ over the full-order model. Energy and enstrophy errors remain below 11\% across all three evaluated architectures. Application to this second regime requires case-specific neural-network retraining and pressure reference-cell parameter retuning. A time-extrapolation test, using training data limited to $t\in[0,3]$\,s and a ROM rollout to $t=6$\,s, shows that the ROM remains bounded over the investigated extrapolation interval. Velocity and pressure errors increase beyond the training window, while the LSTM turbulent-viscosity closure remains robust, identifying long-horizon pressure accuracy as the primary current limitation. These results demonstrate the potential of combining a consistent projection-based formulation with a data-driven turbulence closure for efficient reduced-order simulation.

\end{abstract}


\maketitle

\section{Introduction}     
Fluid dynamics represents one of the most intricate and continuously evolving disciplines, characterized by ongoing advancements and a wide spectrum of theoretical investigations and practical applications. It constitutes a branch of fluid mechanics devoted to the analysis of the behavior of liquids and gases in motion. The resolution of a fluid dynamics problem typically entails solving, either analytically or numerically, highly complex differential equations in order to evaluate physical quantities such as velocity, pressure, density, or temperature as functions of spatial coordinates and time.
At present, numerical simulations, most notably Computational Fluid Dynamics (CFD), have become the standard approach. CFD is now routinely used for the simulation of complex flow problems; however, its application is often restricted by the considerable computational cost and memory requirements associated with high-resolution models. This issue becomes more severe when turbulent flows are considered, or when
the geometry and physics of the problem increase in complexity, as in
time-resolved multiphysics simulations involving moving geometries and
coupled fluid--multibody dynamics
\cite{Bakhshaei2021FlappingWings}. For this reason, reduced-order modeling (ROM) approaches have been widely
investigated as an alternative strategy to reduce computational effort while
still capturing the main flow behavior. The objective of these methods is to
decrease the dimensionality of high-fidelity models, thereby mitigating
computational costs while retaining physical accuracy
\cite{Benner2015,Rozza202059}. Reduced-order methods can be broadly
classified into projection-based and non-projection-based approaches, the
latter including truncation methods, goal-oriented techniques, and
low-dimensional models constructed from input--output data
\cite{Willcox2002,ROWLEY2005,BuiThanh2007,Ravindran2000}.

Beyond conventional forward simulations, computational efficiency is also
critical in inverse and data-assimilation problems, where the governing model
must be evaluated repeatedly to estimate uncertain parameters, states, or
boundary conditions. Stochastic data-assimilation techniques have, for
example, been applied to estimate noisy, time-dependent velocity boundary
profiles in patient-specific cardiovascular flows and unknown heat-flux
boundary conditions in inverse heat-transfer problems
\cite{Bakhshaei2025Cardiovascular,Bakhshaei2026InverseHeat}.

Within this framework, projection-based methods achieve dimensionality reduction by projecting the governing partial differential equations (PDEs) onto a lower-dimensional subspace, referred to as the reduced basis (RB) \cite{Rozza2008, Veroy2003,Hesthaven2016}. This results in a compact model that remains consistent with the governing physical laws. Several strategies are commonly employed to construct the reduced basis, including Dynamic Mode Decomposition (DMD), greedy selection procedures, and Proper Orthogonal Decomposition (POD)~\cite{Prudhomme2001,SCHMID2010,Kutz2016,Tissot2014, hajisharifi2024comparison}. These approaches enable the identification of low-dimensional subspaces that retain the dominant features of the original high-fidelity solution. 
To alleviate the computational effort required by nonlinear operators,
several ROM techniques rely on hyper-reduction methods, including gappy
POD and discrete empirical interpolation
\cite{Holmes2012,Sirovich1987,Galbally2009,Chaturantabut2010,
Amsallem2012}; deep-learning-enhanced variants of Gappy POD have also
been developed for sparse full-field reconstruction in conjugate
heat-transfer applications
\cite{hajisharifi2026deep}.



A classical example of a projection-based technique is the 
POD--Galerkin method \cite{Quarteroni2016}. Initially, POD extracts 
the most significant modes from a dataset (e.g., snapshots of 
velocity, pressure, or temperature fields), thereby isolating the 
principal coherent structures of the flow while discarding less 
energetic contributions. Its ability to identify dominant coherent 
structures has also been demonstrated in complex turbulent wake 
configurations involving flexible structures and fluid--structure 
interaction \cite{Salavatidezfouli2024WindTurbine}. Subsequently, 
Galerkin projection is employed to project the governing equations 
onto the reduced subspace, resulting in a model defined exclusively 
by the coefficients of the dominant modes 
\cite{Bergmann2009,Berkooz1993}. This procedure drastically lowers 
the number of unknowns compared with the full-order model. While 
POD--Galerkin methods were first introduced within the context of 
finite element (FE) discretizations, increasing interest has recently 
been directed toward their extension to finite volume methods (FVMs), 
owing to their robustness and their ability to guarantee both local 
and global conservation properties 
\cite{chakir:hal-01420726,2014,Lorenzi2016,Stabile2018,Carlberg2018,
Haasdonk2008,Eymard2000,2012,Syrakos2017,VersteegH.K2007AItc,
fletcher1991computational}.

In this work, simulations are carried out within the OpenFOAM framework. The solver relies on a finite volume (FV) discretization on collocated grids, which is a common choice in large-scale CFD applications. From a practical standpoint, collocated layouts simplify the implementation of the numerical algorithms and typically require less memory than staggered-grid approaches \cite{ferziger2012computational, Harlow1965,Vasilyev2000}.

Despite their widespread use, FV-based POD--Galerkin reduced-order models may still suffer from numerical instabilities. In incompressible flow settings, these instabilities typically originate from the interaction between convection-dominated transport and the pressure--velocity coupling introduced by the discretization of the Navier--Stokes equations~\cite{Patankar2018,Piller2004,Peri1988,Zang1994,Jia2025,Rezaei2025}. If not properly addressed, such effects can adversely impact both the accuracy and robustness of the reduced-order approximation.



Several studies on POD--Galerkin ROMs have shown that enforcing a divergence-free reduced velocity space leads to the elimination of the pressure gradient contribution in the reduced momentum formulation. Constructing a velocity-only reduced-order formulation with stable numerical behavior on collocated grids, however, is not straightforward because the discrete divergence and gradient operators do not satisfy the necessary compatibility conditions. As a result, stability usually requires the combined use of Rhie–Chow interpolation at the full-order level and additional pressure stabilization techniques at the reduced level \cite{Kalashnikova2011,Iollo2000,Balajewicz2013,Balajewicz2016,MA2002,Sanderse2020,Rhie1983}.

To obtain reliable velocity and pressure reconstructions, both in FE and FV reduced-order frameworks, two strategies are commonly adopted. One is the supremizer enrichment of the velocity basis, which enforces the inf–sup condition and thereby removes the pressure instabilities that otherwise arise. While effective, this method typically requires a number of supremizer modes equal to or greater than the number of pressure modes, which inevitably increases the dimensionality and online computational cost of the ROM. The second strategy leverages the velocity approximation within the ROM to solve a pressure Poisson equation. Although widely used, this approach raises difficulties in prescribing appropriate boundary conditions for the pressure equation \cite{Ballarin2014,Akhtar2009,NOACK2005,Caiazzo2014,Kean2020,gresho2000incompressible,Liu2010}.

Several studies have shown that stabilization alone does not necessarily guarantee high accuracy in reduced-order velocity and pressure fields. Compared to direct POD projection of full-order solutions, the resulting errors may increase by one to two orders of magnitude. This limitation has been documented even for canonical non-parametric laminar test cases, including the lid-driven cavity problem \cite{Stabile2018,Ballarin2014,Kean2020,Busto2020}.


An alternative pressure recovery approach involves using a pressure reduced basis in conjunction with the associated temporal coefficients originally derived for the velocity field to reconstruct the pressure contribution in the momentum equations. This method offers the computational advantage of only requiring the resolution of the momentum system. However, because the pressure field is reconstructed using the exact same set of coefficients as the velocity, this strategy necessitates an identical number of modes for both the velocity and pressure reduced bases, which can impose constraints on the overall efficiency of the model\cite{Bergmann2009,Lorenzi2016}.


More recent developments in the FE community include stabilization via local projection methods and the introduction of an artificial compressibility condition to replace the incompressibility constraint in the Navier–Stokes equations \cite{Caiazzo2014,Rubino2020,Novo2021,DeCaria2020}.

In this work, building on our earlier study \cite{https://doi.org/10.48550/arxiv.2010.06964}, we extend the discretize--then--project strategy to fully three-dimensional turbulent configurations.We consider a reduced-order modelling (ROM) of
the incompressible Navier--Stokes equations on collocated grids, for which no supremizer enrichment or additional pressure boundary-enforcement procedures are introduced at the reduced level beyond the standard pressure reference constraint (fixing $p = 0$ at a reference cell) inherited from the FOM solver. The approach is formulated directly at the discrete level: the operators of the full-order model (FOM) together with the boundary contributions are projected onto POD-based reduced spaces. This discretize--then--project strategy ensures consistency between the reduced velocity and pressure representations \cite{Kalashnikova2012,https://doi.org/10.48550/arxiv.2010.06964}.


However, applying the same procedure to the turbulence field does not yield
consistent results at the reduced level. To address this limitation, we
propose a hybrid approach in which the turbulent viscosity is incorporated
through a data-driven methodology
\cite{Kalashnikova2012,Sirisup2005,Stabile2017,Ullmann2014,Fick2018,
Carlberg2013,Lee2020,Placzek2011}. Related comparisons between intrusive
POD--Galerkin formulations and recurrent non-intrusive reduced models have
demonstrated the complementary advantages of physics-based projection and
data-driven temporal prediction in unsteady CFD applications
\cite{Halder2026Cerebrovascular}.


Data-driven methods rely on learning from available data rather than
prescribing an explicit physical model, often using machine learning
techniques to represent the system
\cite{hajisharifi2023non}. Broadly speaking, these methods can be divided into either purely data-driven or physics-informed data-driven approaches. In this work, we adopt a purely data-driven strategy to model the non-intrusive path of our hybrid framework, with the objective of evaluating the consistency of the reconstructed turbulence field. Specifically, we investigate and compare three different neural network architectures, Multilayer Perceptron (MLP), Transformers, and Long Short-Term Memory (LSTM) networks, to predict the temporal evolution of the eddy viscosity coefficients. This allows us to identify the most suitable model for capturing the complex, transient nonlinear interactions inherent in 3D turbulent flows\cite{Ivagnes2023}.



This work adopts a POD–Galerkin reduced order modeling framework following a "discretize-then-project" strategy. As detailed in the sources, the strategy approach is fundamentally intrusive, as it requires direct access to the fully discrete system matrices and numerical algorithms of the high-fidelity solver to project the discrete operators onto the reduced basis. The reduced-order model is constructed by directly projecting the fully discrete operators of the full-order model, together with the associated boundary contributions, onto the reduced spaces. This choice allows the reduced system to retain a close correspondence with the original FV formulation, including both linear and nonlinear terms arising from the discretization. The practical implementation of this approach is supported by ITHACA-FV, an open-source library built on top of OpenFOAM that provides direct access to the discrete operators and numerical solvers of the full-order model \cite{Placzek2011,Baiges2013}.

The present work builds on the study of Star et al.\cite{https://doi.org/10.48550/arxiv.2010.06964}
, but departs from it in several important aspects. In that work, the analysis was restricted to two-dimensional laminar configurations, namely the open and lid-driven cavity configurations, with the primary goal of isolating issues related to pressure–velocity coupling. Here, we extend the same framework to a fully three-dimensional setting, thereby increasing both the physical complexity of the problem and the demands placed on the reduced-order model. This transition notably increases the physical realism and the computational demands placed upon the reduced system. Furthermore, our work addresses convection-dominated turbulent flows, a regime explicitly identified in the sources as being beyond the scope of the original study. Instead of relying solely on intrusive projection for all terms, we implement a hybrid data-driven methodology, where the turbulent closure is resolved non-intrusively, thereby bridging the gap between physical consistency and machine learning performance for complex 3D flows.

In this framework, LSTM, Transformer and MLP were investigated and compared in order to identify the most suitable model for accurately reproducing the temporal dynamics and nonlinear interactions of the turbulent flow.

The paper is structured as follows. Section \ref{sec:methodologies} presents the methodology, starting from the mathematical formulation, detailing the FOM discretization, and describing the ROM approach. Section \ref{sec:NumericalResults} reports the numerical results, including the selected test case, the flow behavior obtained with the proposed methodology, and a comparison of different neural network architectures to assess which performs best. Finally, Section \ref{sec:Conclusions} summarizes the main conclusions and outlines possible directions for future research.

\section{Methodology} \label{sec:methodologies}
In this section, we will show the methodology used to obtain the results. 
\subsection{Mathematical model}
In this study, the incompressible Navier--Stokes equations are adopted to describe the flow dynamics.

Let $\Omega$ denote the computational domain corresponding to the fluid region. We consider an incompressible Newtonian fluid characterized by constant density $\rho$ and kinematic viscosity $\nu$, in the absence of gravitational and body force effects. Under these assumptions, the governing equations consist of the conservation laws for mass and momentum.


\begin{equation}\label{eq:continuity_equation}
    \nabla \cdot \mathbf{u}=0 \quad \text{in  } \Omega,
\end{equation} 
\begin{equation}\label{eq:momentum_equation}
    \frac {\partial \mathbf{u}}{\partial t}=-\nabla\cdot \mathbf{(u \otimes u)+\nu \nabla \cdot (\nabla \mathbf{u})-\nabla p} \quad \text{in  } \Omega.
\end{equation}
%

Here, $\mathbf{u}(\mathbf{x},t)$ denotes the velocity field defined over the spatial domain $\Omega \subset \mathbb{R}^d$, with $d=2$ or $3$, while $p(\mathbf{x},t)$ is the pressure field normalized by the constant fluid density $\rho$. The time variable is indicated by $t$. The momentum equation accounts for convective transport, viscous diffusion, and pressure forces through the corresponding differential operators.

To account for the effects of turbulence within the FOM, a Large Eddy Simulation (LES) approach is adopted. In this framework, the governing equations are filtered at the grid scale, and the influence of the unresolved subgrid-scale (SGS) motions is modeled using the Smagorinsky subgrid-scale model. 

Consequently, the molecular kinematic viscosity $\nu$ in \autoref{eq:momentum_equation} is replaced by an effective viscosity $\nu_{eff}$, defined as:

\begin{equation}\label{eq:effective_viscosity}
    \nu_{eff} = \nu + \nu_T,
\end{equation}
where $\nu_T$ represents the turbulent eddy viscosity. In the context of the Smagorinsky model, this eddy viscosity is computed based on the resolved velocity field:
\begin{equation}\label{eq:smagorinsky}
    \nu_T = (C_s \Delta)^2 |\bar{S}|,
\end{equation}
where $C_s$ is the Smagorinsky constant, which is set to $0.2$ in this work, $\Delta$ is the filter width related to the grid size, and $|\bar{S}| = \sqrt{2\bar{S}_{ij}\bar{S}_{ij}}$ is the magnitude of the resolved-scale strain rate tensor. Additionally, the turbulent Prandtl number $\text{Pr}$ is set to 0.9. This formulation ensures that the FOM snapshots used for training accurately capture the energy dissipation and transient dynamics of the 3D lid-driven cavity before the data is passed to the reduced-order framework. This explicit mathematical treatment at the FOM level ensures that the high-fidelity snapshots of velocity, pressure, and viscosity ($\mathbf{u}, p, \nu_T$) contain the necessary physical information for the subsequent training of the data-driven closure.

The governing equations are augmented by the initial velocity field:
\begin{equation*}
    \mathbf{u}(\mathbf{x},0)=\mathbf{u_0}(\mathbf{x})  \quad \text{in  } \Omega,
\end{equation*}
in which this condition is assumed to be divergence-free, $\nabla \cdot \mathbf{u_0} = 0$. 

To solve the problem above, we have to impose boundary conditions. We consider two different time-independent boundary conditions: wall, moving lid; so we have $\partial \Omega = \Gamma_{\text{wall}} \cup \Gamma_{\text{up}}$.

The no-slip condition is satisfied at the wall in case of viscous fluid, so we have to impose the following condition: 

\begin{equation*}
    \mathbf{u}=\mathbf{u}_{\text{wall}}(\mathbf{x}) \quad \text{on  } \quad \Gamma_{\text{wall}} \quad \text{for} \quad \text{t} \geq 0,
\end{equation*}
where $\mathbf{u}_{\text{wall}}$ denotes the prescribed wall velocity. For stationary walls, this quantity is set equal to zero.


A moving lid condition is imposed on the upper boundary according to

\begin{equation*}
    \mathbf{u}=\mathbf{u}_{\text{up}}(\mathbf{x}) \quad \text{on  } \quad \Gamma_{\text{up}} \quad \text{for} \quad \text{t} \geq 0,
\end{equation*}

where $\mathbf{u}_{\text{up}}$ specifies the imposed tangential velocity of the lid, assumed to be uniform and known a priori.

For the pressure field, the boundary condition that we give is the zero-gradient condition on every wall.
The turbulent viscosity has the condition of fixed value at every wall, it is zero.

\subsection{FOM numerical discretization}

This section describes the finite-volume discretization of the governing equations, Equations~\ref{eq:continuity_equation} and~\ref{eq:momentum_equation}, formulated on a collocated grid arrangement (see Figure~\ref{fig:placeholder}).

The level of detail retained throughout this section is necessary: every discrete 
operator introduced here, $\mathbf{M}$, $\mathbf{M}_p$, $\mathbf{D}_p$, 
$\mathbf{C}_p$, $\mathbf{G}_p$, $\mathbf{G}_f$, and $\mathbf{I}_{p \to f}$, is 
projected onto the POD basis in Section~\ref{sec:ROM} to form the reduced operators 
that govern the online ROM simulation.


\begin{figure}[H]
    \centering
    \resizebox{0.5\linewidth}{!}{%
        \begin{tikzpicture}[
            grid line/.style={thick, black},
            main node/.style={circle, fill=black, inner sep=0pt, minimum size=8pt},
            face node/.style={circle, fill=black, inner sep=0pt, minimum size=6pt},
            vec arrow/.style={->, >=latex, thick, black},
            dim arrow/.style={<->, >=latex, thick, red!80!black},
            font=\large
        ]
        
            \def\cw{4} 
            \def\ch{4} 
        
            \fill[cyan!20] (0,0) rectangle (\cw, \ch);
        
            \draw[grid line] (0, -1) -- (0, \ch + 1);
            \draw[grid line] (\cw, -1) -- (\cw, \ch + 1);
            \draw[grid line] (2*\cw, -1) -- (2*\cw, \ch + 1);
            \draw[grid line] (-1, 0) -- (2*\cw + 1, 0);
            \draw[grid line] (-1, \ch) -- (2*\cw + 1, \ch);
        
            \node[main node] (P) at (0.5*\cw, 0.5*\ch) {};
            \node[main node] (N) at (1.5*\cw, 0.5*\ch) {};
            \node[face node] (fw) at (0, 0.5*\ch) {};
            \node[face node] (fe) at (\cw, 0.5*\ch) {};
            \node[face node] (fn) at (0.5*\cw, \ch) {};
            \node[face node] (fs) at (0.5*\cw, 0) {};
            \node[face node] (ff) at (2*\cw, 0.5*\ch) {};
            \node[face node] at (1.5*\cw, \ch) {};
            \node[face node] at (1.5*\cw, 0) {};
        
            \draw[vec arrow] ($(fw) - (0.4, 0)$) -- ($(fw) + (0.4, 0)$);
            \node[above left] at (fw) {$\bm{u}_{f,w}$};
            
            \draw[vec arrow] ($(fe) - (0.4, 0)$) -- ($(fe) + (0.4, 0)$);
            \node[above right] at (fe) {$\bm{u}_{f,e}$};
            
            \draw[vec arrow] ($(fn) - (0, 0.4)$) -- ($(fn) + (0, 0.6)$);
            \node[above right] at (fn) {$\bm{u}_{f,n}$};
            
            \draw[vec arrow] ($(fs) - (0, 0.4)$) -- ($(fs) + (0, 0.6)$);
            \node[below right] at (fs) {$\bm{u}_{f,s}$};
        
            \node[align=center, yshift=0.8cm] at (P) {$\bm{u}_{p,P}$ \\ $\bm{p}_P$};
            \node[align=center, yshift=0.8cm] at (N) {$\bm{u}_{p,N}$ \\ $\bm{p}_N$};
        
            \draw[dim arrow] (0, -1.2) -- (\cw, -1.2) node[midway, below, text=red!80!black] {$\Delta x$};
            \draw[dim arrow] (-1.2, 0) -- (-1.2, \ch) node[midway, left, text=red!80!black] {$\Delta y$};
            \draw[dim arrow] (P) -- (N) node[midway, below, yshift=-2pt] {$\bm{d}$};
        
        \end{tikzpicture}%
    }

    
    \caption{Two-dimensional collocated control-volume layout. Cell-centered velocity and pressure are associated with the control volumes $P$ and $N$. Face fluxes $\mathbf{u}_f$ are indicated on the cell faces, and $\Delta x$, $\Delta y$ denote the grid spacing. The vector $\mathbf{d}$ connects adjacent cell centers.}

    \label{fig:placeholder}
\end{figure}

To maintain local conservation of mass and momentum within the Finite Volume framework, the governing equations (Eqs. \ref{eq:continuity_equation} and \ref{eq:momentum_equation}) are expressed in their integral form over each computational cell:


\begin{equation}\label{eq:int_mass}
    \int_{\partial(\Omega_h)_k} \mathbf{n} \cdot \mathbf{u}dS=0, 
\end{equation} 
\begin{equation}\label{eq:int_momentum}
    \int_{(\Omega_h)_k}\frac {\partial \mathbf{u}}{\partial t} d\Omega = -\int_{\partial(\Omega_h)_k}{(\mathbf{n}\cdot \mathbf{u})\mathbf{u}dS+\nu \int_{\partial(\Omega_h)_k} \mathbf{n}\cdot (\nabla \mathbf{u})dS-\int_{\partial(\Omega_h)_k} \mathbf{n}pdS},
\end{equation}
Here, $k$ denotes a generic control volume of the mesh, with $(\Omega_h)_k$ indicating the associated cell volume and $\partial(\Omega_h)_k$ its boundary surface. The integrals are evaluated over the cell interior and its boundary through the volume and surface measures $d\Omega$ and $dS$, respectively.

For spatial discretization, we rewrite the governing equations (Equations \ref{eq:int_mass} and \ref{eq:int_momentum}) with matrix-vector notation:

\begin{equation}\label{eq:disc_mass}
    \mathbf{M}\mathbf{u_f}=\mathbf{0}, 
\end{equation} 

\begin{equation}\label{eq:disc_momentum}
    \frac {d\mathbf{u}_p}{dt}=-\mathbf{C}_p(\mathbf{u}_f,\mathbf{u}_p)-\mathbf{r}_p^C+\nu \mathbf{D}_p\mathbf{u}_p-\mathbf{G}_p\mathbf{p}_p+\nu \mathbf{r}_p^D,
\end{equation}

Here, $\mathbf{p}_p=(p_{p,1},p_{p,2},...,p_{p,h})^T \in \mathbb{R}^h$ denotes the pressure degrees of freedom associated with the cell centers, while $\mathbf{u}_p \in \mathbb{R}^{dh}$ collects the corresponding cell-centered velocity components. In three spatial dimensions ($d=3$), the velocity vector is organized by grouping the Cartesian components, such that $\mathbf{u}_p = ((\mathbf{u}_p)_1, (\mathbf{u}_p)_2, (\mathbf{u}_p)_3)^T$, with each component containing the values defined over the $h$ control volumes. The face-based velocity vector $\mathbf{u}_f \in \mathbb{R}^{dm}$ is defined on the mesh faces and is obtained from the cell-centered velocity field through a linear interpolation operator $\mathbf{I}_{p\rightarrow f}$, together with prescribed boundary contributions, according to

\begin{equation*}
    \mathbf{u}_f = \mathbf{I}_{p\rightarrow f}\mathbf{u}_p + \mathbf{u}_b .
\end{equation*}
The vector $\mathbf{u}_b \in \mathbb{R}^{dm}$ contains the velocity values imposed on the boundary faces of the computational domain.


The discrete operators appearing in the formulation are defined as follows. The matrix $\mathbf{M} \in \mathbb{R}^{h\times dm}$ maps face-based quantities to cell centers and corresponds to the discrete divergence operator. Diffusive effects are accounted for through the cell-centered Laplacian operator $\mathbf{D}_p \in \mathbb{R}^{dh\times dh}$. Convective transport is described by the nonlinear operator $\mathbf{C}_p(\mathbf{u}_f,\mathbf{u}_p) \in \mathbb{R}^{dh\times dh}$, which depends on both face-centered and cell-centered velocity fields. The pressure gradient contribution is introduced via the discrete gradient operator $\mathbf{G}_p \in \mathbb{R}^{dh\times h}$. In addition, the vectors $\mathbf{r}_p^C \in \mathbb{R}^{dh}$ and $\mathbf{r}_p^D \in \mathbb{R}^{dh}$ collect the convective and diffusive source terms, respectively.



The spatial discretization of the governing equations, written in integral form in Equations~\ref{eq:int_mass} and~\ref{eq:int_momentum}, is described below for a generic control volume $k$. The resulting discrete operators appearing in Equations~\ref{eq:disc_mass} and~\ref{eq:disc_momentum} are detailed accordingly.

The continuity equation is discretized as following:
\begin{equation}\label{eq:int_disc_mass}
    \int_{\partial(\Omega_h)_k}\mathbf{n} \cdot \mathbf{u}dS = \bm{\sum}_{i=1}^{N_f} \int_{S_{f,i}} \mathbf{n} \cdot \mathbf{u}dS \approx \bm{\sum}_{i=1}^{N_f} {\mathbf{S}_{f,i}} \cdot {\mathbf{u}_{f,i}} = \bm{\sum}_{i=1}^{N_f} {\phi_{f,i}} = 0,     
\end{equation}

Here, $N_f$ denotes the number of faces associated with control volume $k$, and $S_f$ represents the area of a generic face. The discrete divergence operator $\mathbf{M}$ introduced in Equation~\ref{eq:disc_mass} is constructed from the outward-oriented face area vectors corresponding to the faces of the computational mesh. As expressed in Equation~\ref{eq:int_disc_mass}, the incompressibility condition is enforced at the level of the face fluxes, defined as $\bm{\phi}_f = \mathbf{S}_f \cdot \mathbf{u}_f$.

For subsequent developments, it is also convenient to introduce a center-to-center divergence operator $\mathbf{M}_p \in \mathbb{R}^{h \times dh}$, obtained by composing the face-based operator with the interpolation from cell centers to faces, namely
\begin{equation*}
    \mathbf{M}_p = \mathbf{M}\mathbf{I}_{p\rightarrow f}.
\end{equation*}


Now we can write the momentum equation (Equation \ref{eq:disc_momentum}) as following:
\begin{equation}\label{eq:rewrite_disc_momentum}
    \mathbf{M}\mathbf{u}_f=\mathbf{M}\mathbf{I}_{p\rightarrow{f}}\mathbf{u}_p+\mathbf{M}\mathbf{u}_b=\mathbf{M}_p\mathbf{u}_p+\mathbf{r}_p^M=\mathbf{0},
\end{equation}

The vector $\mathbf{r}_p^M \in \mathbb{R}^h$ accounts for the boundary contributions associated with the continuity equation. It is defined through the action of the discrete divergence operator on the prescribed boundary velocity field,
\begin{equation}\label{eq:pressure_matrix}
    \mathbf{r}_p^M=\mathbf{M}\mathbf{u}_b,  
\end{equation}


The pressure gradient contribution is obtained by applying the FV formulation to the pressure field over the boundary of a generic control volume. In particular, the surface integral is evaluated by summing the fluxes across the cell faces,
\begin{equation*}
    \int_{\partial(\Omega_h)_k} \mathbf{n} p \, dS
    = \sum_{i=1}^{N_f} \int_{S_{f,i}} \mathbf{n} p \, dS
    \approx \sum_{i=1}^{N_f} \mathbf{S}_{f,i} \, p_{f,i} .
\end{equation*}

Here, $\mathbf{p}_f$ and $\mathbf{p}_p$ denote the pressure values at face centers and cell centers, respectively. The two representations are related through a linear interpolation operator $\mathbf{\Pi}_{p\rightarrow f} \in \mathbb{R}^{m \times h}$, such that
\begin{equation*}
    \mathbf{p}_f = \mathbf{\Pi}_{p\rightarrow f}\mathbf{p}_p .
\end{equation*}
Within this formulation, the discrete pressure gradient operator $\mathbf{G}_p$ is constructed from the face area vectors weighted by the interpolation coefficients contained in $\mathbf{\Pi}_{p\rightarrow f}$.



The diffusive term is discretized as following:
\begin{equation*}
    \int_{\partial(\Omega_h)_k} \mathbf{n} \cdot (\nabla \mathbf{u}) dS = \bm{\sum}_{i=1}^{N_f} \int_{S_{f,i}} \mathbf{n} \cdot \nabla \mathbf{u} dS \approx \bm{\sum}_{i=1}^{N_f} |{\mathbf{S}_{f,i}}| \frac{\mathbf{u}_{p,N} - \mathbf{u}_{p,P}}{|\mathbf{d}|},
\end{equation*}


where $\mathbf{d}$ denotes the vector connecting the centers of the neighboring cells $N$ and $P$. $\mathbf{D}_p$ is the discrete operator of the discrete term, given by the areas of the faces and the distances between the two centers. In case one cell is near to another of the same domain, the discretization changes into:
\begin{equation*}
    |{\mathbf{S}_{f,b}}| \frac{\mathbf{u}_{f,b}-\mathbf{u}_{p,P}}{|\mathbf{d}_n|},
\end{equation*}

which can be decomposed into two contributions:

\begin{equation*}
    |{\mathbf{S}_{f,b}}| \frac{\mathbf{0}-\mathbf{u}_{p,P}}{|\mathbf{d}_n|}+|{\mathbf{S}_{f,b}}| \frac{\mathbf{u}_{f,b}-\mathbf{0}}{|\mathbf{d}_n|},
\end{equation*}


where the diffusion contribution is accounted for through the operator $\mathbf{D}_p$, while boundary-related effects are gathered in the vector $\mathbf{r}_p^D$. Here, $\mathbf{u}_{f,b}$ denotes the velocity prescribed on boundary face $b$, $\mathbf{d}_n$ represents the vector connecting the boundary face to the corresponding cell center, and $\mathbf{S}_{f,b}$ is the associated face area vector.


Finally, the convective term is discretized as following:
\begin{equation*}
    \int_{\partial(\Omega_h)_k} (\mathbf{n} \cdot \mathbf{u}) \mathbf{u} dS = \bm{\sum}_{i=1}^{N_f} \int_{S_{f,i}} (\mathbf{n} \cdot \mathbf{u}) \mathbf{u} dS \approx \bm{\sum}_{i=1}^{N_f} ({\mathbf{S}_{f,i}} \cdot {\mathbf{u}_{f,i}}) = \bm{\sum}_{i=1}^{N_f} \bm{\phi}_{f,i} \mathbf{u}_{f,i}.
\end{equation*}

The convective contribution is represented by the operator $\mathbf{C}_p(\mathbf{u}_f,\mathbf{u}_p)$, whose nonlinear dependence arises from the face fluxes $\bm{\phi}_f$. When a control-volume face lies on the boundary of the computational domain, the contribution $({\mathbf{S}_{f,b}} \cdot {\mathbf{u}_{f,b}})\mathbf{u}_{f,b}$ is incorporated into the boundary vector $\mathbf{r}_p^C \in \mathbb{R}^{dh}$.


In the present formulation, the discretized convection term exhibits a quadratic dependence on the velocity field. This follows from the fact that the face-centered velocity $\mathbf{u}_f$ is obtained through linear interpolation of the cell-centered velocity $\mathbf{u}_p$. Under this assumption, the convective contribution can be expressed in operator form as the action of a velocity-dependent matrix on $\mathbf{u}_p$, namely

%
\begin{equation}\label{eq:convective_operator}
    \tilde{\mathbf{C}}_p(\mathbf{u}_f)\mathbf{u}_p= \mathbf{C}_p(\mathbf{u}_f,\mathbf{u}_p).
\end{equation}

Applying  Equation \ref{eq:rewrite_disc_momentum} into Equation \ref{eq:disc_mass} and Equation \ref{eq:convective_operator} into Equation \ref{eq:pressure_matrix}, we obtain the complete spatioally discretize system of equations:

\begin{equation*} \label{eq:discrete_mass_conservation}
    \mathbf{M}_p\mathbf{u}_p+\mathbf{r}_p^M=\mathbf{0},
\end{equation*}
\begin{equation*} \label{eq:discrete_momentum_conservation}
    \frac{d\mathbf{u}_p}{dt}=-\tilde{\mathbf{C}}_p(\mathbf{u}_f)\mathbf{u}_p+\nu \mathbf{D}_p\mathbf{u}_p-\mathbf{G}_p\mathbf{p}_p+\mathbf{r}_p,
\end{equation*}

Here, $\mathbf{r}_p \in \mathbb{R}^{dh}$ denotes the combined boundary contribution, defined as $\mathbf{r}_p = -\mathbf{r}_p^C + \nu \mathbf{r}_p^D$. All discrete operators are appropriately scaled according to the finite-volume measures.


In the framework of explicit projection methods, a well-documented difficulty is the emergence of spurious pressure oscillations, commonly referred to as the checkerboard problem. This phenomenon typically arises when collocated grids are employed in conjunction with central-difference discretizations of both the velocity divergence $\nabla \cdot U$ and the pressure gradient $\nabla p$.

The underlying cause of this instability is the decoupling between pressure and velocity at adjacent cell centers. Such a decoupling leads to a violation of the compatibility conditions between the discrete divergence and gradient operators, which in turn gives rise to an extended stencil in the Pressure Poisson Equation (PPE) \cite{Stabile2018,Klaij2015,Date1993,Issa1986,Komen2020,hirsch2007numerical}.

A common approach used in practice is Rhie–Chow interpolation. The idea is to correct the face-centered velocities so that they are compatible with the pressure gradient, which avoids the appearance of non-physical oscillations in the pressure field. This correction is part of the PISO (Pressure Implicit with Splitting of Operators) algorithm and is implemented in solvers such as OpenFOAM\cite{Komen2020, Issa1986, Rhie1983}. In contrast, for explicit projection methods on collocated grids combined with explicit time integration schemes, the application of Rhie–Chow interpolation may be circumvented. Two principal strategies are typically employed to treat the fluxes in this context: the inconsistent flux method and the consistent flux method\cite{https://doi.org/10.48550/arxiv.2010.06964}.

Within the inconsistent flux method, the face velocity $\mathbf{u}_f$ is reconstructed through linear interpolation from the cell-centered velocity $\mathbf{u}_p$ together with the pressure field. Conversely, the consistent flux method introduces an auxiliary governing equation for the face-centered velocity $\mathbf{u}_f$, leading to a more robust and accurate formulation. Both approaches are generally implemented with explicit time integration schemes, such as the forward Euler method \cite{KazemiKamyab2015, Vuorinen2014,canuto2007spectral,Chorin1968,Sanderse2020,https://doi.org/10.48550/arxiv.2010.06964}.

The present study will concentrate on the consistent flux method.

\subsubsection{Inconsistent and Consistent flux methods} \label{CFM}

As explained in \cite{https://doi.org/10.48550/arxiv.2010.06964}, the projection method first evaluates an intermediate velocity $\mathbf{u}^*_p$ field by neglecting the pressure gradient contribution in the momentum equations:


\begin{equation*}
    \frac{\mathbf{u}^*_p-\mathbf{u}^n_p}{\Delta t}=-\tilde{\mathbf{C}}_p(\mathbf{u}^n_f)\mathbf{u}^n_p+\nu \mathbf{D}_p\mathbf{u}^n_p+\mathbf{r}_p.
\end{equation*}

The intermediate velocity field $\mathbf{u}_p^*$ is subsequently adjusted through the pressure contribution:


\begin{equation*}
    \mathbf{u}_p^{n+1}=\mathbf{u}_p^*-\Delta t\mathbf{G}_p\mathbf{p}_p^{n+1}.
\end{equation*}

To enforce the incompressibility constraint, the velocity field is expressed as:

\begin{equation*}
    \mathbf{M}_p\mathbf{u}_p^{n+1}+\mathbf{r}_p^M=(\mathbf{M}_p\mathbf{u}_p^*+\mathbf{r}_p^M)-\Delta t\mathbf{M}_p\mathbf{G}_p\mathbf{p}_p^{n+1} = \mathbf{0},
\end{equation*}
and considering:
\begin{equation}\label{eq:L_p}
    \mathbf{L}_p\mathbf{p}_p^{n+1}=\frac{1}{\Delta t}(\mathbf{M}_p\mathbf{u}_p^*+\mathbf{r}_p^M),
\end{equation}

here, $\mathbf{L}_p = \mathbf{M}_p \mathbf{G}_p \in \mathbb{R}^{h \times h}$ denotes the discrete Laplacian operator obtained by mapping cell-centered pressure gradients to the cell faces. In this formulation, the resulting pressure field may lose coupling between neighboring control volumes, leading to the appearance of non-physical checkerboard modes. To address this issue, the pressure gradient is instead evaluated directly at the cell faces by introducing an alternative operator $\mathbf{L}_f \in \mathbb{R}^{h \times h}$ in place of $\mathbf{L}_p$ \cite{Morinishi1998,Felten2006}.


We obtain:
\begin{equation}\label{eq:L_f}
    \mathbf{L}_f\mathbf{p}_p^{n+1}=\frac{1}{\Delta t}(\mathbf{M}_p\mathbf{u}_p^n+\mathbf{r}_p^M)+\mathbf{M}_p(-\tilde{\mathbf{C}}_p(\mathbf{u}_f^n)\mathbf{u}_p^n+\nu \mathbf{D}_p\mathbf{u}_p^n+\mathbf{r}_p),
\end{equation}
\begin{equation}\label{eq:u_p_n_1}
   \mathbf{u}_p^{n+1}=\mathbf{u}_p^n+\Delta t(-\tilde{\mathbf{C}}_p(\mathbf{u}_f^n)\mathbf{u}_p^n+\nu \mathbf{D}_p\mathbf{u}_p^n+\mathbf{r}_p)-\Delta t\mathbf{G}_p\mathbf{p}_p^{n+1}.
\end{equation}


In the inconsistent flux formulation, face-centered velocities $\mathbf{u}_f$ are obtained by interpolating the cell-centered field through the operator $\mathbf{I}_{p\rightarrow f}$. In contrast, the consistent flux approach incorporates the pressure field resulting from the solution of the pressure Poisson equation (Equation~\ref{eq:L_p}) in order to adjust the face fluxes accordingly.


Firstly we apply the linear interpolation operator onto the Equation \ref{eq:L_f}:
\begin{equation}\label{eq:u_f_G_p}
    \mathbf{u}_f^{n+1}=\mathbf{I}_{p\rightarrow{f}}[\mathbf{u}_p^n+\Delta t(-\tilde{\mathbf{C}}_p(\mathbf{u}_f^n)\mathbf{u}_p^n+\nu \mathbf{D}_p\mathbf{u}_p^n+\mathbf{r}_p)-\Delta t\mathbf{G}_p\mathbf{p}_p^{n+1}].
\end{equation}

Rather than computing the pressure gradient at cell centers, it is evaluated directly on the cell faces using a discrete operator $\mathbf{G}_f \in \mathbb{R}^{dm \times h}$, yielding


\begin{equation}\label{eq:u_f}
    \mathbf{u}_f^{n+1}=\mathbf{I}_{p\rightarrow{f}}[\mathbf{u}_p^n+\Delta t(-\tilde{\mathbf{C}}_p(\mathbf{u}_f^n)\mathbf{u}_p^n+\nu \mathbf{D}_p\mathbf{u}_p^n+\mathbf{r}_p)]-\Delta t\mathbf{G}_f\mathbf{p}_p^{n+1},
\end{equation}

for a generic cell $k$, the pressure gradient contribution appearing in the last term of Equation \ref{eq:u_f} is evaluated as

\begin{equation*}
    \sum_{i=1}^{N_f}{\mathbf{S}_{f,i}} \frac{p_{p,N}-p_{p,P}}{|\mathbf{d}|}.
\end{equation*}

The operator $\mathbf{G}_f$ is constructed from coefficients involving the face-normal vectors and inverse center-to-center distances. In this formulation, pressure gradients are computed directly from cell-centered pressure values, whereas the operator $\mathbf{G}_p$ relies on interpolated pressure fields evaluated at the cell centers.


By taking the divergence of Equation \ref{eq:u_f_G_p}, according to Equation \ref{eq:rewrite_disc_momentum}, we obtain:
\begin{equation}\label{eq:M_u_f}
    \mathbf{M}\mathbf{u}_f^{n+1}=(\mathbf{M}_p\mathbf{u}_p^n+\mathbf{r}_p^M)+\mathbf{M}_p[\Delta t(-\tilde{\mathbf{C}}_p(\mathbf{u}_f^n)\mathbf{u}_p^n+\nu \mathbf{D}_p\mathbf{u}_p^n+\mathbf{r}_p)]-\Delta t\mathbf{M}\mathbf{G}_f\mathbf{p}_p^{n+1}.
\end{equation}
So, substituting the pressure, we can see that the $\mathbf{u}_f$ is divergence free.

Finally, the system of equations for the CFM is formed by Equations \ref{eq:L_f}, \ref{eq:u_p_n_1} and \ref{eq:u_f}, which are solved in this particular order to obtain $\mathbf{u}_f, \mathbf{u}_p $ and $\mathbf{p}_p$ at time $t^{n+1}$.


\subsection{ROM approach}\label{sec:ROM}

The ROM framework adopted in this work follows the classical \textit{offline--online} decomposition, a foundational strategy in projection-based reduced-order modelling. The two phases have fundamentally different computational characters and are executed independently.

The \textbf{offline phase} is performed once, prior to any reduced-order simulation. It consists of three main steps: (i)~running the FOM (OpenFOAM via ITHACA-FV) to collect high-fidelity snapshots of the velocity $\mathbf{u}$, pressure $p$, and turbulent viscosity $\nu_T$ fields; (ii)~computing a set of global spatial basis functions via Proper Orthogonal Decomposition (POD) and projecting the discrete FOM operators onto these bases to obtain a compact set of precomputed algebraic tensors (stored in memory for later reuse); (iii)~extracting the temporal POD coefficients from the snapshot database and using them to train the data-driven closure model (LSTM, MLP, or Transformer). Although the offline phase is computationally intensive, it is executed only once and its cost is amortized over any number of subsequent online evaluations.

The \textbf{online phase} exploits all the precomputed reduced operators and the trained neural network to advance the reduced state in time. At each time step, the ROM assembles and solves a low-dimensional algebraic system using only matrix--vector operations on the precomputed tensors, while the turbulent viscosity coefficients are supplied by the trained network through a single forward pass. Because all expensive FOM-level computations are confined to the offline stage, the online simulation is orders of magnitude cheaper than the corresponding FOM solve.

ROM in this study is constructed through a hybrid strategy designed to handle the complexities of turbulent flows. As established in previous studies \cite{https://doi.org/10.48550/arxiv.2010.06964}, the standard 'discretize-then-project' approach provides high accuracy for laminar regimes but encounters a significant physical inconsistency when applied to turbulent viscosity field $\nu_T$.

Specifically, while the intrusive Galerkin projection accurately resolves the temporal coefficients for velocity $\mathbf{u}$ and pressure $p$, it fails to yield a physically coherent representation of the turbulent viscosity in the reduced space. To address this, we adopt a dual-path framework: 
\begin{itemize}
    \item \textbf{Intrusive Path:} The fields of pressure and velocity are resolved via the consistent flux discretize-then-project methodology to ensure mass conservation.
    \item \textbf{Non-Intrusive Path:} The turbulent viscosity field, originally modeled in the FOM using the Smagorinsky LES framework, is reconstructed using a data-driven closure (specifically an LSTM network).
\end{itemize}
This hybrid formulation enables the efficient simulation of unsteady turbulent flows by combining the mathematical rigor of projection-based ROMs with the predictive adaptability of deep learning.

\subsubsection{POD-GALERKIN Method}

The consistent flux formulation is projected onto a reduced space by applying the POD--Galerkin procedure to the fully discrete equations given in Equations~\ref{eq:L_f}, \ref{eq:u_p_n_1} and \ref{eq:u_f}. This formulation relies on the assumption that the full-order solution admits a reduced representation in terms of orthonormal spatial basis functions with coefficients that vary in time.



These modes, $\mathbf{\Phi}$ (cell-centered velocity), $\mathbf{X}$ (cell-centered pressure), $\mathbf{\Psi}$ (face velocity),  and $\mathbf{\Xi}$ (turbulent viscosity) are defined as global because they are computed across the entire discrete spatial domain $\Omega_h$ or the total set of face areas $\Sigma$. Unlike local nodal values, these basis functions capture the dominant energy-containing coherent structures of the flow during high-fidelity simulations. The  $\mathbf{u}_p$ as discrete cell-centered velocity fields can be approximated by:

\begin{equation}\label{eq:u_p}
\mathbf{u}_p \approx \mathbf{u}_{p,r}=\mathbf{\Phi} \mathbf{a},
\end{equation}
where $\mathbf{\Phi} = (\bm{\phi}_1,...,\bm{\phi}_{N_r^u}) \in \mathbb{R}^{dh \times {N_r^u}}$  contains the global velocity modes $\bm{\phi_i} \in \mathbb{R}^{dh}$ and $\mathbf{a}=(a^1,...,a^{N_r^u})\in \mathbb{R}^{N_r^u}$ are the temporal coefficients with $N_r^u$ the number of velocity modes. Similar approximations are used for the pressure field $\mathbf{p}_p$ and the face velocity field $\mathbf{u}_f$.

The pressure field in discrete form is represented as

\begin{equation}\label{eq:p_p}
    \mathbf{p}_p \approx \mathbf{p}_{p,r}=\mathbf{X}\mathbf{b},
\end{equation}

here, $\mathbf{X} = (\bm{\chi}_1,\ldots,\bm{\chi}_{N_r^p}) \in \mathbb{R}^{h \times N_r^p}$ collects the pressure basis functions defined at the cell centers, with each mode $\bm{\chi} = (\chi_1,\ldots,\chi_h)^T \in \mathbb{R}^h$. The vector $\mathbf{b}^n=(b^1,...,b^{N_r^p})\in \mathbb{R}^{N_r^p}$ are the corresponding time-dependent coefficients with $N_r^p$ the number of pressure modes.


Within the reduced-order framework, the face velocity field is expressed as

\begin{equation}\label{eq:u_f_approx}
    \mathbf{u}_f \approx \mathbf{u}_{f,r}=\mathbf{\Psi}\mathbf{c},
\end{equation}

where $\mathbf{\Psi}= (\mathbf{\psi}_1,...,\psi_{N_r^u}) \in \mathbb{R}^{dm \times N_r^u}$ is a matrix containing the face velocity modes $\mathbf{\psi_i} \in \mathbb{R}^m$ and $\mathbf{c}(t)=(c^1,...,c^{N_r^u}) \in \mathbb{R}^{N_r^u}$ are the the associated time-varying coefficients.

The turbulent viscosity field is represented as

\begin{equation}
\label{eq:nut_pod}
\nu_T \approx \nu_{T,r} = \mathbf{\Xi}\mathbf{d},
\end{equation}

where
$\mathbf{\Xi}=(\boldsymbol{\xi}_1,\ldots,
\boldsymbol{\xi}_{N_r^\nu})
\in \mathbb{R}^{h\times N_r^\nu}$
collects the turbulent-viscosity POD modes,
and
$\mathbf{d}(t)=
(d^1,\ldots,d^{N_r^\nu})
\in \mathbb{R}^{N_r^\nu}$
contains the corresponding time-dependent modal coefficients.

The computation of these global modes follows the Method of Snapshots. Taking the velocity field as a representative example, the procedure involves:

\begin{enumerate}
    \item \textbf{Snapshots Collection:} A set of $N_s^u$ high-fidelity solutions is collected: $\mathcal{S} = \{\mathbf{u}_p^1, \dots, \mathbf{u}_p^{N_s^u}\}$.
    \item \textbf{Correlation Matrix Construction:} A matrix $\mathbf{C}^u \in \mathbb{R}^{N_s^u \times N_s^u}$ is computed using the discrete $L_2$-inner product over the entire domain:
    \begin{equation*}
        C_{i,j}^u = (\mathbf{u}_p^i, \mathbf{u}_p^j)_{L_2(\Omega_h)} = \sum_{k=1}^h \mathbf{u}_{p,k}^i \cdot \mathbf{u}_{p,k}^j (\Omega_h)_k,
    \end{equation*}
    where $(\Omega_h)_k$ is the volume of the $k$-th mesh cell.
    \item \textbf{Eigenvalue Problem Resolution:} We solve $\mathbf{C}^u \mathbf{Q}^u = \mathbf{Q}^u \mathbf{\lambda}^u$ to determine the eigenvectors $\mathbf{Q}^u$ and eigenvalues $\lambda^u$.
    \item \textbf{Global Mode Reconstruction:} The modes $\bm{\phi}_i$ are constructed as an orthonormal linear combination of the snapshots:
    \begin{equation*}
        \bm{\phi}_i = \frac{1}{N_s^u \sqrt{\lambda_{ii}^u}} \sum_{n=1}^{N_s^u} \mathbf{u}_p^n Q_{in}^u \quad \text{for} \quad i = 1, \dots, N_r^u.
    \end{equation*}
\end{enumerate}

The cell-centered velocity modes $\bm{\phi}$ are only approximately divergence-free, matching the behaviour of the cell-centered velocity field $\mathbf{u}_p$. Consequently, the pressure contribution must be retained in the reduced-order formulation.

We now derive the Galerkin projection corresponding to the consistent flux method. To this end, the approximations of the cell-centered velocity $\mathbf{u}_p$, the cell-centered pressure $\mathbf{p}_p$, and the face velocity $\mathbf{u}_f$ (Equations 
\ref{eq:u_p}, \ref{eq:p_p}, and \ref{eq:u_f_approx}) are substituted into Equations \ref{eq:L_f}, \ref{eq:u_p_n_1} and \ref{eq:u_f}. The Pressure Poisson Equation (PPE) is subsequently formulated in the reduced pressure space by left-multiplying both sides of the equation by $\mathbf{X}^T \mathbf{V}$:


\begin{equation*}
    \mathbf{X}^T \mathbf{V} \mathbf{L}_f \mathbf{X} \mathbf{b}^{n+1} = \frac{1}{\Delta t} (\mathbf{X}^T \mathbf{V} \mathbf{M}_p \bm{\Phi} \mathbf{a}^n + \mathbf{X}^T \mathbf{V} \mathbf{r}_p^M) + \mathbf{X}^T \mathbf{V} \mathbf{M}_p(-\tilde{\mathbf{C}}_p (\mathbf{I}_{p\rightarrow{f}} \bm{\Phi} \mathbf{a}^n) \bm{\Phi} \mathbf{a}^n + \nu \mathbf{D}_p \bm{\Phi} \mathbf{a}^n + \mathbf{r}_p),
\end{equation*}
where $\mathbf{V} \in \mathbb{R}^{h \times h}$ is a diagonal matrix collecting the cell-centered control-volume measures.


We can rewrite this equation as follows:
\begin{equation*}
    \mathbf{L}_r \mathbf{b}^{n+1} = \frac{1}{\Delta t} (\mathbf{M}_r \mathbf{a}^n + \mathbf{q}_r^M) -{\hat{\mathbf{A}}_r} (\mathbf{a}^n) \mathbf{a}^n + \nu \mathbf{B}_r \mathbf{a}^n + \mathbf{q}_r,
\end{equation*}
where
\begin{itemize}
    \item $\mathbf{L}_r = \mathbf{X}^T \mathbf{V} \mathbf{L}_f \mathbf{X} \in \mathbb{R}^{N_r^p \times N_r^p}$,
    \item $\mathbf{M}_r = \mathbf{X}^T \mathbf{V} \mathbf{M}_p\mathbf{\Phi} \in \mathbb{R}^{N_r^p \times N_r^u}$,
    \item $\mathbf{B}_r = \mathbf{X}^T \mathbf{V} \mathbf{M}_p \mathbf{D}_p \mathbf{\Phi} \in \mathbb{R}^{N_r^p \times N_r^u}$,
    \item $\mathbf{q}_r = \mathbf{X}^T \mathbf{V} \mathbf{M}_p \mathbf{r}_p \in \mathbb{R}^{N_r^p}$,
    \item $\hat{\mathbf{A}}_r= \mathbf{X}^T \mathbf{V} \mathbf{M}_p \tilde{\mathbf{C}}_p (\mathbf{I}_{p\rightarrow{f}} \mathbf{\Phi}_i) \mathbf{\Phi} \in \mathbb{R}^{N_r^p \times N_r^u N_r^p \times N_r^u}$.
\end{itemize}
The first three are linear terms, the fourth is a reduced vector and the fifth term is the non-linear convection term; they are all determined during the offline stage. Therefore, $\hat{\mathbf{A}}_{r} (\mathbf{a}^n) \mathbf{a}^n$ during the online stage is evaluated as:
\begin{equation*}
    \bm{\sum}_{i=1}^{N_r^u}(\mathbf{a}^n)^T \hat{\mathbf{A}}_{r,i} \mathbf{a}^n.
\end{equation*}
The discretized momentum Equation \ref{eq:u_p_n_1} is projected onto the reduced base obtained from the velocity modes, multiplying the left hand of the equation for $\mathbf{\Phi}^T \mathbf{V}$:

\begin{equation*}
    \mathbf{\Phi}^T \mathbf{V} \mathbf{\Phi} \mathbf{a}^{n+1} = \mathbf{\Phi}^T \mathbf{V} \mathbf{\Phi} \mathbf{a}^n + \Delta t \bm{\Phi}^T \mathbf{V} (-\tilde{\mathbf{C}}_p (\mathbf{I}_{p\rightarrow{f}} \mathbf{\Phi} \mathbf{a}^n) \mathbf{\Phi} \mathbf{a}^n + \nu \mathbf{D}_p \mathbf{\Phi} \mathbf{a}^n + \mathbf{r}_p) - \Delta t \mathbf{\Phi}^T \mathbf{V} \mathbf{G}_p \mathbf{X} \mathbf{b}^{n+1}.
\end{equation*}
We can rewrite it as follows:
\begin{equation*}
    \mathbf{a}^{n+1} = \mathbf{a}^n + \Delta t (-\hat{\mathbf{C}}_r (\mathbf{a}^n) \mathbf{a}^n + \nu \mathbf{D}_r \mathbf{a}^n + \mathbf{r}_r) - \Delta t \hat{\mathbf{G}_r} \mathbf{b}^{n+1},
\end{equation*}
with 
\begin{itemize}
    \item $\mathbf{D}_r = \mathbf{\Phi}^T \mathbf{V} \mathbf{D}_p \mathbf{\Phi} \in \mathbb{R}^{N_r^u \times N_r^u}$,
    \item $\hat{\mathbf{G}}_r = \mathbf{\Phi}^T \mathbf{V} \mathbf{G}_p \mathbf{X} \in \mathbb{R}^{N_r^u \times N_r^p}$,
    \item $\mathbf{r}_r = \mathbf{\Phi}^T \mathbf{V} \mathbf{r}_p \mathbf{\Phi} \in \mathbb{R}^{N_r^u}$,
    \item $\hat{\mathbf{C}}_r (\mathbf{a}^n) \in \mathbb{R}^{N_r^u \times N_r^u \times N_r^u}$.
\end{itemize}
The first two are reduced matrices, the third is a reduced vector and the fourth is a tensor; they are all determined during the offline stage.

To derive the governing system associated with the consistent flux formulation, the reduced representations of the cell-centered velocity, face-centered velocity, and pressure fields Equations \ref{eq:u_p}, \ref{eq:p_p}, and \ref{eq:u_f_approx} are inserted into the corresponding full-order discrete Equations \ref{eq:L_f}, \ref{eq:u_p_n_1}, and \ref{eq:u_f}, yielding:


%
%

\begin{equation}\label{eq:rom_p1}
    \mathbf{L}_r \mathbf{b}^{n+1} = \frac{1}{\Delta t} (\mathbf{M}_r \mathbf{a}^n + \mathbf{q}_r^M) - \mathbf{A}_r (\mathbf{c}^n) \mathbf{a}^n + \nu \mathbf{B}_r \mathbf{a}^n + \mathbf{q}_r,
\end{equation}
Equation~\eqref{eq:rom_p1} is the reduced pressure Poisson equation:
it produces the new pressure coefficients $\mathbf{b}^{n+1}$ from the
current velocity state $\mathbf{a}^n$ and face-flux coefficients
$\mathbf{c}^n$. The term $\mathbf{M}_r\,\mathbf{a}^n$ is the reduced
divergence of the velocity field; it acts as the source term that drives
the pressure correction enforcing the incompressibility constraint at the
next step.

In incompressible flow, pressure is determined only up to an arbitrary 
constant since only pressure gradients appear in the momentum equation. 
The FOM solver fixes this gauge ambiguity by imposing $p = 0$ at a 
reference cell. At the reduced level, this same constraint must 
be enforced explicitly.

Following the residual form of Equation~\eqref{eq:rom_p1}, the pressure reference-cell condition is enforced by
adding the term $\tau \mathbf{K}\mathbf{b}^{n+1}$ to the residual of the reduced pressure Poisson equation:
\begin{equation}
\mathbf{L}_r\mathbf{b}^{n+1} + \tau \mathbf{K}\mathbf{b}^{n+1} = \mathbf{f}^n,
\label{eq:reference_pressure_constraint}
\end{equation}
or equivalently,
\[
(\tau \mathbf{K} + \mathbf{L}_r)\mathbf{b}^{n+1} = \mathbf{f}^n .
\]
Here, $\mathbf{f}^n$ denotes the right-hand side of Equation~\eqref{eq:rom_p1}, 
$\mathbf{K} = \mathrm{diag}(\mathbf{v})$ is a diagonal matrix whose entries
$v_i$ are the values of the $i$-th pressure POD mode evaluated at the reference cell
(mesh point index 0), and $\tau \geq 0$ is a scalar parameter controlling the strength
of the reduced pressure reference-cell constraint. This term enforces at the reduced level
the same pressure gauge-fixing condition imposed in the FOM solver, namely $p=0$ at the
reference cell. When $\tau = 0$, the reference-cell constraint is not enforced through this
additional residual term. The value of $\tau$ is selected through a direct reconstruction-error
sweep as described in Section~\ref{sec:Re3200}.

\begin{equation} \label{eq:rom_u1}
    \mathbf{a}^{n+1} = \mathbf{a}^n + \Delta t (-\mathbf{C}_r (\mathbf{c}^n) \mathbf{a}^n + \nu \mathbf{d}_r \mathbf{a}^n + \mathbf{r}_r) - \Delta t \hat{\mathbf{G}_r} \mathbf{b}^{n+1},
\end{equation}
Equation~\eqref{eq:rom_u1} advances the cell-centred velocity
coefficients $\mathbf{a}^{n+1}$. The last term $-\Delta t \hat{\mathbf{G}_r} \mathbf{b}^{n+1}$ is
the pressure gradient correction: it uses the pressure coefficients obtained from the reduced pressure
solve, Equation~\eqref{eq:rom_p1} together with the reference-cell constraint in
Equation~\eqref{eq:reference_pressure_constraint}, to push $\mathbf{a}^{n+1}$
toward the divergence-free constraint. The convective term uses $\mathbf{c}^n$, 
the face-flux coefficients from the previous step, rather than an
interpolation of $\mathbf{a}^n$ directly; this is the defining feature
of the consistent flux approach, which tracks face fluxes as independent
reduced unknowns rather than reconstructing them from cell-centred
quantities.

\begin{equation} \label{eq:rom_uf1}
    \mathbf{W}_r \mathbf{c}^{n+1} = \mathbf{N}_r \mathbf{a}^n + \Delta t (-\mathbf{K}_r (\mathbf{c}^n) \mathbf{a}^n + \nu \mathbf{P}_r \mathbf{a}^n + \mathbf{s}_r) - \Delta t \mathbf{G}_r \mathbf{b}^{n+1},
\end{equation}
Equation~\eqref{eq:rom_uf1} is unique to the consistent flux method:
a standard Galerkin ROM without face-flux unknowns would not include it.
It advances the face-flux POD coefficients $\mathbf{c}^{n+1}$ by solving
a linear system whose left-hand side involves the face-flux mass matrix
$\mathbf{W}_r = \mathbf{\Psi}^T \mathbf{\Sigma}\,\mathbf{\Psi}$.
The term $\mathbf{N}_r\,\mathbf{a}^n =
\mathbf{\Psi}^T\mathbf{\Sigma}\,\mathbf{I}_{p\rightarrow f}\,\mathbf{\Phi}\,\mathbf{a}^n$
maps the cell-centred velocity onto the face-flux POD space, ensuring
that the face fluxes remain consistent with the cell-centred velocity
throughout the simulation. The same pressure correction
$-\Delta t\,\mathbf{G}_r\,\mathbf{b}^{n+1}$ that appears in
equation~\eqref{eq:rom_u1} is applied here at the face level, coupling
all three equations through the shared pressure $\mathbf{b}^{n+1}$.

with
\begin{itemize}
    \item $\mathbf{W}_r = \mathbf{\Psi}^T \mathbf{\Sigma} \mathbf{\Psi} \in \mathbb{R}^{N_r^u \times N_r^u}$,
    \item $\mathbf{N}_r = \mathbf{\Psi}^T \mathbf{\Sigma} \mathbf{I}_{p\rightarrow{f}} \mathbf{\Phi} \in \mathbb{R}^{N_r^u \times N_r^u}$,
    \item $\mathbf{P}_r = \mathbf{\Psi}^T \mathbf{\Sigma} \mathbf{D}_p \mathbf{\Phi} \in \mathbb{R}^{N_r^u \times N_r^u}$, 
    \item $\mathbf{G}_r = \mathbf{\Psi}^T \mathbf{\Sigma} \mathbf{G}_f \mathbf{X} \in \mathbb{R}^{N_r^u \times N_r^p}$,
    \item $\mathbf{s}_r = \mathbf{\Psi}^T \mathbf{\Sigma} \mathbf{r}_p \in \mathbb{R}^{N_r^u}$.
\end{itemize}
where $\mathbf{\Sigma} \in \mathbb{R}^{dm \times dm}$ is the matrix that stores the face-area information associated with the computational cells.

The last one is the reduced vector and the other are matrices. The reduced convection terms $\mathbf{A}_r (\mathbf{a}^n) \in \mathbb{R}^{N_r^u \times N_r^u \times N_r^u}$, $\mathbf{C}_r (\mathbf{a}^n) \in \mathbb{R}^{N_r^u \times N_r^u \times N_r^u}$ and $\mathbf{K}_r (\mathbf{a}^n) \in \mathbb{R}^{N_r^u \times N_r^u \times N_r^u}$ are determined, respectively, by 
\begin{equation*}
    \mathbf{A}_{r,i} = \mathbf{X}^T \mathbf{V} \mathbf{M}_p \tilde{\mathbf{C}_p} (\mathbf{\Psi}_i) \mathbf{\Phi},
\end{equation*}
\begin{equation*}
    \mathbf{C}_{r,i} = \mathbf{\Phi}^T \mathbf{V} \tilde{\mathbf{C}_p} (\mathbf{\Psi}_i) \mathbf{\Phi},
\end{equation*}
\begin{equation*}
    \mathbf{K}_{r,i} = \mathbf{\Psi}^T \mathbf{\Sigma} \tilde{\mathbf{C}_p} (\mathbf{\Psi}_i) \mathbf{\Phi}.
\end{equation*}
In this way, the face-flux equation enforces the discrete mass-conservation constraint at the
reduced level through the face-centered variables. The cell-centered velocity POD basis is not
assumed to be exactly divergence-free on the collocated grid; therefore, the reduced pressure
equation and the face-flux constraint are retained in the ROM formulation to preserve the
pressure--velocity coupling and control the residual divergence.


The algebraic construction of the discretize-then-project framework is not
restricted to the cavity geometry: the POD basis is extracted from discrete
snapshot degrees of freedom, and the Galerkin projection operates on
finite-volume operators assembled at cell and face level. The same formulation
can therefore, in principle, be applied to structured, unstructured, and
body-fitted meshes. Demonstrating its accuracy and robustness on complex
geometries nevertheless requires case-specific validation, in the same spirit as POD--Galerkin approaches that have been successfully demonstrated on non-trivial geometries, including compressible and incompressible flows around airfoils~\cite{BourguetBrazaDervieux2011,ColaneraEtAl2024Airfoil} 
and three-dimensional patient-specific vascular 
networks~\cite{BallarinEtAl2016CABG}.

\subsection{Closure Model}
In recent years, the evolution of machine learning techniques has increasingly highlighted the role of data-driven methods. This term refers to approaches in which decisions, predictions, and analyses are not guided by predefined theoretical models or manually designed rules, but instead emerge directly from the available data. In other words, rather than describing in advance how a system should behave, the machine is allowed to “learn” from past observations, discovering patterns and regularities that can be used to tackle new cases. This shift in perspective has revolutionized fields ranging from medicine to finance, from robotics to natural language processing.

Data-driven approaches rely on the availability of data rather than on explicitly prescribed physical models. In practice, this often leads to the use of artificial neural networks, which are able to approximate complex input–output relationships when sufficient data are available. Over the years, a variety of network architectures have been proposed, and their suitability depends strongly on the structure of the problem under consideration. In the present work, attention is restricted to three commonly used architectures, namely MLPs, LSTMs, and Transformer-based models.

An MLP is used as a simple neural baseline. The model is implemented as a sequence of fully connected layers with nonlinear activations and is trained to directly relate the reduced input variables to the target quantities. This architecture is known to work reasonably well for static regression problems. However, since no explicit mechanism is included to represent temporal structure, its effectiveness decreases when the solution depends strongly on the system dynamics over time, as is the case for unsteady flow fields.

To model the temporal evolution of the turbulent viscosity coefficients,
recurrent neural network architectures are a natural choice, as the
reduced-order dynamics exhibit strong temporal correlations
\cite{hajisharifi2024lstm}. Among these architectures, LSTM networks are particularly suitable due to their ability to propagate information over extended time horizons in a stable manner. This property is essential in the present context, where the turbulent closure depends not only on the instantaneous reduced state but also on its recent temporal history. Compared to simpler feedforward models, the LSTM formulation provides improved robustness when learning long-term dependencies, which is a known limitation of standard recurrent networks.


Transformer architectures are also considered in this work as a possible alternative for learning temporal dependencies in the reduced-order setting. Their attention-based formulation allows interactions between distant time instants to be represented without relying on recurrence, which can be beneficial when long-range correlations are present in the reduced dynamics. At the same time, this design leads to highly parallel training procedures, although its effectiveness in the present turbulent viscosity prediction task is not guaranteed a priori. For this reason, Transformers are included in the comparative analysis alongside recurrent and feedforward models.



\subsubsection{Neural Networks}
In this work, as previously discussed, a data-driven strategy is employed to address the physical inconsistency encountered when applying standard projection techniques to the turbulent flows. While the 'discretize-then-project' methodology ensures consistency for velocity and pressure, it has been observed that standard Galerkin projection fails to yield a physically coherent representation of the turbulent viscosity field in the reduced space. This occurs because the intrusive projection is unable to accurately resolve the complex temporal evolution and non-linear interactions of the eddy viscosity within the reduced basis. To overcome this limitation, we adopt a hybrid approach: velocity and pressure are resolved via intrusive projection, while a LSTM neural network is used as a non-intrusive closure to predict the turbulent viscosity coefficients. The network is trained using the temporal coefficients of velocity $u$ and pressure $p$ obtained from high-fidelity FOM simulations. More specifically, the neural network learns the nonlinear relationship between the reduced velocity and pressure states and the corresponding turbulent viscosity coefficients extracted from the FOM snapshot database.

Following the same logic, we also evaluated two other architectures, MLP and Transformer, to provide a comprehensive comparison. For the configurations and datasets investigated here, the LSTM provided
the most consistent turbulent-viscosity reconstruction across the tested
mode counts. This empirical behavior is compatible with the use of temporal
memory, although the present comparisons do not isolate the architecture
feature responsible for the observed differences.; these results are detailed in the following chapter. All hyperparameters and architectural details for these three networks are reported in the Appendix, see \autoref{tab:hyperparameters} and \autoref{tab:layer_features}.

\subsubsection{Data driven method for turbulent viscosity prediction}

As discussed in the previous subsection, the turbulent viscosity field cannot be accurately represented through a purely intrusive projection strategy. Consequently, a data-driven closure is introduced to reconstruct the turbulent viscosity coefficients within the reduced-order framework. The proposed closure is designed for eddy-viscosity turbulence models, such as the Smagorinsky LES model adopted in this work, and is based on the Boussinesq hypothesis, according to which the effect of unresolved turbulent fluctuations can be represented through an eddy viscosity $\nu_T$. Before describing the offline and online stages of the methodology, it is useful to recall how the turbulent viscosity enters the governing equations.

Using the FOM, we obtain the three primary fields of interest: velocity $\mathbf{u}$, pressure $p$, and turbulent viscosity $\nu_T$. The role of the turbulent viscosity within the projection framework can be summarized as follows.

For clarity, we briefly recall the corresponding full-order formulation in continuous notation.

First, we define a forcing term, $F_{aux}$, which represents the variation in the velocity field due to convective and diffusive effects over a time step $\Delta t$:
\begin{equation*}
    F_{aux} = \Delta t [-\nabla \cdot (\varphi U) + \nabla \cdot (\nu_{eff}\nabla U)];
\end{equation*}

where $-\nabla \cdot (\varphi u)$ represents the nonlinear convective transport and $\nabla \cdot (\nu_{eff}\nabla u)$ denotes the viscous diffusive term. In this formulation, $\nu_{eff} = \nu + \nu_T$ is the effective viscosity, which combines the molecular kinematic viscosity $\nu$ and the turbulent eddy viscosity $\nu_T$.

In the present formulation, the turbulent viscosity is explicitly retained in the diffusive term of the momentum equation. This differs from earlier discretize-then-project approaches developed for laminar flows, such as the framework presented in \cite{https://doi.org/10.48550/arxiv.2010.06964}
, where the turbulent contribution was not incorporated at the projection level. As a result, the present approach provides a more consistent treatment of diffusive effects when turbulent dynamics are considered.

The velocity field is then updated to obtain an intermediate state, $U_{aux}^{(1)}$, which accounts for forcing but has not yet been corrected for pressure effects:

\begin{equation*}
    U_{aux}^{(1)} = U + F_{aux};
\end{equation*}

Subsequently, a pressure Poisson equation (PPE) is solved to enforce the incompressibility constraint:

\begin{equation*}
    \nabla ^{2} p = \frac{1}{\Delta t} \nabla \cdot U_{aux}^{(1)};
\end{equation*}

\begin{equation*}
    p(cellRef) = p_{ref};
\end{equation*}

Once the pressure field is obtained, the final velocity is corrected:

\begin{equation*}
    U = U_{aux}^{(1)} -\Delta t \nabla p;
\end{equation*}

This correction step enforces the incompressibility constraint by removing the divergence of the velocity field, while retaining the contribution of the turbulent eddy viscosity provided by the data-driven model.

The previous equations highlight that the turbulent viscosity directly influences the momentum update through the effective viscosity $\nu_{eff}$. Therefore, inaccuracies in the reconstruction of $\nu_T$ may propagate to the velocity and pressure fields. Rather than deriving and solving an additional reduced transport equation for the turbulent viscosity, the present work adopts a non-intrusive strategy in which the turbulent viscosity coefficients are inferred from the reduced flow state using a neural network. More specifically, the network receives the reduced velocity and pressure coefficients as inputs and predicts the corresponding turbulent viscosity coefficients. The predicted coefficients are subsequently used to reconstruct the turbulent viscosity field through the POD expansion introduced in Eq.~\ref{eq:nut_pod}, thereby providing the turbulent contribution required to evaluate the effective viscosity $\nu_{eff}=\nu+\nu_T$. The resulting methodology follows the standard offline--online decomposition adopted throughout this work.

The overall methodology can be divided into two stages. First, the full-order model is advanced using the consistent flux formulation described above in order to generate high-fidelity snapshots of velocity, pressure, and turbulent viscosity. The main operations involved in the generation of these snapshots are summarized in Steps 1--4 below. The resulting dataset is then employed during the offline stage of the hybrid ROM framework, where POD coefficients are extracted and used to train the neural-network closure (Steps 5--7). Finally, the trained model is deployed during the online ROM simulation (Step 8).

\begin{enumerate}
    \item \textbf{Forcing Field Computation:} Calculate the
    forcing term $F_{aux}$ by evaluating convective and diffusive
    transport, explicitly incorporating the effective viscosity
    $\nu_{eff} = \nu + \nu_T$ to account for turbulent stresses.

    \item \textbf{Intermediate Velocity Prediction:} Update
    the current velocity field to a provisional state
    $U_{aux}^{(1)}$ by applying the forcing term, initially
    neglecting pressure contributions.

    \item \textbf{Pressure Poisson Resolution:} Solve the
    Poisson equation for the normalized scalar pressure field
    $p$ to enforce the incompressibility constraint.

    \item \textbf{Velocity and Flux Correction:} Apply the
    computed pressure gradient $\nabla p$ to correct the
    intermediate velocity field, ensuring a divergence-free
    solution.

    \item \textbf{Modal Coefficient Extraction:} Employ
    \texttt{ITHACAutilities::getCoeffs} to extract the temporal
    POD coefficients for velocity, pressure, and turbulent
    viscosity from the FOM snapshot dataset.

    \item \textbf{Data Preprocessing and Normalization:}
    Normalize the extracted coefficients using a StandardScaler,
    assemble input sequences (velocity and pressure coefficients)
    and target sequences (turbulent viscosity coefficients), and
    store the scaler for later inverse transformation.

    \item \textbf{Neural Network Training:} Execute
    \texttt{model.fit} to train the selected architecture
    (LSTM, MLP, or Transformer), learning the nonlinear temporal
    mapping from the reduced flow state to the eddy viscosity
    coefficients.
\end{enumerate}

The single \textbf{online step} is then:

\begin{enumerate}
    \setcounter{enumi}{7}
    \item \textbf{Prediction and Reconstruction:} At each
    time step of the ROM simulation, deploy the trained model
    via \texttt{model.predict}, apply the inverse StandardScaler
    transform, and reconstruct the turbulent viscosity field
    according to Eq.~\ref{eq:nut_pod}. The reconstructed
    turbulent viscosity field is combined with the molecular
    viscosity through $\nu_{eff}=\nu+\nu_T$ before being used
    in the reduced momentum equations. The ROM then advances
    the velocity and pressure coefficients using only
    precomputed reduced operators, without any further
    FOM evaluation.
\end{enumerate}

This procedure integrates the consistent flux formulation used for velocity and pressure with a non-intrusive, data-driven treatment of the turbulent closure. A schematic overview of the workflow is provided in Figure~\ref{graphical abstract}.

\begin{figure}[H]
    \centering
    \resizebox{0.95\linewidth}{!}{%
    \begin{tikzpicture}[
        node distance=0.7cm and 1cm,
        font=\sffamily\footnotesize,
        base box/.style={
            rectangle, 
            rounded corners=3pt, 
            minimum width=3.8cm, 
            minimum height=0.9cm, 
            align=center, 
            draw=darkGray, 
            thick,
            fill=white,
            inner sep=5pt
        },
        blue box/.style={base box, fill=fomBlue!20, draw=fomBlue},
        orange box/.style={base box, fill=offlineOrange!20, draw=offlineOrange},
        green box/.style={base box, fill=nnGreen!20, draw=nnGreen},
        purple box/.style={base box, fill=onlinePurple!20, draw=onlinePurple},
        group box/.style={
            rectangle, 
            rounded corners=5pt, 
            draw=darkGray!50, 
            dashed, 
            inner sep=8pt
        },
        arrow/.style={
            ->, 
            >={Latex[length=3mm, width=2mm]}, 
            thick, 
            darkGray, 
            rounded corners=5pt
        }
    ]

    \node[blue box] (openfoam) {\textbf{OpenFOAM}};
    \node[blue box, right=1.0cm of openfoam] (ithaca) {\textbf{ITHACA-FV}\\ROM Library};

    \coordinate (top_center) at ($(openfoam.east)!0.5!(ithaca.west)$);
    \node[blue box, below=0.6cm of top_center] (fom_solver) {
        \textbf{FOM Solver (CFD)}\\
        LES: 3D Navier-Stokes \ref{eq:continuity_equation}, \ref{eq:momentum_equation}\\
        Smagorinsky SGS Model \ref{eq:effective_viscosity}, \ref{eq:smagorinsky}\\
        \tikz{
            \node[draw=fomBlue, fill=fomBlue!20, rounded corners=2pt, align=center, inner sep=2pt, font=\scriptsize]
            {
                Spatial discretization: FVM on Collocated Grids\\
                Time integration / pressure-velocity coupling: Projection Method + CFM
            };
        }
    };
    
    \begin{scope}[on background layer]
        \node[group box, draw=darkGray, thick, fill=fomBlue!5, fit=(openfoam)(ithaca)(fom_solver), label={[text=darkGray, font=\bfseries]north:FOM}] (top_group) {};
    \end{scope}


    \node[orange box, cylinder, shape border rotate=90, aspect=0.25,
    below=1.6cm of fom_solver] (snapshots)
    {\textbf{Snapshot Collection}\\
    CFM Eqs.~\ref{eq:L_f}, \ref{eq:u_p_n_1}, \ref{eq:u_f}\\
    Smagorinsky Eq.~\ref{eq:smagorinsky}\\
    $u_p,\;p_p,\;u_f,\;\nu_T$};

    \node[orange box, below=of snapshots] (pod) {
    \textbf{Reduction}\\
    POD Basis Generation\\
    (Method of Snapshots / SVD)\\
    Generate Spatial Modes:\\
    $\Phi$: cell-centered velocity $u_p$\\
    $X$: cell-centered pressure $p_p$\\
    $\Psi$: face velocity        $u_f$\\
    $\Xi$: turbulent viscosity $\nu_T$
    };
    

    \node[orange box, below=of pod] (projection) {\textbf{Project Snapshots onto Basis}\\
    Orthogonal Projection\\
    $u_p \approx \Phi a$ \ref{eq:u_p}\\
    $p_p \approx X b$ \ref{eq:p_p}\\
    $u_f \approx \Psi c$ \ref{eq:u_f_approx}\\
    $\nu_T \approx \Xi d$ \ref{eq:nut_pod}\\
    Extract Temporal Coefficients $(a,b,c,d)$};


    \node[orange box, below=of projection] (galerkin) {
    \textbf{Galerkin Projection}\\
    Insert Eqs.~\ref{eq:u_p}, \ref{eq:p_p}, \ref{eq:u_f_approx}\\
    into CFM Eqs.~\ref{eq:L_f}, \ref{eq:u_p_n_1}, \ref{eq:u_f}\\
    Project discrete FOM operators in POD Bases\\
    $\Phi,\;X,\;\Psi$\\
    Strategy: discretize-then-Project using CFM};
    

    \node[orange box, below=of galerkin] (rom) {
    \textbf{Build Intrusive ROM}\\
    Solve reduced CFM System for temporal coefficients\\
    $(a,b,c)$ via Eqs.~\ref{eq:rom_p1}, \ref{eq:reference_pressure_constraint}, \ref{eq:rom_u1}, \ref{eq:rom_uf1}\\
    };
    \node[above=0.2cm of snapshots, font=\bfseries\color{offlineOrange}] (offline_title) {Intrusive Path (Physics-Based)};
    
    \begin{scope}[on background layer]
        \node[group box, draw=red, thick, fill=offlineOrange!5, fit=(offline_title)(rom)] (offline_group) {};
    \end{scope}


    \node[green box, chamfered rectangle, chamfered rectangle xsep=5pt,
    right=3.5cm of projection] (nn_training)
    {\textbf{Neural Network Training}\\
    Non-Intrusive Turbulence Closure\\
    Input: POD Coefficients $(a,b)$\\
    Target: Turbulent Viscosity Coefficients $d$\\
    Learn Mapping $\mathcal{F}_{NN}:(a,b)\rightarrow d$\\
    Architectures: LSTM, Transformer, MLP};

    \node[above=0.2cm of nn_training, font=\bfseries\color{nnGreen}] (nn_title) {Non-Intrusive Path (Training)};
    
    \begin{scope}[on background layer]
        \node[group box, draw=red, thick, fill=nnGreen!5, fit=(nn_title)(nn_training)] (nn_group) {};
    \end{scope}

    \node[fit=(offline_group)(nn_group), inner sep=0pt] (shared_offline) {};
    
    \node[text=red, font=\bfseries\large] (offline_text) at ($(nn_group.north) + (0, 1.2cm)$) {Offline Phase};
    
    \draw[->, red, thick, >=Latex] (offline_text.west) -- (offline_text.west -| offline_group.east);
    \draw[->, red, thick, >=Latex] (offline_text.south) -- (nn_group.north);

    \coordinate (middle_center) at ($(shared_offline.south)$);


    \node[purple box, below=2.0cm of rom] (integration) {
    \textbf{Hybrid Mechanism Coupling}\\
    ROM $(a,b,c)$\\
    + NN Closure $(a,b)\rightarrow d$\\
    Reconstruct $\nu_T \approx \Xi d$\\
    \textit{Inject $\nu_T$ into reduced momentum model}
    };


    \node[green box, right=1.0cm of integration, align=center] (nn_mapping)
    {\textbf{Neural Network Inference}\\
    $\mathcal{F}_{NN}:(a,b)\rightarrow d$};

    \node[purple box, rounded rectangle, below=of integration] (prediction) {\textbf{Prediction and Reconstruction}\\
    $u_p \approx \Phi a, \ p_p \approx Xb, \ u_f \approx \Psi c, \nu_T \approx \Xi d$\\
    (Fast \& Consistent Prediction)};

    \coordinate (online_center) at ($(integration.north)!0.5!(nn_mapping.north)$);
    \node[above=0.2cm of online_center, font=\bfseries\color{onlinePurple}] (online_title) {Testing/Inference};
    
    \begin{scope}[on background layer]
        \node[group box, draw=onlinePurple, thick, fill=onlinePurple!5, fit=(online_title)(prediction)(integration)(nn_mapping), label={[text=onlinePurple, font=\bfseries]north:Online Phase}] (online_group) {};
    \end{scope}

    \draw[arrow] (openfoam.south) -- (fom_solver.north -| openfoam.south);
    \draw[arrow] (ithaca.south) -- (fom_solver.north -| ithaca.south);
    \draw[arrow] (fom_solver.south) -- (snapshots.north);
    \draw[arrow] (snapshots) -- (pod);
    \draw[arrow] (pod) -- (projection);
    \draw[arrow] (projection) -- (galerkin);
    \draw[arrow] (galerkin) -- (rom);
    
    \draw[arrow] (projection.east) -- (nn_training.west) node[midway, above, font=\scriptsize, darkGray] {Temporal Coeffs ($a,b,d$)};
    
    \draw[arrow] (rom.south) -- (integration.north);
    
    \path ([xshift=0.6cm]nn_mapping.east) coordinate (target_entry);
    \draw[arrow, dashed] (nn_training.south) -- ++(0,-0.5) coordinate (drop)
        -- (drop -| target_entry) 
        -- (target_entry) node[midway, right, align=left] {Transfer\\Weights}
        -- (nn_mapping.east);

    \draw[arrow] ([yshift=-0.15cm]nn_mapping.west) -- ([yshift=-0.15cm]integration.east) node[midway, below, font=\scriptsize] {$d$};
    
    \draw[arrow] ([yshift=0.15cm]integration.east) -- ([yshift=0.15cm]nn_mapping.west) node[midway, above, font=\scriptsize] {$a,b$};
    
    \draw[arrow] (integration.west) .. controls +(-1,0.5) and +(-1,-0.5) .. (integration.west) node[midway, left, font=\scriptsize, align=right] {$t \to t+1$\\Iterate};
    
    \draw[arrow] (integration.south) -- (prediction.north);

    \end{tikzpicture}
    }
    \caption{Detailed workflow of the proposed hybrid methodology combining a physics-based ”discretize-
then-project” ROM with a data-driven turbulence closure which is trained offline but integrated directly into the online time-stepping loop.}
    \label{graphical abstract}
\end{figure}

\section{Numerical Results} \label{sec:NumericalResults}

\subsection{Test case description}
This section outlines the numerical setup adopted in the present study. The benchmark problem considered is the classical lid-driven cavity, modeled in a three-dimensional domain. The simulations are performed using the consistent flux method, already available in ITHACA-FV, an open-source $C^{++}$ library built on top of OpenFOAM. All numerical solvers rely on the OpenFOAM framework. The objective is to analyze the turbulent regime of this configuration, focusing on the reconstruction of velocity and pressure fields through ROM techniques, while turbulent viscosity is recovered via a data-driven approach.

\begin{figure}[H]
\centering
\begin{tikzpicture}[scale=0.65, line join=round]
  \definecolor{CubeBlue}{RGB}{70,130,180}
  \definecolor{LidBlue}{RGB}{173,216,230}
  \definecolor{ArrowRed}{RGB}{180,40,40}

  \def\a{3.2}   
  \def\d{1.2}   

  \coordinate (A)  at (0,0);
  \coordinate (B)  at (\a,0);
  \coordinate (C)  at (\a,\a);
  \coordinate (D)  at (0,\a);
  \coordinate (A2) at (\d,\d);
  \coordinate (B2) at (\a+\d,\d);
  \coordinate (C2) at (\a+\d,\a+\d);
  \coordinate (D2) at (\d,\a+\d);

  \fill[CubeBlue!8]  (B)--(B2)--(C2)--(C)--cycle;
  \fill[CubeBlue!5]  (D)--(D2)--(C2)--(C)--cycle;
  \fill[CubeBlue!4]  (A)--(B)--(C)--(D)--cycle;

  \draw[gray!50, dashed, thin] (A)--(A2);
  \draw[gray!50, dashed, thin] (A2)--(B2);
  \draw[gray!50, dashed, thin] (A2)--(D2);

  \draw[CubeBlue, thick] (A)--(B)--(C)--(D)--cycle;
  \draw[CubeBlue, thick] (B)--(B2)--(C2)--(C);
  \draw[CubeBlue, thick] (C2)--(D2)--(D);
  \draw[CubeBlue, thick] (B2)--(C2);
  \draw[CubeBlue, thick] (D2)--(D);

  \fill[LidBlue, opacity=0.55] (D)--(C)--(C2)--(D2)--cycle;
  \draw[CubeBlue!80, thick]    (D)--(C)--(C2)--(D2)--cycle;

  \draw[->, ArrowRed, very thick, >=latex]
    ($(D)!0.15!(C)+(0,0.25)$) -- ($(D)!0.82!(C)+(0,0.25)$)
    node[midway, above, font=\small, text=ArrowRed]
    {$\mathbf{u}_{\mathrm{lid}} = (1,\,0,\,0)$\,m/s};

  \node[font=\small, text=gray!70!black]
    at (\a*0.5, \a*0.3) {no-slip walls};

  \draw[<->, gray!60, thin] (0,-0.35) -- (\a,-0.35);
  \node[font=\small, text=gray!60!black, fill=white, inner sep=1pt]
    at (\a*0.5,-0.35) {$L = 1$\,m};

  \node[font=\small, text=gray!60!black, anchor=west]
    at (\a+\d+0.15, \d*0.5) {$50^3$ cells};

  \node[font=\small, anchor=north west]
    at (\a+\d+0.1, 0.0)
    {\begin{tabular}{@{}l@{}}
       $Re = 100$: $\nu = 0.01$\,m$^2$/s\\[2pt]
       $Re = 3200$: $\nu = 3.125\!\times\!10^{-4}$\,m$^2$/s
     \end{tabular}};

\end{tikzpicture}
\caption{\small Three-dimensional lid-driven cavity. The top face 
(shaded) moves with a uniform tangential velocity 
$\mathbf{u}_{\mathrm{lid}} = (1,0,0)$\,m/s; all other walls are 
subject to no-slip conditions. The domain is a unit cube ($L = 1$\,m) 
discretized with a structured $50\times50\times50$ hexahedral mesh. 
Two Reynolds numbers are considered: $Re = 100$ 
($\nu = 0.01$\,m$^2$/s) and $Re = 3200$ 
($\nu = 3.125\times10^{-4}$\,m$^2$/s).}
\label{fig:cavity_schematic}
\end{figure}

The computational domain corresponds to a unit cube with edge length equal to 1 m (Figure~\ref{fig:cavity_schematic}). A structured mesh of $50 \times 50 \times 50$ hexahedral cells is employed. On the top boundary, a uniform tangential velocity $U=(1,0,0)$ m/s is prescribed, whereas no-slip boundary conditions are imposed on the remaining walls. The kinematic viscosity is set to $\nu = 0.01$ m$^2$/s, the cavity edge length is $L = 1$ m, giving a Reynolds number of $Re = UL/\nu = 100$. Turbulence is modeled through a turbulent viscosity $\nu_T$, using a LES simulation with the Smagorinsky turbulence model: a uniform zero condition is enforced at the walls, and a zero-gradient condition is prescribed for pressure. In particular, the Smagorinsky constant is $C_s =0.2$ and the turbulent Prandtl number is $\text{Pr} = 0.9$. 

For the $Re = 100$ case, the simulation is advanced for 2000 time steps 
($\Delta t = 0.005$\,s, total time $T = 10$\,s). The discretization 
schemes are chosen as follows: Euler forward for the temporal derivative, 
Gauss linear for gradients, Gauss linearUpwind for the divergence and 
Gauss linear orthogonal for the Laplacian. In the reduced-order framework, 
10 POD modes are retained for each variable, namely velocity, pressure, 
and turbulent viscosity. The same discretization schemes are retained for the 
$Re = 3200$ case; the number of POD modes is examined through a dedicated sensitivity study in Section~\ref{sec:Re3200}.

The second test case, introduced in Section~\ref{sec:Re3200}, is conducted 
at $Re = 3200$ ($\nu = 3.125 \times 10^{-4}$\,m$^2$/s), a regime in 
which the lid-driven cavity exhibits sustained unsteady dynamics with 
non-negligible turbulent viscosity, in contrast to the near-steady behavior 
at $Re = 100$. According to Leriche and Gavrilakis~\cite{Leriche2000}, 
turbulent characteristics in the cubical cavity emerge near 
$Re \approx 3000$, placing the present second test case in a genuinely 
unsteady and partially turbulent flow regime. The simulation is advanced 
for $T = 4$\,s with internal time step $\Delta t = 10^{-4}$\,s and 
snapshots stored every $\Delta t_{\text{write}} = 2 \times 10^{-3}$\,s, 
yielding 2000 snapshots for the ROM. The same mesh, boundary conditions, turbulence
model, and discretization schemes are retained for both cases. The number of POD
modes for the $Re = 3200$ case is selected through the sensitivity study reported
in Section~\ref{sec:Re3200}.


\subsection{Results}
In the first stage, the FOM was solved in order to generate the high-fidelity data required for the subsequent offline phase of the ROM process. Following the discretize-then-project approach, POD was applied to the FOM solutions to extract the dominant spatial modes associated with velocity, pressure, and turbulent viscosity fields. The corresponding temporal coefficients were also computed, representing the time evolution of each mode. Together, these components form the basis of the ROM, which provides a compact yet accurate representation of the original high-dimensional system. In the subsequent online phase, the velocity and pressure fields were reconstructed from the ROM.

Based on the reduced-order results, three neural network architectures  (LSTM, MLP, and Transformer) were trained and validated using the  temporal POD coefficients extracted from the FOM snapshots. For the  $Re = 100$ case, the first 1800 snapshots were employed for training,  while the remaining 200 snapshots were reserved for validation. The  objective was to establish a data-driven strategy for predicting the  turbulent viscosity field $\nu_T$, thereby enabling the resolution of  turbulent flow problems through the hybrid ROM--data-driven framework:  intrusive Galerkin projection for velocity and pressure, and a  non-intrusive data-driven closure for turbulent viscosity.

\subsubsection{Visualization, plots and errors}
To quantitatively validate the proposed approach, the obtained results are here presented and discussed. As a first step, the error between the FOM and the ROM was evaluated for different numbers of POD modes, with the aim of identifying the optimal number of modes to be used for this specific case study.

\begin{figure}[H]
\centering
\includegraphics[width=0.4\linewidth]{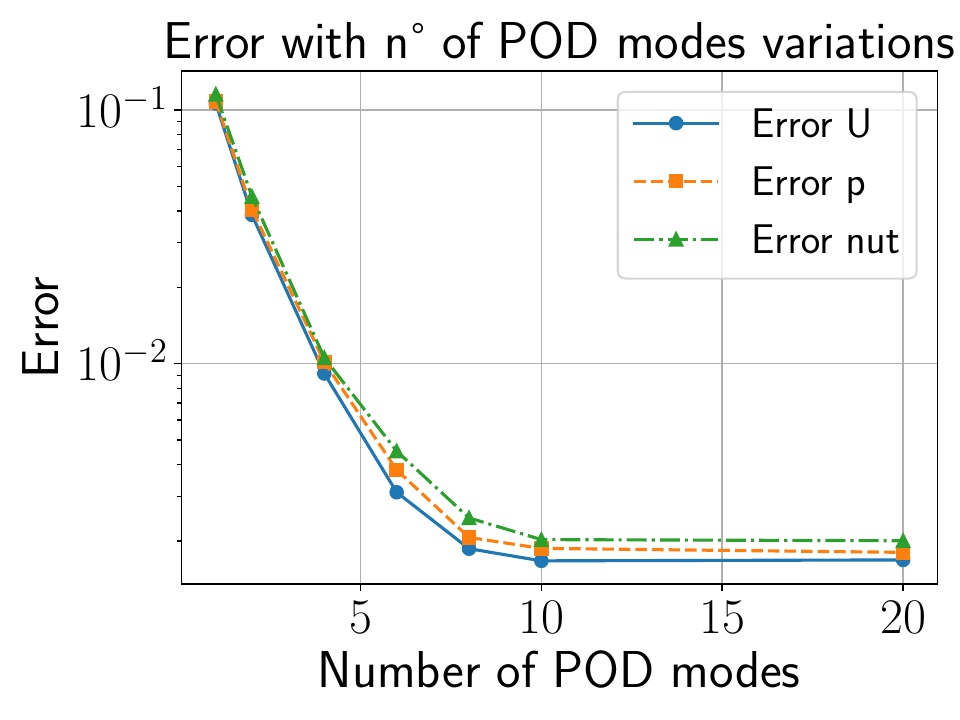}
\caption{\small Time-averaged relative $L^2$ error between 
the FOM and the ROM solutions as a function of the number 
of retained POD modes, for velocity $\mathbf{u}$, pressure 
$p$, and turbulent viscosity $\nu_T$. The ROM error does 
not necessarily decrease monotonically with the number of 
modes, as increasing the mode count may introduce 
instabilities in the reduced system.}
\label{errors_POD_modes_variation}
\end{figure}

From this analysis, shown in Figure~\ref{errors_POD_modes_variation}, it can be observed that the ROM error decreases rapidly with the number of modes up to $N = 10$ for all three fields, beyond which further improvement is marginal. 
It is important to note that the quantity shown is the ROM online error, the error obtained by running the full ROM solver with different numbers of retained POD modes, rather than the pure POD projection error. Unlike the projection error, which decreases monotonically by construction, the ROM online error does not necessarily decrease monotonically with the number of modes, since increasing the mode count may introduce instabilities in the reduced system and degrade accuracy beyond a certain threshold.

On the basis of this analysis, $N = 10$ modes per field 
was selected as the optimal value, providing an excellent 
balance between accuracy and stability of the reduced 
system. The total reduced-space dimension is $3 \times 
10 = 30$, corresponding to 10 POD modes retained for 
each of the three fields (velocity $\mathbf{u}$, pressure 
$p$, and turbulent viscosity $\nu_T$). For the $Re = 3200$ case, the POD truncation error (not the ROM 
online error) as a function of mode count is shown separately in 
Figure~\ref{fig:pod_errors_Re3200}.

\begin{equation}
\frac{
    \int_0^T \int_\Omega (F_{\text{FOM}} - F_{\text{ROM}}) dx dt}{\int_0^T \int_\Omega F_{\text{FOM}}dx dt}
\end{equation}


\begin{figure}[H]
\centering
\includegraphics[width=0.4\linewidth]{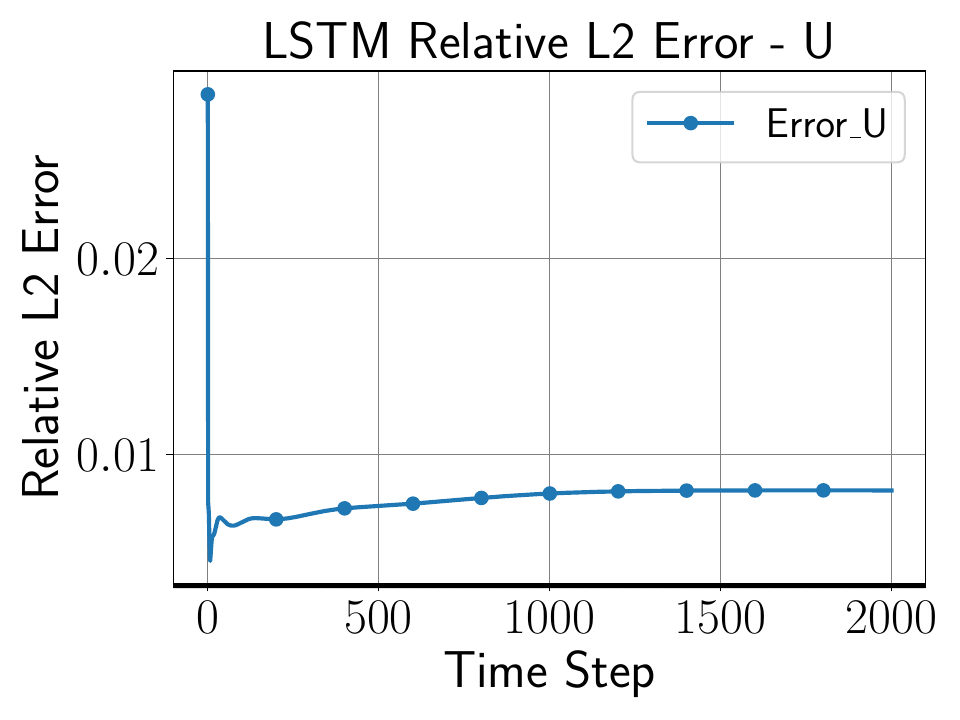}
\hfill
\includegraphics[width=0.4\linewidth]{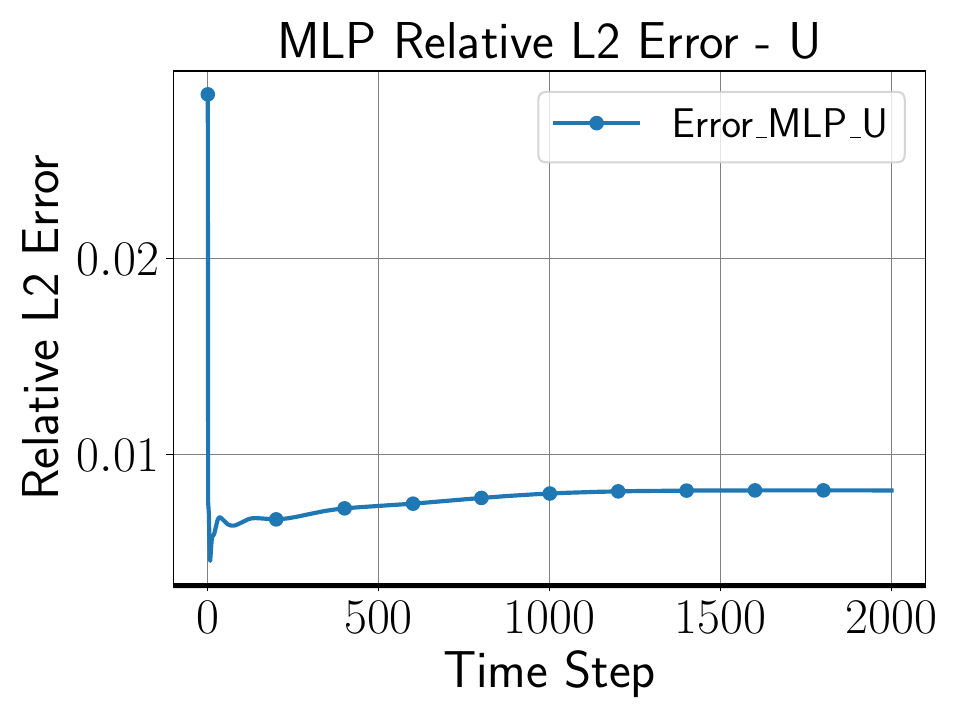}
\hfill
\includegraphics[width=0.4\linewidth]{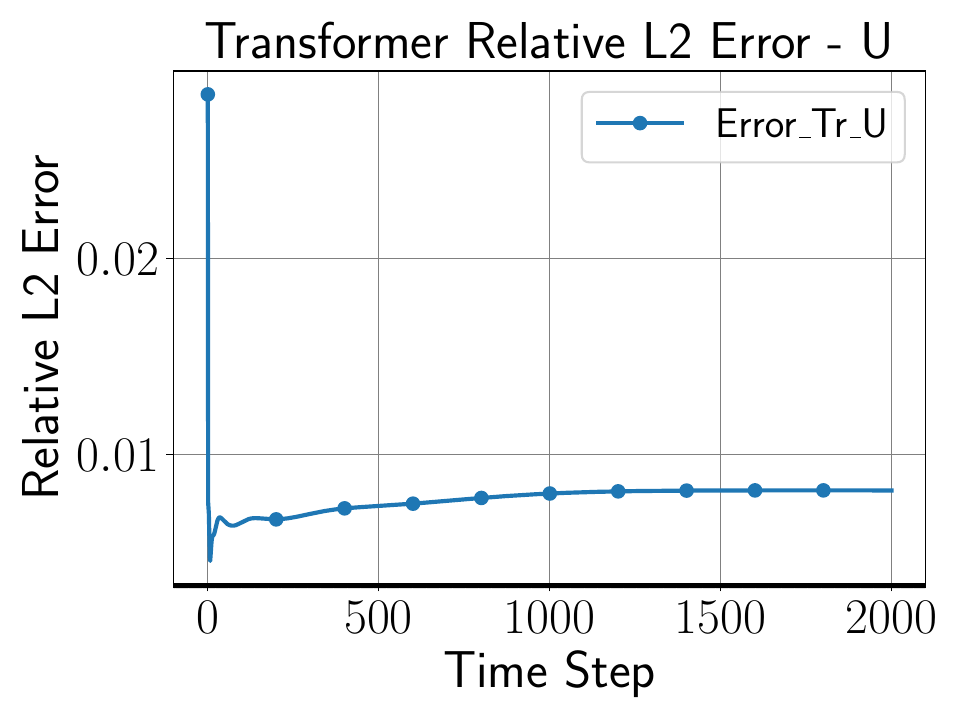}
\caption{\small Relative $L^2$ error between the FOM and ROM 
velocity field $\mathbf{u}$ as a function of time ($Re = 100$, 
$N = 10$ POD modes). Results shown for LSTM (left), MLP 
(centre), and Transformer (right).}
\label{FOMvROM_time_vel}
\end{figure}

\begin{figure}[H]
\centering
\includegraphics[width=0.4\linewidth]{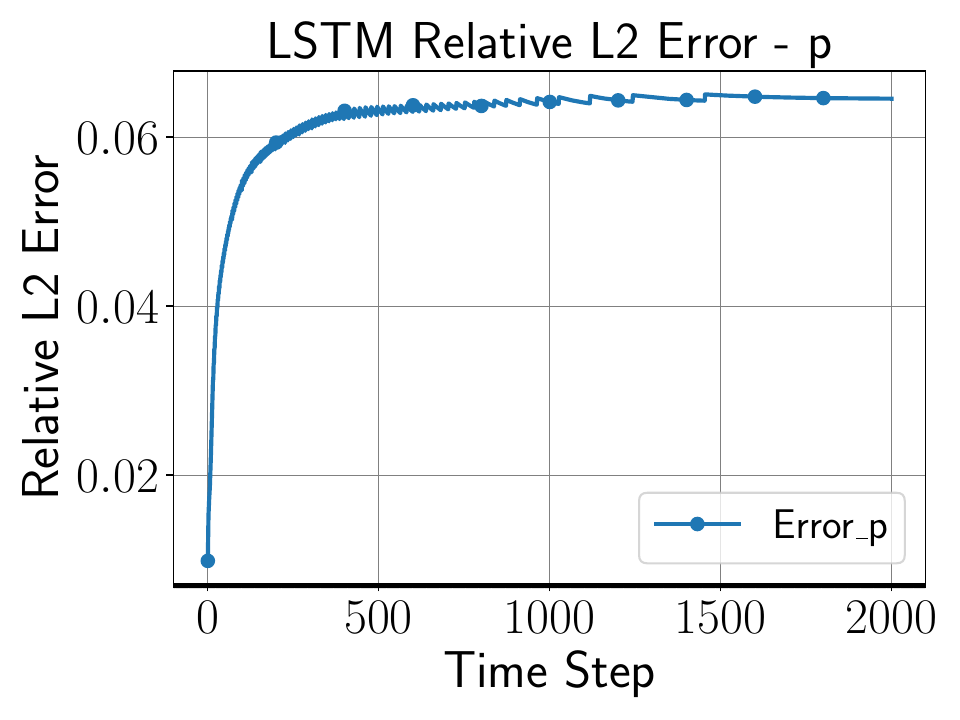}
\hfill
\includegraphics[width=0.4\linewidth]{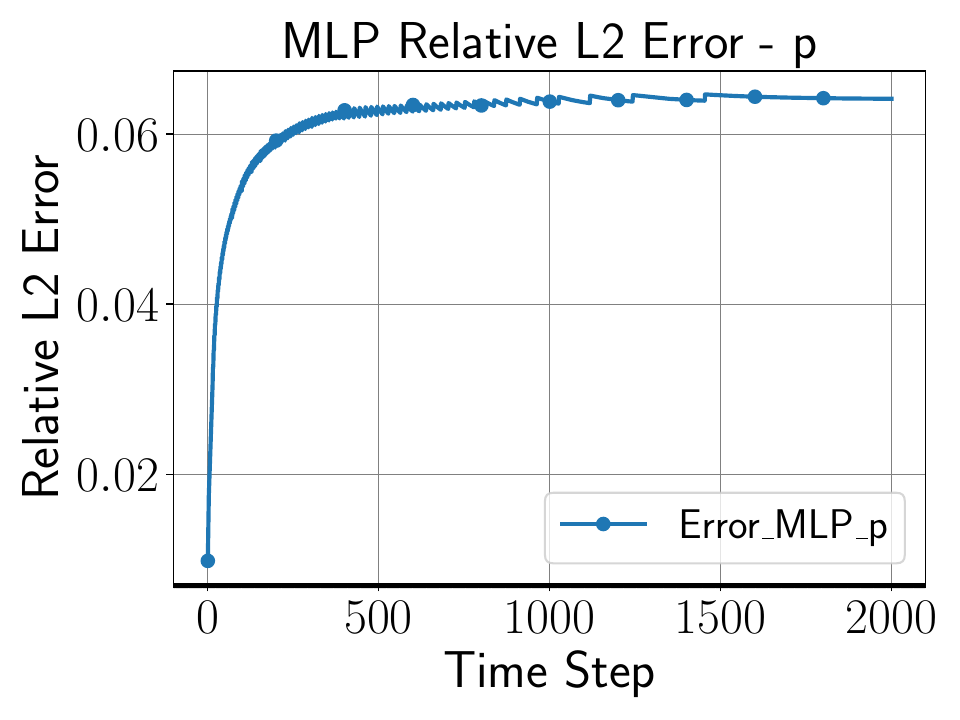}
\hfill
\includegraphics[width=0.4\linewidth]{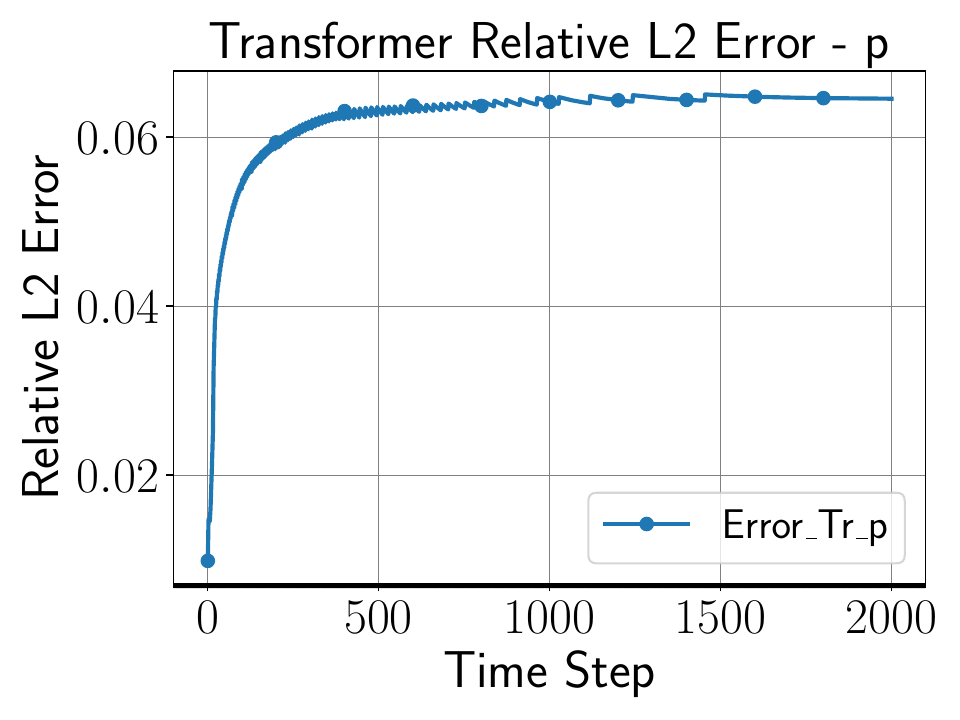}
\caption{\small Relative $L^2$ error between the FOM and ROM 
pressure field $p$ as a function of time ($Re = 100$, 
$N = 10$ POD modes). Results shown for LSTM (left), MLP 
(centre), and Transformer (right).}
\label{FOMvROM_time_pressure}
\end{figure}

\begin{figure}[H]
\centering
\includegraphics[width=0.4\linewidth]{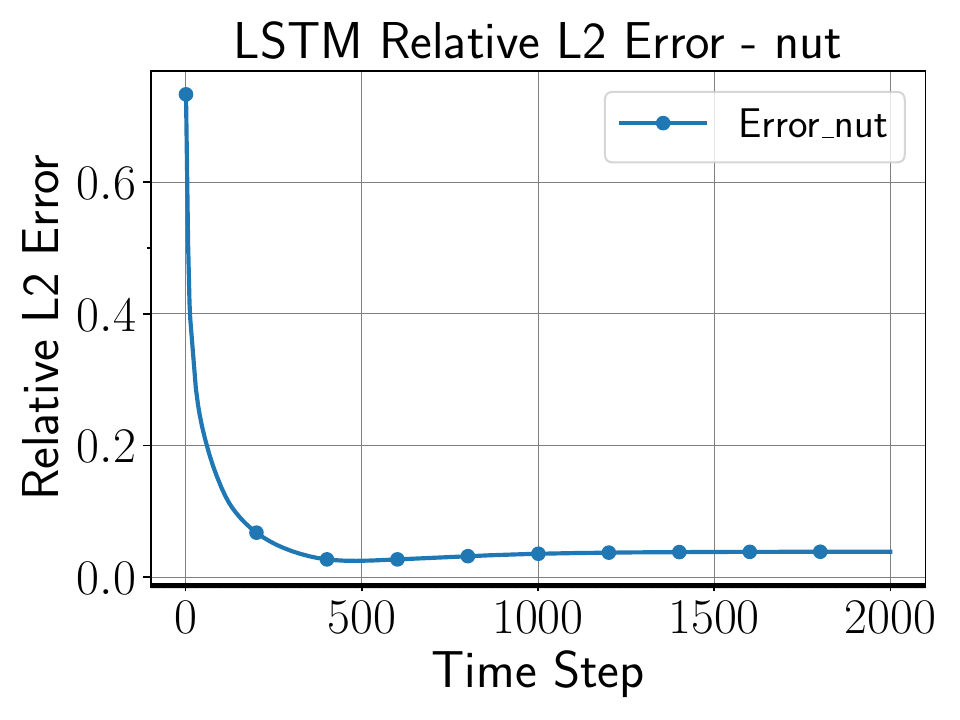}
\hfill
\includegraphics[width=0.4\linewidth]{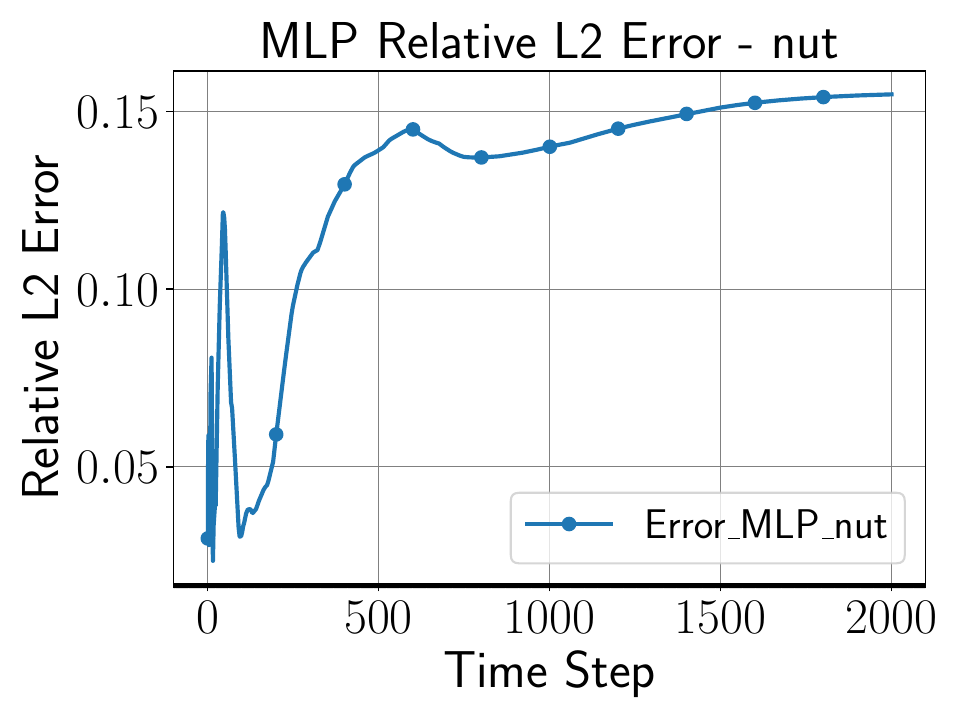}
\hfill
\includegraphics[width=0.4\linewidth]{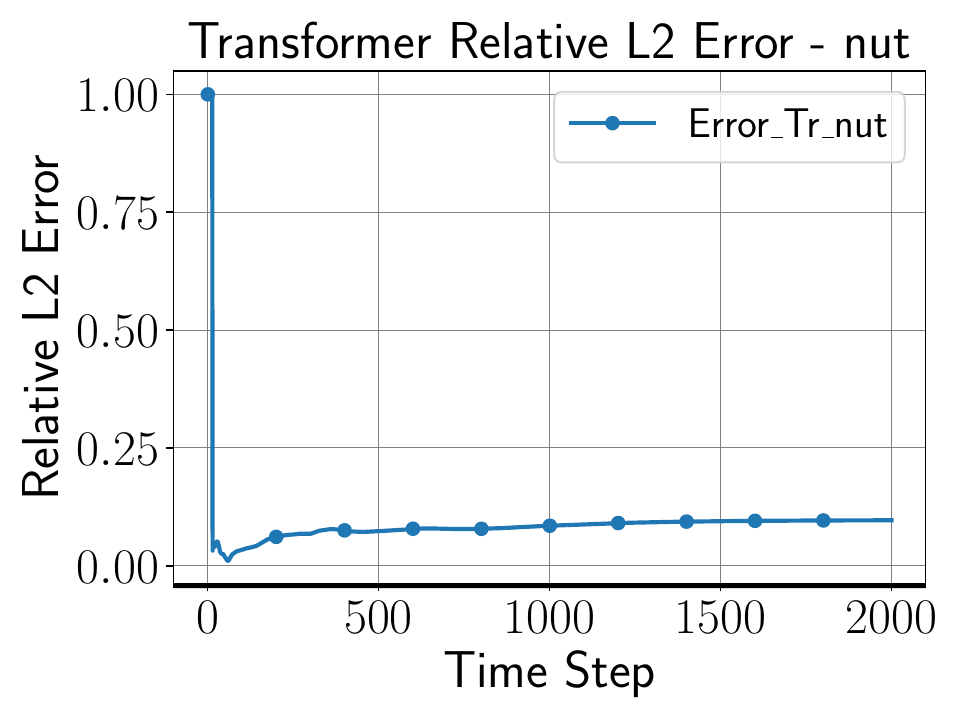}
\caption{\small Relative $L^2$ error between the FOM and ROM 
turbulent viscosity field $\nu_T$ as a function of time 
($Re = 100$, $N = 10$ POD modes). Results shown for LSTM 
(left), MLP (centre), and Transformer (right).}
\label{FOMvROM_time_turbulent}
\end{figure}

Having selected $N = 10$ POD modes per field as the optimal 
value (Figure~\ref{errors_POD_modes_variation}), 
Figures~\ref{FOMvROM_time_vel}, \ref{FOMvROM_time_pressure}, 
and~\ref{FOMvROM_time_turbulent} show the relative $L^2$ 
reconstruction errors between the FOM and ROM solutions as 
a function of time for velocity, pressure, and turbulent 
viscosity respectively, using all three neural network 
architectures (LSTM, MLP, Transformer) at $Re = 100$.

\begin{equation}
    \frac{\int_\Omega (F_{\text{FOM}} - F_{\text{ROM}}) dx}{\int_\Omega F_{\text{FOM}}dx}
\end{equation}

The error plots confirm that for the $Re = 100$ case the LSTM 
neural network provides the most accurate reconstruction across 
all three fields, consistent with its ability to capture temporal 
dependencies through its gated memory mechanism.

\begin{figure}[H] 
\centering
\includegraphics[width=0.5\linewidth]{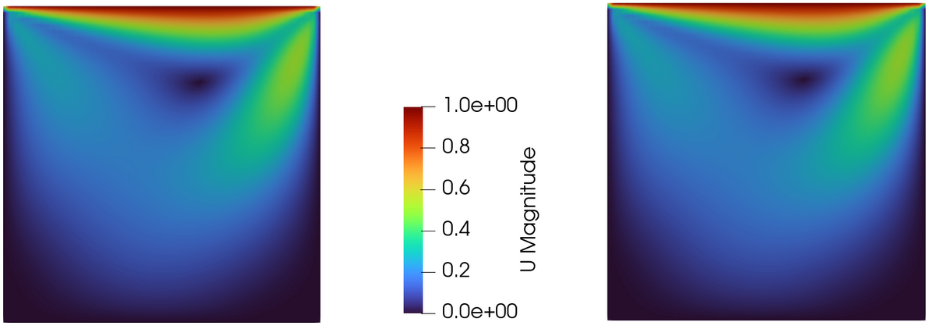}
\caption{\small Velocity magnitude $|\mathbf{u}|$ on a cross-section of the cavity at a representative time instant ($Re = 100$, LSTM, $N = 10$ POD modes). Left: FOM. Right: ROM.}
\label{FOMvROM_velocity}
\end{figure}

\begin{figure}[H] 
\centering
\includegraphics[width=0.5\linewidth]{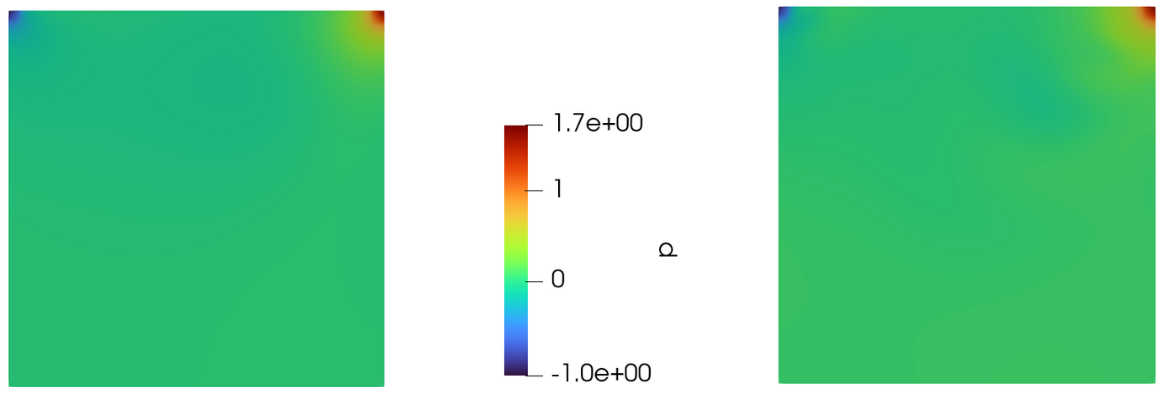}
\caption{\small Pressure field $p$ on a cross-section of the cavity at a representative time instant ($Re = 100$, LSTM, $N = 10$ POD modes). Left: FOM. Right: ROM.}

\label{FOMvROM_pressure}
\end{figure}

\begin{figure}[H]
\centering
\includegraphics[width=0.5\linewidth]{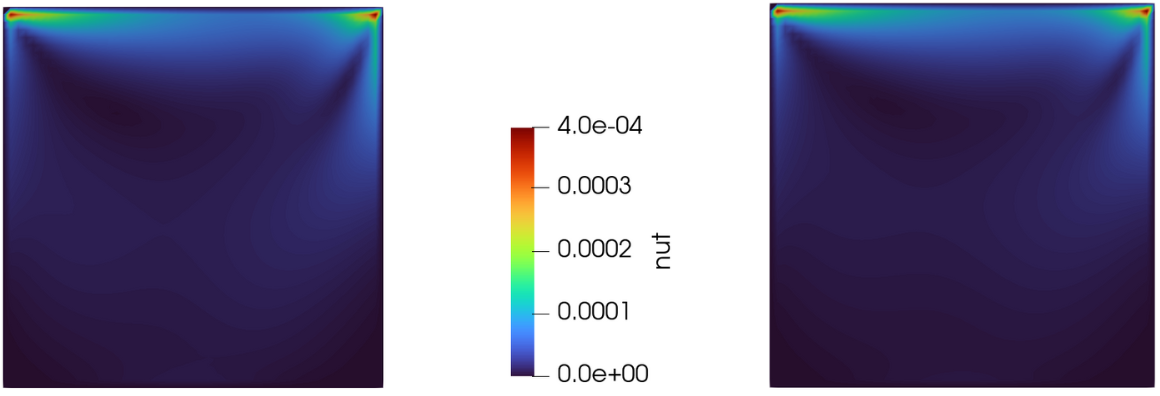}
\caption{\small Turbulent viscosity field $\nu_T$ on a cross-section of the cavity at a representative time instant ($Re = 100$, LSTM, $N = 10$ POD modes). Left: FOM. Right: ROM.}

\label{FOMvROM_turbulent_viscosity}
\end{figure}

To further illustrate the quality of the reconstruction, 
Figures~\ref{FOMvROM_velocity}, \ref{FOMvROM_pressure}, and
\ref{FOMvROM_turbulent_viscosity} show a comparison show a comparison 
between the FOM and ROM solutions at a representative 
time instant ($Re = 100$, LSTM closure, $N = 10$ POD modes).

The two solutions are in close visual agreement, confirming 
that the hybrid discretize-then-project ROM with data-driven 
closure accurately reproduces the FOM velocity, pressure, 
and turbulent viscosity fields at $Re = 100$.

\begin{figure}[H] 
\centering
\includegraphics[width=0.6\linewidth]{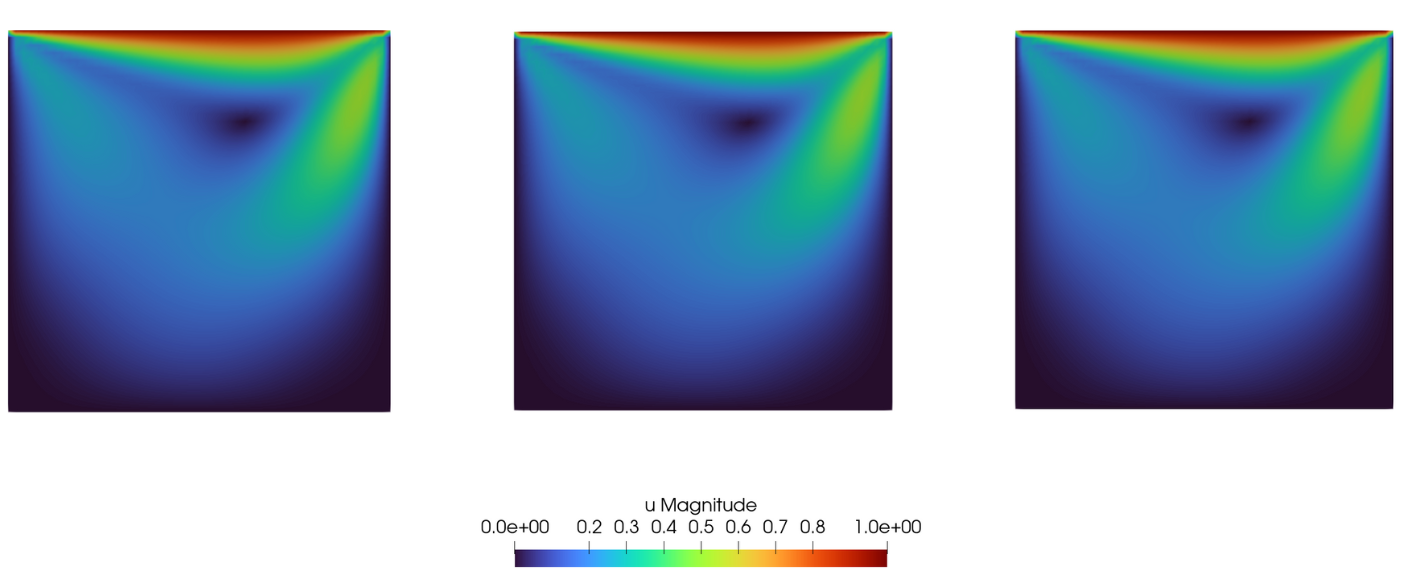}
\caption{\small Velocity magnitude $|\mathbf{u}|$ on the test-case section at a representative time instant (Re = 100, \(N = 10\) POD modes): LSTM (left), MLP (center), and Transformer (right).}
\label{confronto_NN_vel}
\end{figure}

\begin{figure}[H]
\centering
\includegraphics[width=0.6\linewidth]{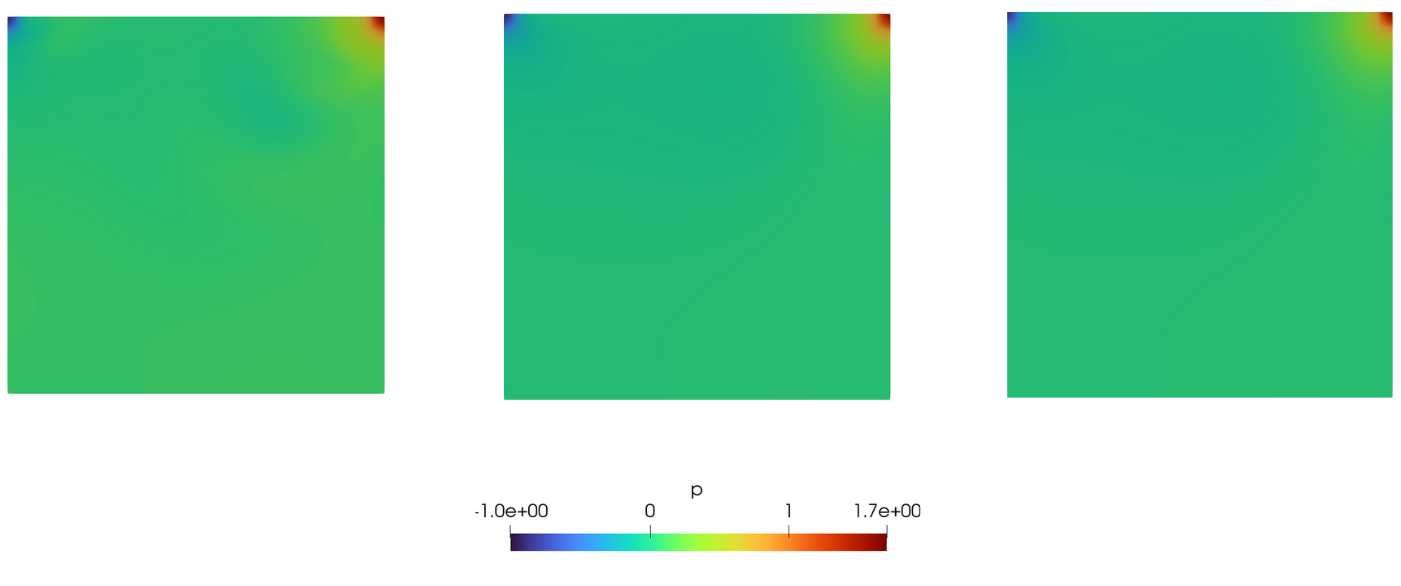}
\caption{\small Pressure field \(p\) on the test-case section at a representative time instant (Re = 100, \(N = 10\) POD modes): LSTM (left), MLP (center), and Transformer (right).}
\label{confronto_NN_pressure}
\end{figure}

\begin{figure}[H] 
\centering
\includegraphics[width=0.6\linewidth]{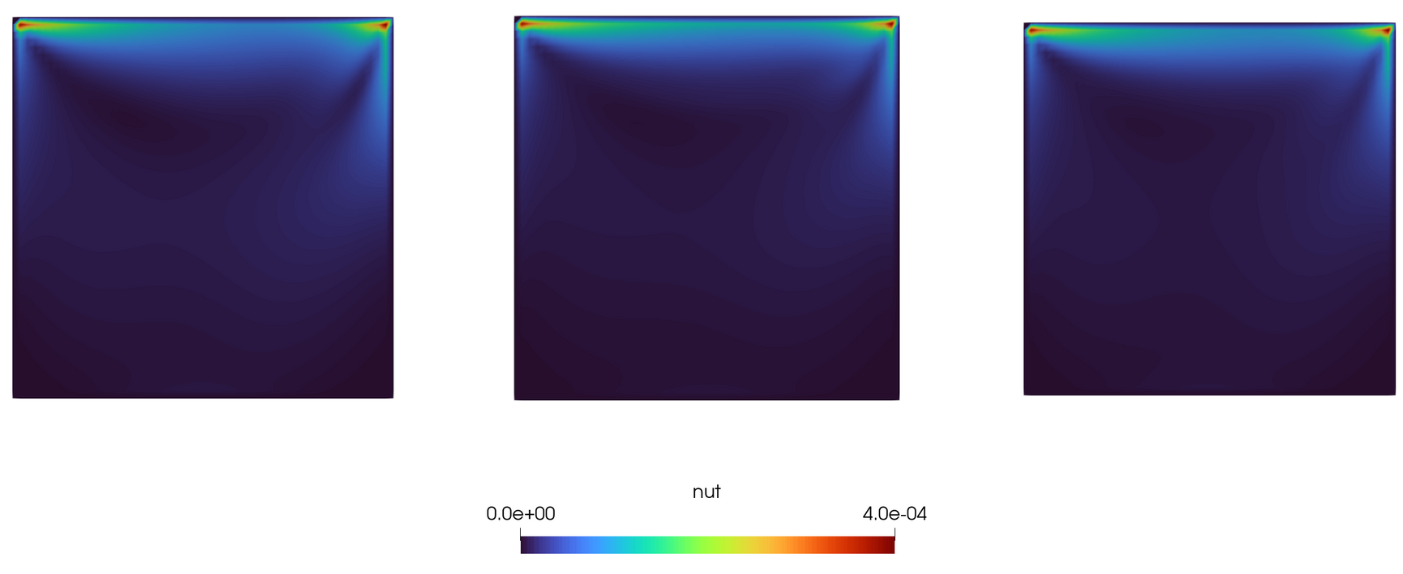}
\caption{\small Turbulent viscosity field \(\nu_T\) on the test-case section at a representative time instant (Re = 100, \(N = 10\) POD modes): LSTM (left), MLP (center), and Transformer (right).}
\label{confronto_NN_turbulent}
\end{figure}

Figures~\ref{confronto_NN_vel}, \ref{confronto_NN_pressure}, and~\ref{confronto_NN_turbulent} compare the field reconstructions obtained with the three neural network architectures at a representative time instant ($Re = 100$, $N = 10$ POD modes): LSTM, MLP, and Transformer.

\begin{figure}[H]
\centering
\includegraphics[width=0.4\linewidth]{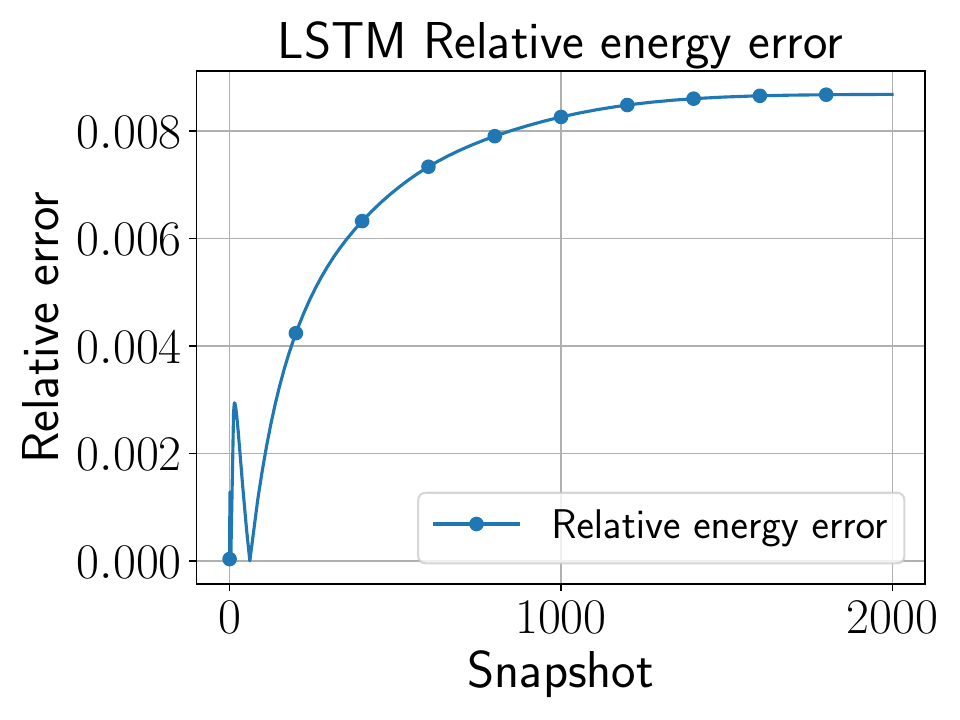}
\hfill
\includegraphics[width=0.4\linewidth]{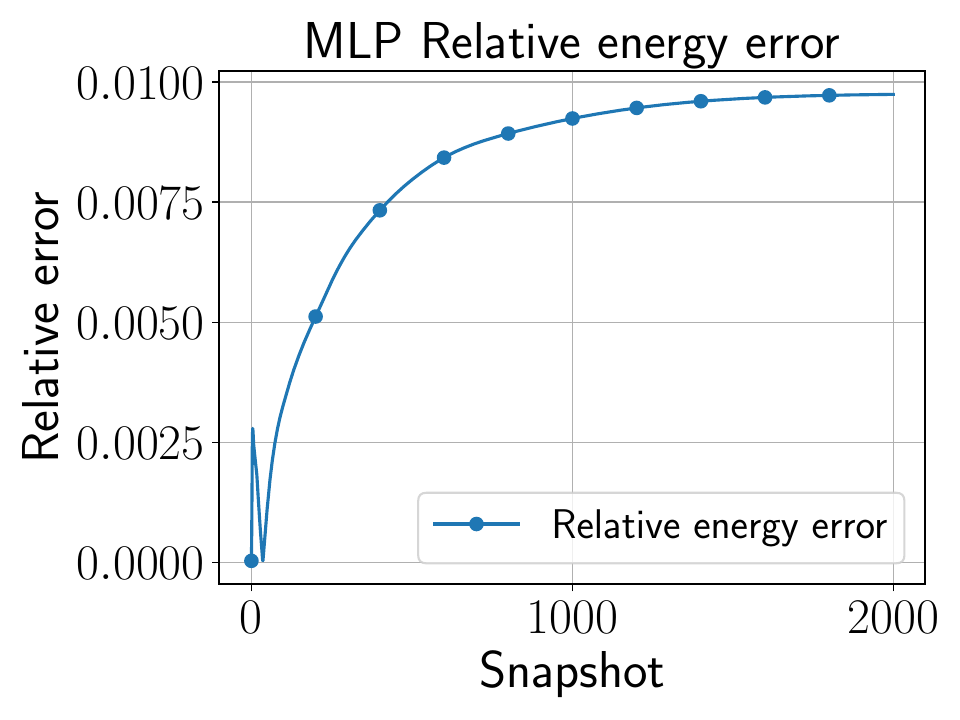}
\hfill
\includegraphics[width=0.4\linewidth]{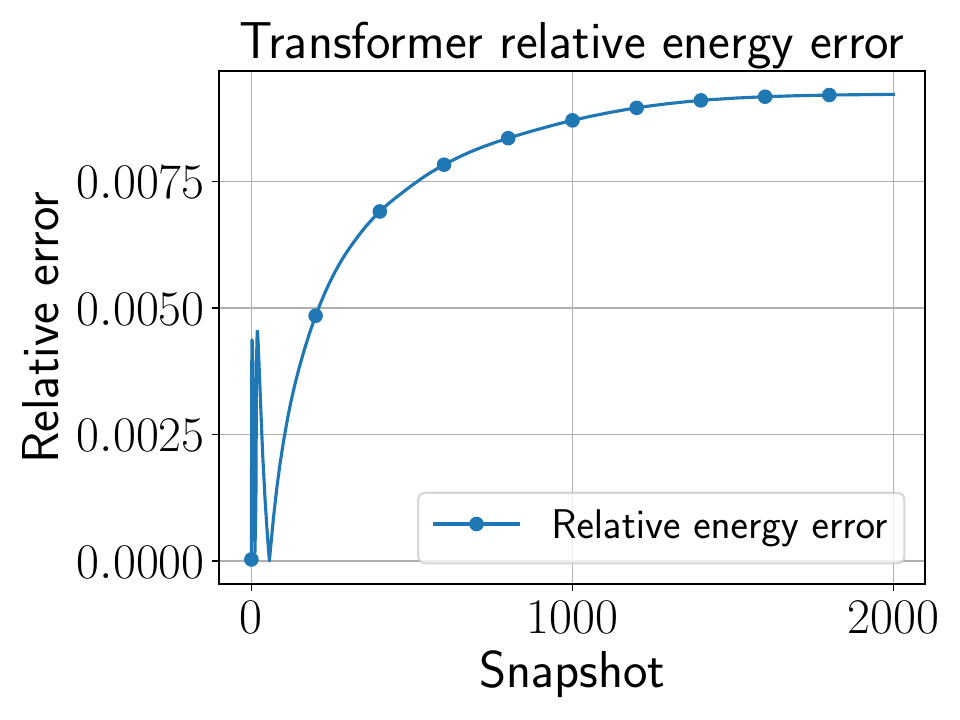}
\caption{\small Relative error (FOM vs ROM) in terms of kinetic energy as a function of time (Re = 100, \(N = 10\) POD modes): LSTM (left), MLP (center), and Transformer (right).}
\label{energy}
\end{figure}

\begin{figure}[H]
\centering
\includegraphics[width=0.4\linewidth]{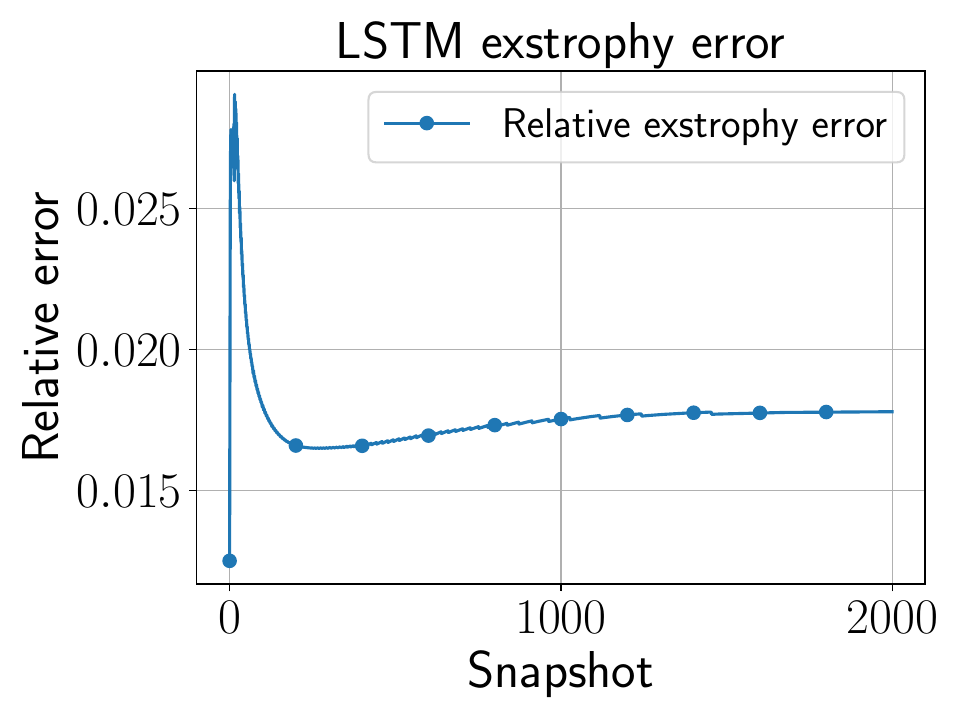}
\hfill
\includegraphics[width=0.4\linewidth]{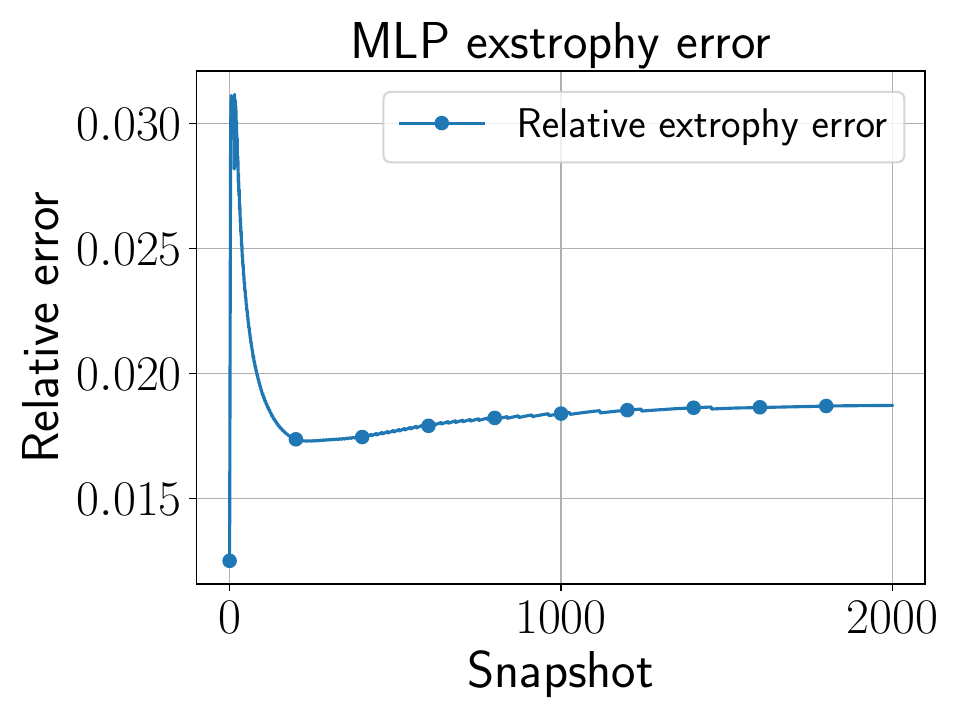}
\hfill
\includegraphics[width=0.4\linewidth]{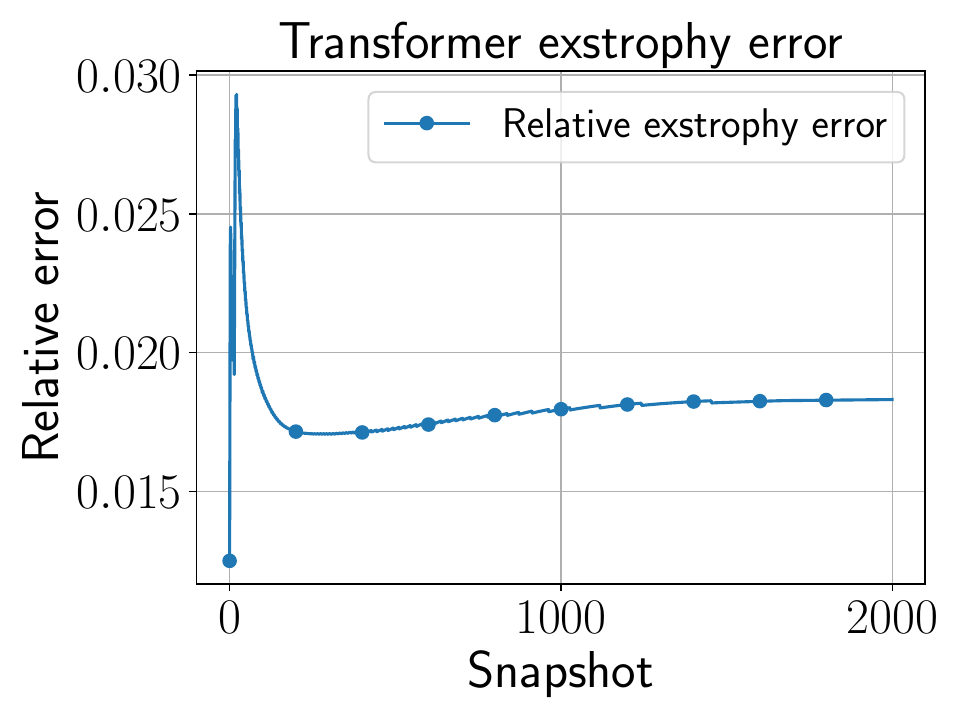}
\caption{\small Relative error (FOM vs ROM) in terms of enstrophy as a function of time (Re = 100, \(N = 10\) POD modes): LSTM (left), MLP (center), and Transformer (right).}
\label{enstrophy}
\end{figure}

Finally, we report in Figure~\ref{energy} the kinetic energy error and in Figure~\ref{enstrophy} the enstrophy error, obtained for all three architectures at $Re = 100$ with $N = 10$ POD modes. The enstrophy error, in particular, quantifies how well the ROM retains the vortical structures of the flow compared to the FOM.

The kinetic energy and enstrophy are defined, respectively, as
\begin{equation}
E(t) = \frac{1}{2}\int_\Omega |\mathbf{u}(t)|^2 \, dx,
\label{eqkineticenergy}
\end{equation}
\begin{equation}
\mathcal{E}(t) = \frac{1}{2}\int_\Omega |\boldsymbol{\omega}(t)|^2 \, dx,
\label{eqenstrophy}
\end{equation}
where $\boldsymbol{\omega} = \nabla \times \mathbf{u}$ is the vorticity field. The corresponding relative errors between the FOM and ROM solutions, reported in Figures \ref{energy} and \ref{enstrophy}, are computed as
\begin{equation}
\varepsilon_E(t) = \frac{\left| E^{FOM}(t) - E^{ROM}(t) \right|}{E^{FOM}(t)}, \qquad
\varepsilon_{\mathcal{E}}(t) = \frac{\left| \mathcal{E}^{FOM}(t) - \mathcal{E}^{ROM}(t) \right|}{\mathcal{E}^{FOM}(t)}.
\label{eqenergyenstrophyerror}
\end{equation}

In terms of energy and enstrophy, the three architectures exhibit similar performance at $Re = 100$, with the LSTM showing a marginal advantage consistent with its superior temporal modeling capability.

To quantify the computational advantage of the proposed ROM framework, wall-clock
timings were recorded for both the full-order and reduced-order simulations. The FOM (OpenFOAM LES, $50^3$ mesh, 2000 time steps at $\Delta t = 0.005$\,s) required
\textbf{458\,s} of wall-clock time on the workstation used in this study. The ROM online phase was timed for all three neural network architectures over 2000 time steps: 55.3\,s (LSTM), 53.7\,s (MLP), and 54.4\,s (Transformer), corresponding to speedups of $8.3\times$, $8.5\times$, and $8.4\times$ respectively over the FOM (458\,s). This speedup is realised entirely in the online phase and applies to every repeated evaluation of the same trained configuration, for example in optimization loops, sensitivity analyses, or ensemble post-processing, without any further FOM solve, which is the setting for which the offline/online decomposition is designed. To provide a complete accounting of the offline computational cost, we
report the measured neural-network training times for the final
$Re = 3200$ configuration ($N = 15$ modes and 1{,}400 training
samples), which is the higher-dimensional training configuration
considered in this study. Training was performed entirely on an
AMD Ryzen~9 7950X3D CPU workstation with 16 physical cores and
187\,GiB of system RAM; no GPU or other specialised accelerator was
used. The LSTM required approximately $167$\,s for $1{,}200$ epochs,
the MLP approximately $60$\,s for $1{,}000$ epochs, and the Transformer
approximately $203$\,s for $1{,}000$ epochs. The LSTM was trained separately, while the MLP and Transformer were
trained concurrently. The complete comparative training campaign
therefore required approximately $167 + \max(60, 203) = 370$\,s
($6.2$\,min) of elapsed wall-clock time.

For the final LSTM-based configuration, the complete one-time offline
cost was approximately $3702$\,s: $2935$\,s for the FOM simulation
(which includes a brief preliminary ROM verification in the same
binary run), approximately $600$\,s for POD basis construction
and reduced-operator assembly (estimated, performed in a separate
run), and approximately $167$\,s for LSTM training. Thus, LSTM
training accounts for approximately $4.5\%$ of the total offline
wall-clock cost and is approximately eighteen times shorter than
the FOM simulation. Training was completed on the same CPU
workstation used for the other computations and did not require
specialised hardware. Peak memory consumption during training was
not separately instrumented; however, no memory-capacity limitation
was encountered during any stage of the pipeline.

For $Q$ evaluations of the same trained configuration, the cumulative
FOM and ROM costs are respectively $Q T_{\mathrm{FOM}}$ and
$T_{\mathrm{offline}} + Q T_{\mathrm{online}}$. Their equality gives
the break-even number of evaluations
\begin{equation}\label{eq:breakeven}
Q^*
 = \frac{T_{\mathrm{offline}}}
        {T_{\mathrm{FOM}}-T_{\mathrm{online}}}
   = \frac{3702}{2935-58.9}
 \approx 1.29.
\end{equation}
Consequently, although a single FOM evaluation is less expensive than
constructing and running the ROM, the cumulative ROM workflow becomes
less expensive from the second evaluation onward.

All the results for the $Re = 100$ case are summarized in 
Table~\ref{tab:error_summary}. Wall-clock timings for the 
$Re = 3200$ generalization case are reported in 
Section~\ref{sec:Re3200} and Table~\ref{tab:Re3200_errors}.

\begin{table}[H]
    \centering
    \caption{Summary of relative errors, energy/enstrophy metrics, and computational cost between the FOM and ROM across different neural network architectures ($Re = 100$).}
    \label{tab:error_summary}

    \begin{tabular}{lccc} 
        \toprule
        & \textbf{LSTM} & \textbf{MLP} & \textbf{Transformer} \\
        \midrule
        \multicolumn{4}{l}{\textit{Reconstruction error at final time step ($N=10$ modes)}} \\
        \hspace{3mm} Velocity            & 0.7\% & 0.8\%  & 0.8\% \\
        \hspace{3mm} Pressure            & 6.0\% & 6.0\%  & 6.0\% \\
        \hspace{3mm} Turbulent viscosity & 2.0\% & 10.0\% & 7.0\% \\
        
        \midrule
        \multicolumn{4}{l}{\textit{Time-averaged reconstruction error ($N=10$ modes)}} \\
        \hspace{3mm} Velocity            & 0.7\% & 0.8\%  & 0.8\% \\
        \hspace{3mm} Pressure            & 6.0\% & 6.0\%  & 6.0\% \\
        \hspace{3mm} Turbulent viscosity & 4.0\% & 14.0\% & 9.0\% \\
        
        \midrule
        \multicolumn{4}{l}{\textit{Global Errors}} \\
        \hspace{3mm} Energy              & 0.7\% & 0.8\%  & 0.8\% \\
        \hspace{3mm} Enstrophy           & 1.7\% & 1.8\%  & 1.8\% \\

        \midrule
        \multicolumn{4}{l}{\textit{Computational Cost (2000 steps)}} \\
        \hspace{3mm} FOM wall-clock (s)  & \multicolumn{3}{c}{458} \\
        \hspace{3mm} ROM online wall-clock (s)  & 55.3 & 53.7 & 54.4 \\
        \hspace{3mm} Speedup             & $8.3\times$ & $8.5\times$ & $8.4\times$ \\
        
        \bottomrule
\end{tabular}%
\end{table}

\begin{figure}[H]
\centering
\includegraphics[width=0.70\linewidth]{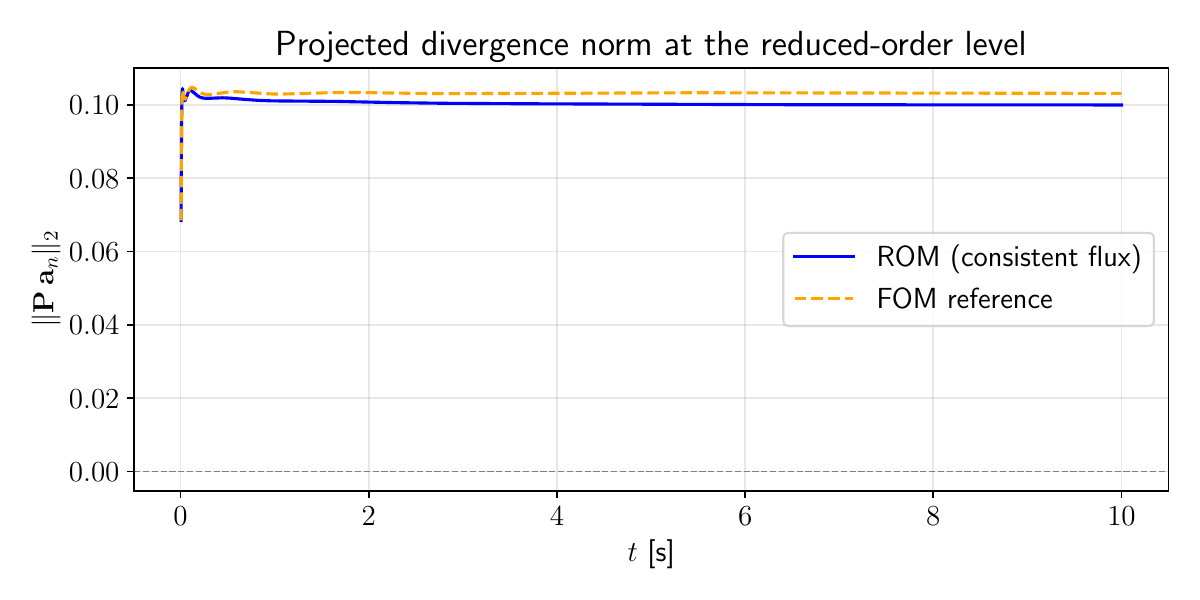}\\[4pt]
\includegraphics[width=0.70\linewidth]{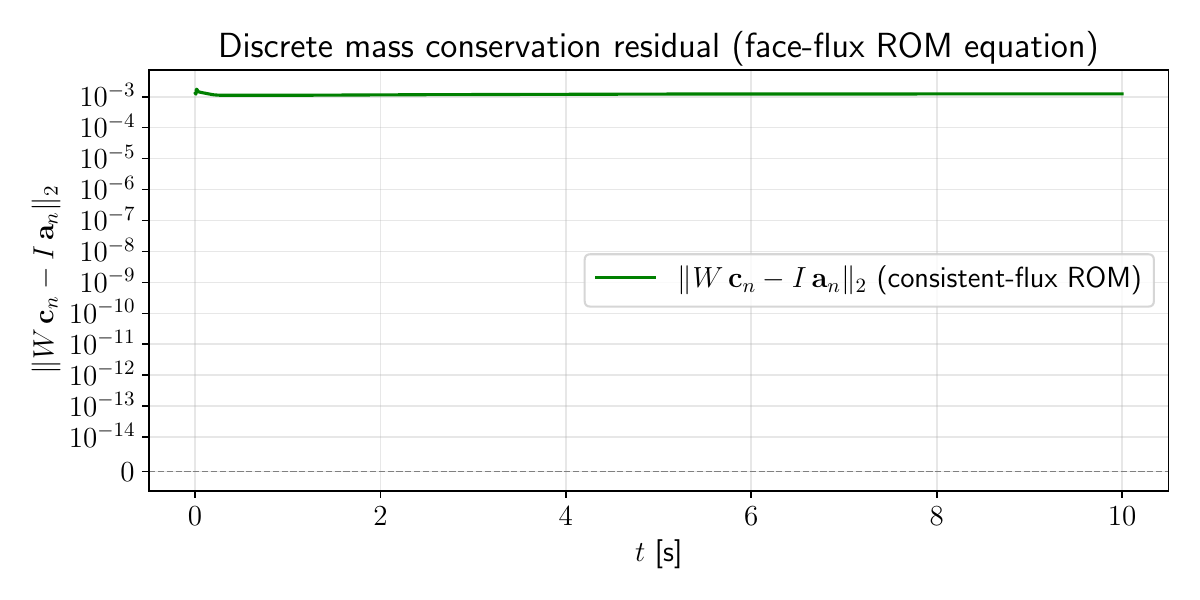}\\[4pt]
\includegraphics[width=0.70\linewidth]{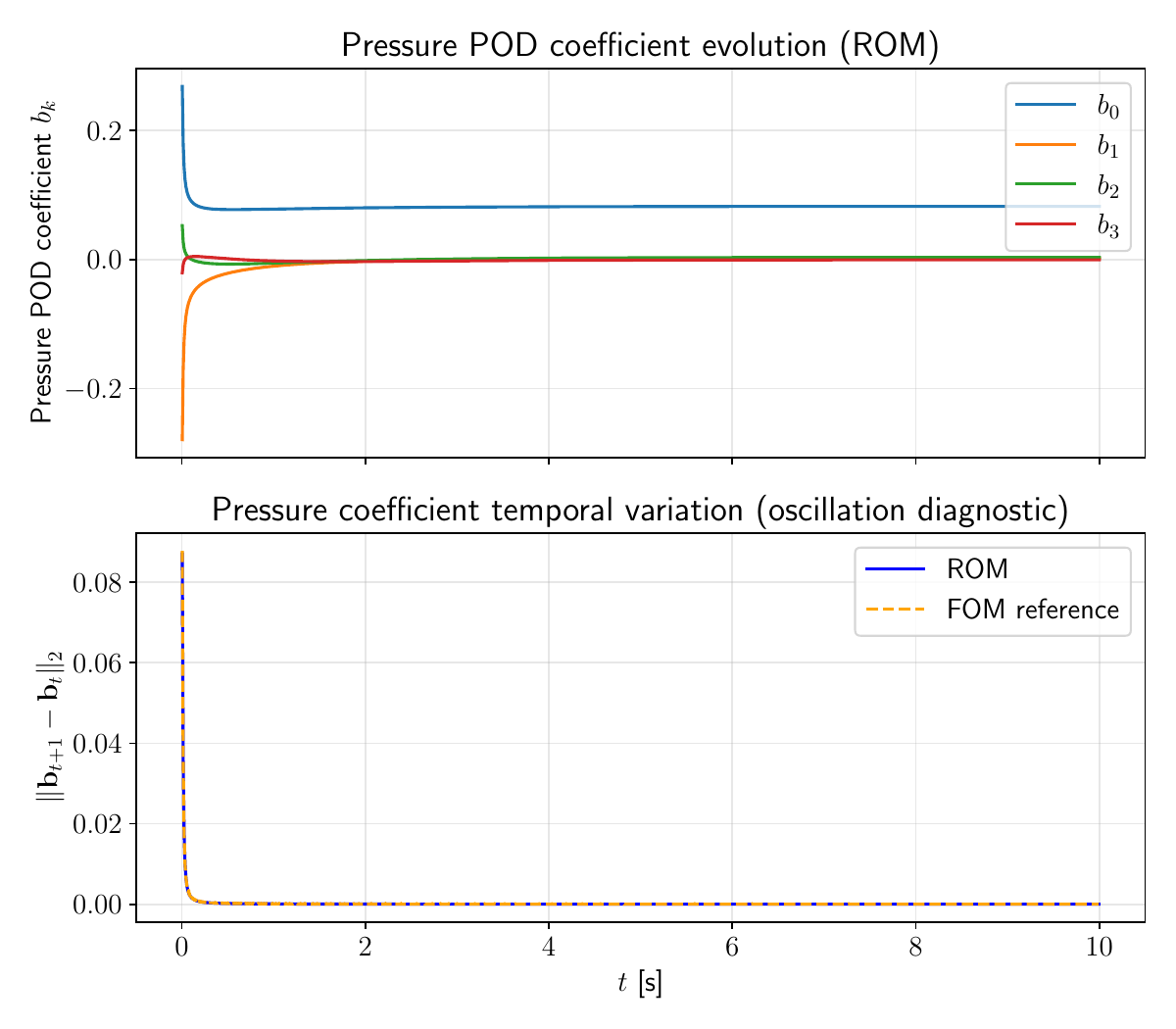}
\caption{\small Divergence and pressure diagnostics over 2000 time steps
($Re = 100$, $\Delta t = 0.005$\,s, $N = 10$ POD modes).
\textit{Top:} projected divergence norm $\|\mathbf{M}_r\,\mathbf{a}^n\|_2$
for ROM (solid blue) and FOM (dashed orange); curves within 3\%, confirming
no divergence amplification by the ROM time-marching.
\textit{Middle:} face-flux residual
$\|\mathbf{W}_r\,\mathbf{c}^n - \mathbf{N}_r\,\mathbf{a}^n\|_2$,
bounded and non-growing, confirming stable discrete mass conservation.
\textit{Bottom (upper):} pressure POD coefficients $b_0$--$b_3$;
smooth evolution rules out checkerboard instability.
\textit{Bottom (lower):} step-to-step pressure variation
$\|\mathbf{b}^{n+1}-\mathbf{b}^n\|_2$ for ROM and FOM;
statistically identical signals confirm no pressure oscillation amplification.}
\label{fig:div_diagnostics}
\end{figure}

\paragraph{Divergence and pressure diagnostics.}
To demonstrate that the discretize-then-project formulation controls
pressure--velocity decoupling at the reduced level, three complementary
diagnostics were computed over all 2000 time steps of the online simulation
and are reported in Figure~\ref{fig:div_diagnostics}.

The \emph{projected divergence norm} $\|\mathbf{M}_r\,\mathbf{a}^n\|_2$
(equivalently $\|\mathbf{P}\,\mathbf{a}^n\|_2$ in
Figure~\ref{fig:div_diagnostics}, where $\mathbf{P} \equiv \mathbf{M}_r$
is the notation used in the diagnostic script)
monitors the driving term of the reduced pressure Poisson
equation~\eqref{eq:rom_p1}. Expanding the definition,
\begin{equation*}
    \mathbf{M}_r\,\mathbf{a}^n
    = \mathbf{X}^T\mathbf{V}\,\mathbf{M}_p\,
      \underbrace{(\mathbf{\Phi}\,\mathbf{a}^n)}_{\approx\,\mathbf{u}_p^n},
\end{equation*}
where $\mathbf{\Phi}\,\mathbf{a}^n$ reconstructs the cell-centred velocity
field from the POD coefficients, $\mathbf{M}_p$ applies the discrete
divergence operator to that reconstructed field, and $\mathbf{X}^T\mathbf{V}$
projects the result onto the pressure POD basis.
In a perfectly divergence-free ROM this quantity would be identically zero;
any non-zero value is the residual divergence that the pressure Poisson
solver must correct at the next step.
At each time step the solver first evaluates $\mathbf{M}_r\,\mathbf{a}^n$,
uses it as the source term to solve for the new pressure $\mathbf{b}^{n+1}$
via equation~\eqref{eq:rom_p1}, and then applies the resulting pressure
gradient to correct the velocity to $\mathbf{a}^{n+1}$;
the diagnostic therefore captures the divergence \emph{before} the
pressure correction, at the moment when it is largest.
Because $\mathbf{M}_r\,\mathbf{a}^n \in \mathbb{R}^{N_r^p}$ is a
10-dimensional vector, the $L_2$ norm reduces it to a single scalar per
time step, so Figure~\ref{fig:div_diagnostics} (top panel) shows one
number per step over the full 2000-step simulation.
As shown there, the ROM curve (consistent-flux) tracks the FOM reference
throughout the simulation to within 3\%
(ROM mean $1.00\times10^{-1}$, FOM mean $1.03\times10^{-1}$).
This residual is not introduced by the ROM time-marching; it is an intrinsic
property of the truncated POD basis, ten velocity modes cannot span a
perfectly divergence-free subspace on a collocated finite-volume grid, and the ROM inherits this level from the FOM basis without amplifying it.

The \emph{face-flux mass conservation residual}
$\|\mathbf{W}_r\,\mathbf{c}^n - \mathbf{N}_r\,\mathbf{a}^n\|_2$
monitors equation~\eqref{eq:rom_uf1}. Physically, the condition
$\mathbf{W}_r\,\mathbf{c} = \mathbf{N}_r\,\mathbf{a}$ is the
reduced-level statement of discrete mass conservation for the
face fluxes: it expresses that the face-flux POD coefficients
$\mathbf{c}^n$ and the cell-centred velocity coefficients
$\mathbf{a}^n$ must remain mutually consistent.
A naive Galerkin ROM without the face-flux equation would have
no mechanism to enforce this condition, allowing the residual
to grow freely and eventually destabilise the pressure field.
The residual here remains bounded in
$[1.1\times10^{-3},\,1.7\times10^{-3}]$ with no growth trend
across all 2000 steps. The $\mathcal{O}(10^{-3})$ magnitude is
physically expected: $\mathbf{c}^n$ is computed using the
velocity at the previous step $\mathbf{a}^{n-1}$, while the
residual is evaluated against the pressure-corrected
$\mathbf{a}^n$, so the gap is
$\mathcal{O}(\Delta t) = \mathcal{O}(5\times10^{-3})$.
The absence of growth confirms that the consistent-flux formulation
stably enforces discrete mass conservation throughout the simulation.

The \emph{pressure oscillation diagnostic} uses two complementary
views of the pressure POD coefficient vector $\mathbf{b}^n$, which
is the output of equation~\eqref{eq:rom_p1} at every step.

The upper sub-panel of the bottom plot in Figure~\ref{fig:div_diagnostics} shows the time
evolution of the first four individual pressure coefficients
$b_0, b_1, b_2, b_3$ over all 2000 steps. In the presence of a
checkerboard instability, these curves would exhibit rapid
high-frequency alternations between successive steps, a direct
visual signature of pressure--velocity decoupling at adjacent
degrees of freedom. The smooth, monotonically settling behaviour
observed here confirms that no such instability is present in the
consistent-flux ROM.

The bottom panel provides the corresponding quantitative measure:
the step-to-step variation $\|\mathbf{b}^{n+1}-\mathbf{b}^n\|_2$,
which collapses all 10 pressure coefficients into a single scalar
per step. A growing or persistently large signal here would be the
reduced-level signature of the checkerboard instability, in which
pressure and velocity become decoupled at adjacent degrees of
freedom and the pressure field develops spurious oscillations.

The ROM and FOM signals remain closely matched throughout
(ROM mean $1.80\times10^{-4}$, FOM mean $1.86\times10^{-4}$;
ROM max $8.71\times10^{-2}$, FOM max $8.78\times10^{-2}$),
confirming that the ROM introduces no pressure oscillation
amplification relative to the FOM reference. The shared initial
transient peak is present in both FOM and ROM and corresponds
to the early cavity dynamics, ruling out a numerical artefact
of the reduced-order scheme.

\subsection{Generalization to Higher Reynolds Number ($Re = 3200$)}
\label{sec:Re3200}
To further validate the proposed framework and address the question of 
near-steady behavior observed at $Re = 100$, the full offline--online 
workflow was repeated for a second test case at $Re = 3200$ 
($\nu = 3.125 \times 10^{-4}$\,m$^2$/s, $U_{\text{lid}} = 1$\,m/s), 
representing a flow regime above the onset of turbulent characteristics 
in the three-dimensional lid-driven cavity~\cite{Leriche2000}, exhibiting 
sustained unsteady dynamics and non-negligible turbulent viscosity $\nu_T$ 
throughout the entire simulation window. The same $50\times50\times50$ 
structured mesh, boundary conditions, and Smagorinsky LES model 
($C_s = 0.2$) are retained from the $Re = 100$ case. At $Re = 3200$ 
the subgrid-scale model accounts for a larger fraction of the turbulent 
dissipation than at $Re = 100$, which is consistent with the engineering 
LES approach adopted throughout this work and does not affect the 
validity of the ROM framework being evaluated.

\begin{figure}[H]
\centering
\includegraphics[width=0.6\linewidth]{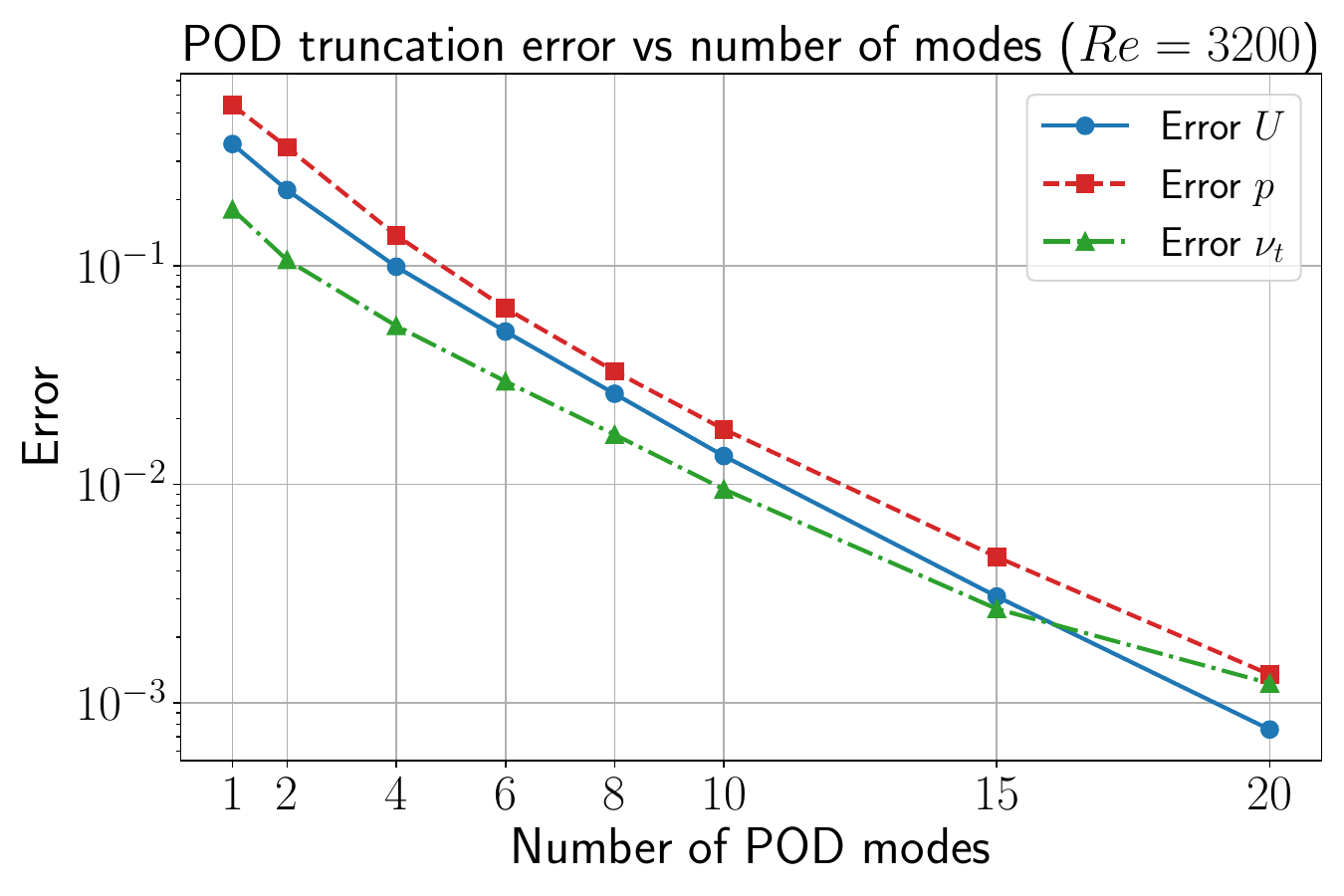}
\caption{\small POD truncation error as a function of the number of 
retained modes for the $Re = 3200$ test case. The error metric is the 
time-averaged relative $L^2$ norm:
$\varepsilon(n) = \langle \|\mathbf{u}_{\text{FOM}} - 
\mathbf{u}_{\text{POD},n}\|_2 / \|\mathbf{u}_{\text{FOM}}\|_2 \rangle_t$.
At $N = 10$ modes the truncation errors are $1.35\%$ (velocity), 
$1.78\%$ (pressure), and $0.95\%$ (turbulent viscosity). The curve is 
still declining at $N = 10$, reflecting the greater complexity of the 
$Re = 3200$ flow compared to $Re = 100$.}
\label{fig:pod_errors_Re3200}
\end{figure}

The LES simulation was advanced for $T = 4$\,s with internal solver step 
$\Delta t = 10^{-4}$\,s and snapshots stored every 
$\Delta t_{\text{write}} = 2 \times 10^{-3}$\,s, yielding $N_s = 2000$ 
snapshots. The FOM wall-clock time on the same workstation was 
\textbf{2935\,s}. POD was then applied to extract $N = 10$ modes per 
field. The time-averaged relative $L^2$ POD truncation errors at 
$N = 10$ are:
\begin{equation}
    \varepsilon_U = 1.35 \times 10^{-2}, \quad
    \varepsilon_p = 1.78 \times 10^{-2}, \quad
    \varepsilon_{\nu_T} = 9.49 \times 10^{-3},
\end{equation}
as shown in Figure~\ref{fig:pod_errors_Re3200}. 

We note that the POD truncation error curve for $Re = 3200$ is still declining at $N = 10$ modes,
reflecting the greater complexity of the flow at higher Reynolds number. By contrast, the
ROM online error at $Re = 100$ had effectively plateaued at $N = 10$ (Figure~\ref{errors_POD_modes_variation}).
Therefore, $N = 10$ is first evaluated as a baseline for direct comparison with the $Re = 100$ case.
A systematic mode sensitivity study is
then performed for $Re = 3200$ later in this section, showing that $N = 15$ is the best overall configuration and
is therefore adopted for the final $Re = 3200$ results.

The relative $L^2$ velocity and pressure reconstruction errors at the 
final time step are summarised in Table~\ref{tab:Re3200_errors}. The 
velocity error reaches approximately $13.9\%$ (LSTM), $13.8\%$ (MLP), 
and $13.7\%$ (Transformer) at $t = 4$\,s. The pressure reconstruction error is approximately $40\%$ across all three architectures, larger than at $Re = 100$ ($\sim 6\%$), reflecting the increased complexity of the pressure field at higher Reynolds number. Importantly, the pressure error is nearly identical across architectures because it is governed primarily by the reduced pressure Poisson solve rather than directly by the neural-network closure. We note that the pressure POD coefficients have a very small $L^2$ norm (mean $\approx 0.013$ at $Re = 3200$, compared to $\approx 0.080$ at $Re = 100$), which amplifies the relative error metric; the absolute pressure reconstruction error remains small in both cases (below $0.25\%$ when normalized by a unit-offset denominator), and the elevated relative figures reflect normalization by a near-zero denominator rather than poor absolute accuracy.

The pressure reference-cell parameter $\tau$ in 
Equation~\eqref{eq:reference_pressure_constraint} was optimized via a 
direct reconstruction-error sweep over a range of candidate values. This 
reduced the pressure error from $65\%$ for $\tau = 100{,}000$ to $40\%$ 
for $\tau = 15{,}000$, showing that the appropriate strength of the 
reference-cell constraint must be selected based on the closed-loop ROM 
reconstruction error rather than on conditioning alone. The POD truncation 
error for pressure at $N = 10$ modes is only $1.78\%$, confirming that 
the dominant error source is the online ROM solve rather than the basis 
representation.

The mean energy and enstrophy errors remain below $8\%$ across all three 
architectures, indicating that the ROM correctly captures the dominant 
flow physics despite the elevated pressure reconstruction error. Notably, 
despite exhibiting the largest overfitting ratio between training and 
validation loss, the LSTM architecture achieves the lowest turbulent 
viscosity reconstruction error among the three models, suggesting that 
temporal memory provides a physical advantage in the ROM time-stepping 
that is not captured by standard loss metrics alone.

At $N = 15$ modes, the validation losses increase further for all 
three architectures (LSTM: $1.44$, MLP: $0.69$, Transformer: $0.86$), 
consistent with the higher-dimensional closure mapping being harder 
to learn with the same network capacity. Nevertheless, the LSTM 
reconstruction accuracy remains comparable ($\nu_T$ error $16.94\%$ 
versus $17.05\%$ at $N = 10$), further confirming that raw loss metrics 
do not reliably predict ROM reconstruction accuracy. This motivates 
future work on neural-network regularization, larger training datasets, 
and sequence-aware validation strategies. Nevertheless, the final assessment 
in this work is based on closed-loop ROM rollout against the FOM, rather 
than on one-step neural-network validation loss alone.

\begin{figure}[H]
\centering
\begin{subfigure}[b]{0.48\textwidth}
    \centering
    \includegraphics[width=\linewidth]{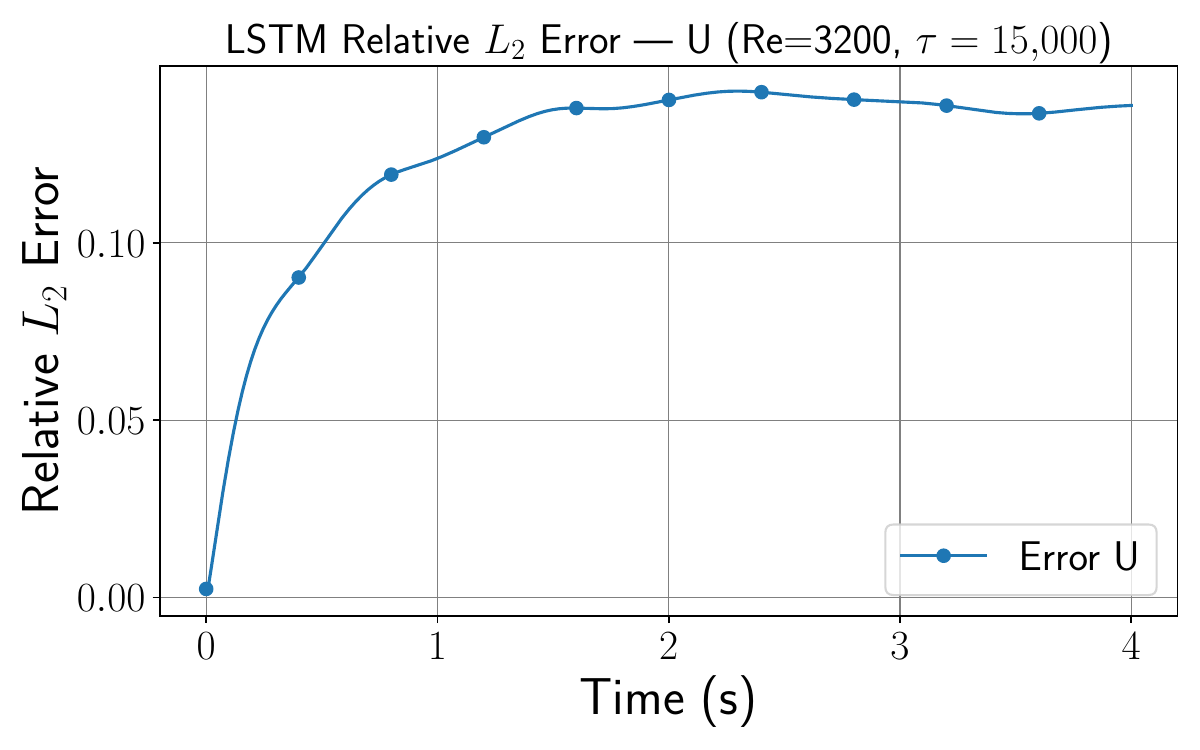}
    \caption{Velocity error $e_U(t)$}
    \label{fig:error_U_Re3200_python}
\end{subfigure}
\hfill
\begin{subfigure}[b]{0.48\textwidth}
    \centering
    \includegraphics[width=\linewidth]{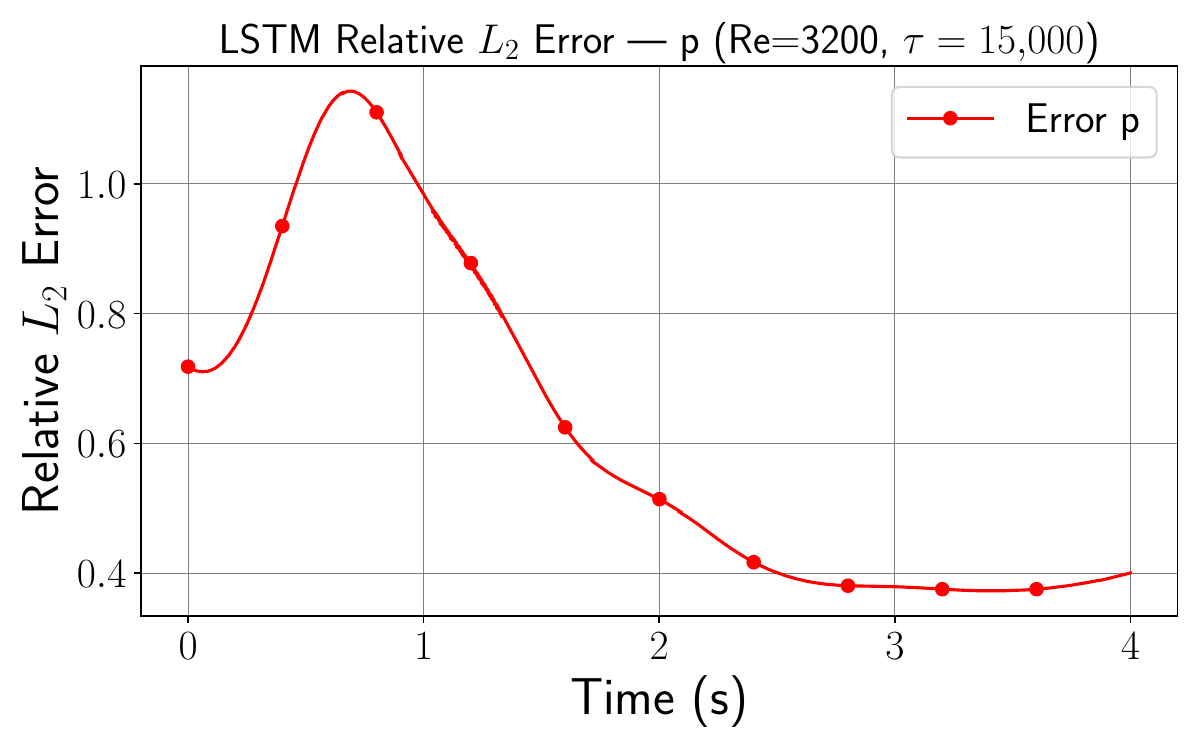}
    \caption{Pressure error $e_p(t)$}
    \label{fig:error_p_Re3200_python}
\end{subfigure}
\caption{\small Relative $L^2$ errors between FOM and ROM solutions
over 2000 time steps for the $Re = 3200$ test case ($N = 10$ POD modes, LSTM closure, $\tau = 15{,}000$). Left: velocity error rises from zero and plateaus
at approximately $13.9\%$ by $t \approx 1.5$\,s, remaining bounded
and stable throughout. Right: pressure error peaks at approximately
$114\%$ during the initial transient ($t \approx 0.8$\,s), then
decreases steadily to approximately $40\%$ at $t = 4$\,s. The
pressure error is governed by the online ROM solve rather than the
POD basis (POD truncation error: $1.78\%$) and remains bounded
with no sign of divergence.}
\label{fig:errors_Re3200_python}
\end{figure}

\begin{table}[H]
\centering
\caption{Reconstruction errors and online-phase wall-clock timings for the $Re = 3200$ test case at $N = 10$ ($\tau = 15{,}000$) and $N = 15$ ($\tau = 1{,}000$) POD modes. U error and p error are relative $L^2$ 
errors at the final time step ($t = 4$\,s). Energy and enstrophy errors 
are time-averaged means. FOM wall-clock: 2935\,s.}
\label{tab:Re3200_errors}
\begin{tabular}{llccccccc}
\toprule
$N$ & Architecture & U error & p error & $\nu_T$ error & 
Energy & Enstrophy & ROM online time (s) & Speedup \\
\midrule
\multirow{3}{*}{10} 
 & LSTM        & 13.88\% & 40.04\% & 17.05\% & 4.55\% & 7.42\% & 66.8 & $43.9\times$ \\
 & MLP         & 13.79\% & 40.18\% &  9.48\% & 4.76\% & 7.63\% & 66.7 & $44.0\times$ \\
 & Transformer & 13.65\% & 40.21\% & 12.43\% & 4.74\% & 7.62\% & 69.1 & $42.5\times$ \\
\midrule
\multirow{3}{*}{\textbf{15}}
 & \textbf{LSTM}        & \textbf{12.34\%} & \textbf{29.76\%} & \textbf{16.94\%} & \textbf{8.13\%} & \textbf{10.17\%} & \textbf{58.9} & $\mathbf{49.9\times}$ \\
 & \textbf{MLP}         & \textbf{12.77\%} & \textbf{29.66\%} & \textbf{34.22\%} & \textbf{8.04\%} & \textbf{10.06\%} & \textbf{61.7} & $\mathbf{47.5\times}$ \\
 & \textbf{Transformer} & \textbf{12.37\%} & \textbf{29.75\%} & \textbf{21.43\%} & \textbf{8.15\%} & \textbf{10.20\%} & \textbf{67.2} & $\mathbf{43.7\times}$ \\
\bottomrule
\end{tabular}
\end{table}

Comparing the POD truncation errors (Figure~\ref{fig:pod_errors_Re3200})
with the ROM online errors (Figure~\ref{fig:rom_error_vs_modes_Re3200})
reveals that the truncation error at $N = 10$ ($1.78\%$ for pressure)
is more than twenty times smaller than the corresponding ROM online
error ($40\%$), confirming that the dominant source of reconstruction
error at $Re = 3200$ is the online ROM solve rather than the POD basis.
Increasing the mode count from $N = 10$ to $N = 15$ reduces the online
pressure error by approximately $10$ percentage points (from $40\%$ to
$30\%$), while the truncation error decreases only marginally (from
$1.78\%$ to $0.46\%$), further confirming that the online solve
remains the bottleneck throughout the mode count range studied. 
While the energy and enstrophy errors increase slightly at $N = 15$
compared to $N = 10$ (see Table~\ref{tab:Re3200_errors}), the velocity and pressure reconstruction errors
improve substantially, and $N = 15$ remains the preferred configuration
given that velocity and pressure accuracy are the primary metrics for
ROM validation in this context.

\subsubsection*{Sensitivity to the Number of POD Modes}
To determine the optimal number of POD modes for the $Re = 3200$ case,
a systematic study was performed by varying $N \in \{10, 15, 20, 25\}$
and running the full ROM online phase (with LSTM, MLP, and Transformer
closures) for each value. For each $N$, the pressure reference-cell parameter
$\tau$ in Equation~\eqref{eq:reference_pressure_constraint} was re-optimized via a
direct reconstruction-error sweep, since the appropriate strength of the
reference-cell constraint depends on both the Reynolds number and the retained pressure basis. The resulting optimal
values exhibit a non-monotonic trend across $N$:

\begin{table}[H]
\centering
\caption{Optimal pressure reference-cell parameter $\tau$ and LSTM reconstruction errors as a function of the number of
retained POD modes ($Re = 3200$).}
\label{tab:tau_vs_N}
\begin{tabular}{ccccc}
\toprule
$N$ & Optimal $\tau$ & LSTM $e_U$ & LSTM $e_p$ & LSTM $e_{\nu_T}$ \\
\midrule
10 & 15,000 & 13.88\% & 40.04\% & 17.05\% \\
\textbf{15} & \textbf{1,000} & \textbf{12.34\%} & \textbf{29.76\%} & \textbf{16.94\%} \\
20 & 0     & 12.47\% & 30.73\% & 17.11\% \\
25 & 10    & 11.44\% & 32.55\% & 17.89\% \\
\bottomrule
\end{tabular}
\end{table}

The optimal $\tau$ decreases from $15{,}000$ at $N = 10$ to $1{,}000$
at $N = 15$ and to $0$ at $N = 20$, indicating that the strength needed
to enforce the pressure reference-cell constraint depends on the retained
pressure basis. At $N = 25$, however, setting $\tau = 0$ leads to an
ill-conditioned reduced pressure solve (pressure error $122\%$), and a
small nonzero value ($\tau = 10$) is required to impose the pressure
gauge consistently. This reveals a non-monotonic relationship between
mode count and the strength of the reference-cell constraint.

\begin{figure}[H]
\centering
\includegraphics[width=0.7\linewidth]{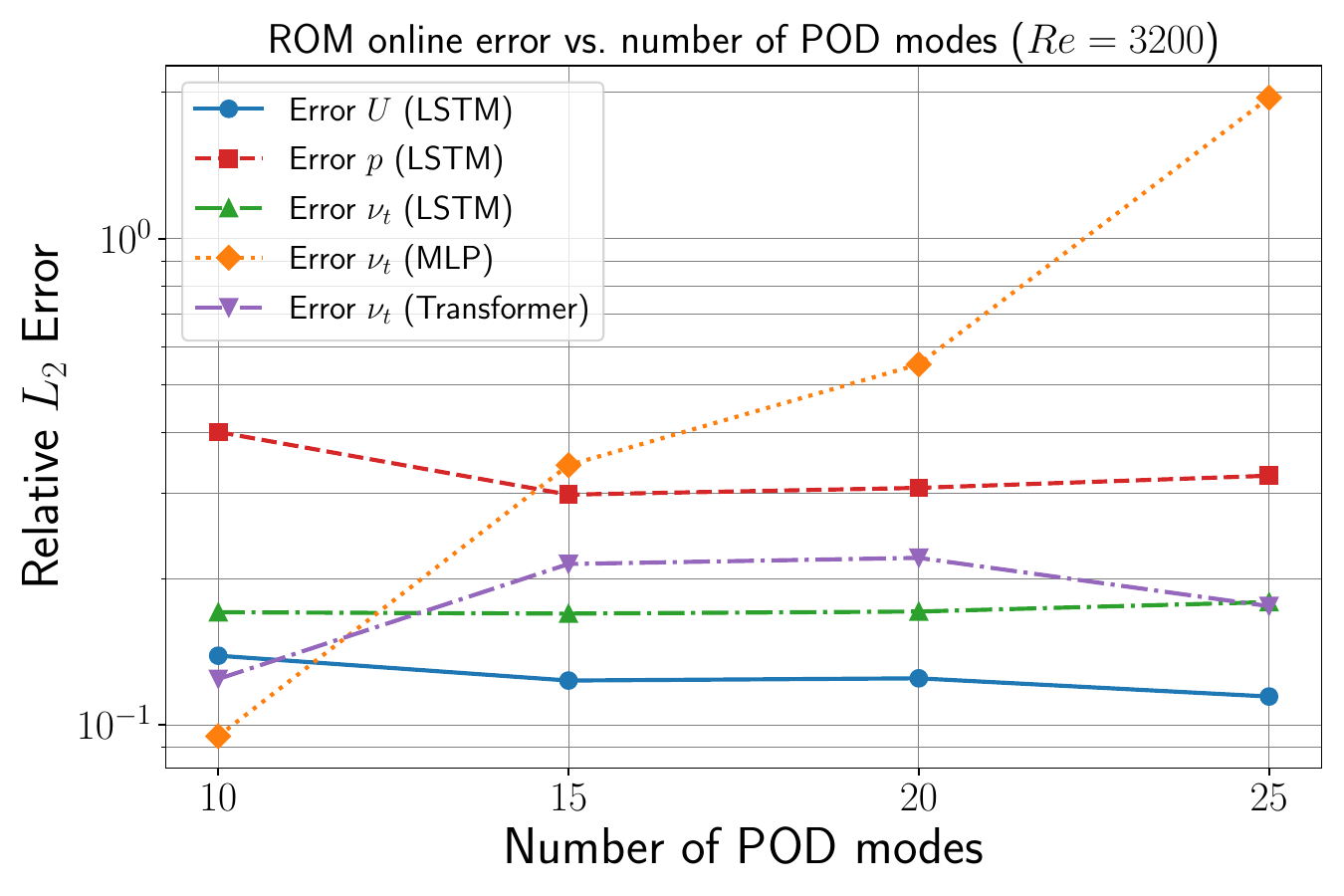}
\caption{\small ROM online error as a function of the number of retained 
POD modes for the $Re = 3200$ test case. The relative $L^2$ error at the 
final time step is shown for velocity $\mathbf{u}$ and pressure $p$ 
(LSTM closure) and for turbulent viscosity $\nu_T$ across all three 
architectures. The LSTM $\nu_T$ error remains stable across the tested mode
counts ($\sim 17\%$), whereas the MLP $\nu_T$ error increases sharply
at $N = 25$ ($195\%$). This result indicates a loss of robustness for
the present MLP configuration at the highest tested dimension, although
the available experiment does not isolate whether the cause is network
capacity, optimization sensitivity, input scaling, or the absence of
temporal context. The Transformer exhibits non-monotonic intermediate
behavior. The velocity and pressure errors (LSTM) improve from $N = 10$ to $N = 15$ and remain
stable thereafter, identifying $N = 15$ as the best overall mode count for the $Re = 3200$ case with the current
architecture.}
\label{fig:rom_error_vs_modes_Re3200}
\end{figure}

The sensitivity of the ROM online error to the number of retained POD modes is shown in Figure~\ref{fig:rom_error_vs_modes_Re3200},
confirming $N = 15$ as the best overall choice for the $Re = 3200$ case.

The mode sensitivity study identifies $N = 15$ as the best overall 
configuration. Although $N = 25$ gives the lowest velocity error 
($11.44\%$), $N = 15$ gives the lowest pressure error ($29.76\%$), 
nearly the same velocity accuracy ($12.34\%$), a lower reduced dimension, 
and avoids the closure degradation observed at higher $N$. Relative to 
the $N = 10$ baseline, $N = 15$ improves the velocity error by $1.54$ 
percentage points and the pressure error by $10.28$ percentage points. Beyond $N = 15$, further increasing the mode count does not improve the overall error balance:
although the velocity error decreases slightly at $N = 25$, the pressure error increases at
$N = 20$ and $N = 25$. This indicates that the online ROM solve, rather than the POD basis alone,
remains the dominant limitation for the pressure reconstruction.

A striking architecture-dependent behaviour is also observed in the
turbulent viscosity closure: the LSTM $\nu_T$ error remains stable
across all mode counts ($\sim 17\%$), while the MLP error increases
monotonically from $9.48\%$ at $N = 10$ to $195\%$ at $N = 25$,
effectively diverging. The Transformer exhibits intermediate behaviour,
degrading at $N = 15$--$20$ before recovering at $N = 25$ ($17.55\%$).
The MLP divergence at $N = 25$ ($195\%$ $\nu_T$ error) can be 
attributed to an unfavorable ratio between model capacity and 
output dimensionality. The MLP has $11{,}594$ trainable parameters 
and must learn a mapping from $2N = 50$ inputs to $N = 25$ outputs 
— a significantly harder regression problem than at $N = 10$ 
($20$ inputs, $10$ outputs), yet with the same network capacity. 
Without temporal context to exploit structure in the coefficient 
sequences, the MLP must learn this mapping from a relatively small 
training set (1{,}400 samples), leading to effective underfitting 
of the higher-dimensional closure problem and eventual divergence.

The Transformer's recovery at $N = 25$ ($17.55\%$, comparable to 
the LSTM) is notable given its degradation at $N = 15$--$20$. 
This suggests that the self-attention mechanism, which captures 
pairwise relationships across the lookback window of 15 time steps, 
becomes more effective as the output dimensionality grows — 
potentially because the higher-dimensional $\nu_T$ field exhibits 
richer temporal correlations that the attention mechanism can 
exploit, while the fixed-capacity MLP cannot. This non-monotonic 
behaviour warrants further investigation in future work. Overall, the LSTM provides the most consistent behaviour across the full mode-count study, while
the Transformer shows partial recovery at higher mode count and the MLP becomes unreliable as the
closure dimension increases. Based on this analysis, $N = 15$ is adopted as
the final configuration for the $Re = 3200$ case.

\begin{figure}[H]
\centering
\begin{subfigure}[b]{0.48\textwidth}
    \centering
    \includegraphics[width=\linewidth]{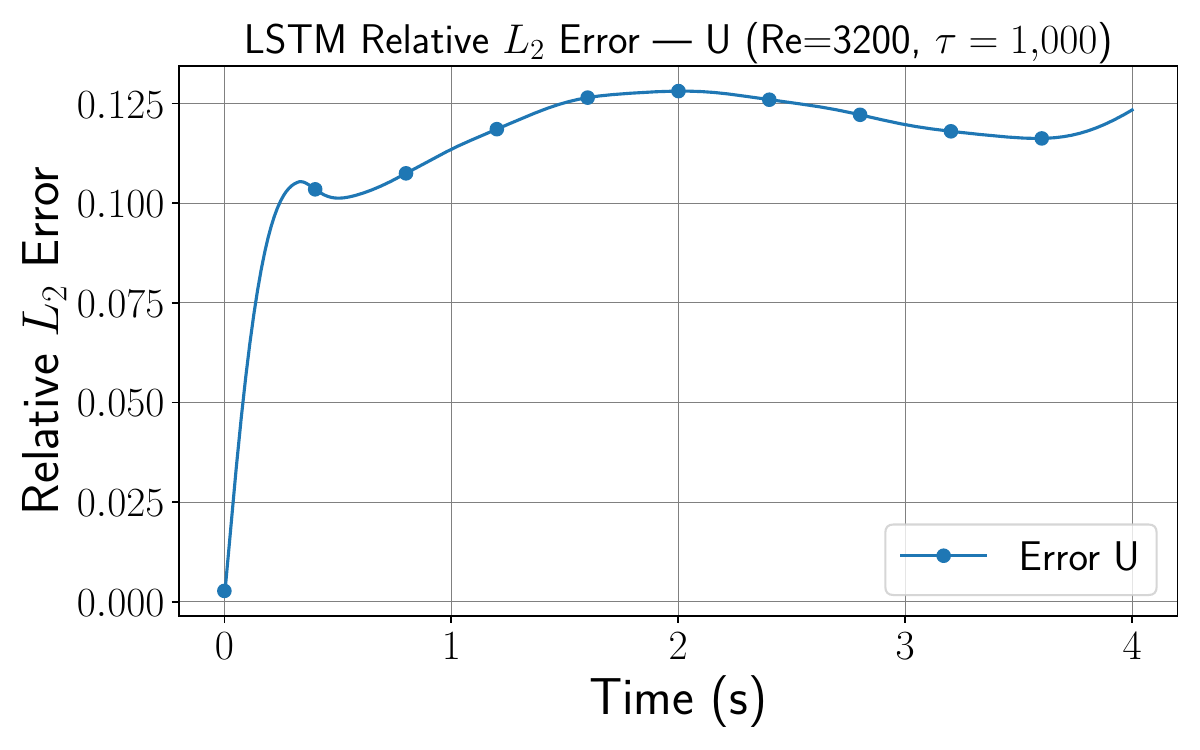}
    \caption{Velocity error $e_U(t)$}
    \label{fig:error_U_Re3200_N15}
\end{subfigure}
\hfill
\begin{subfigure}[b]{0.48\textwidth}
    \centering
    \includegraphics[width=\linewidth]{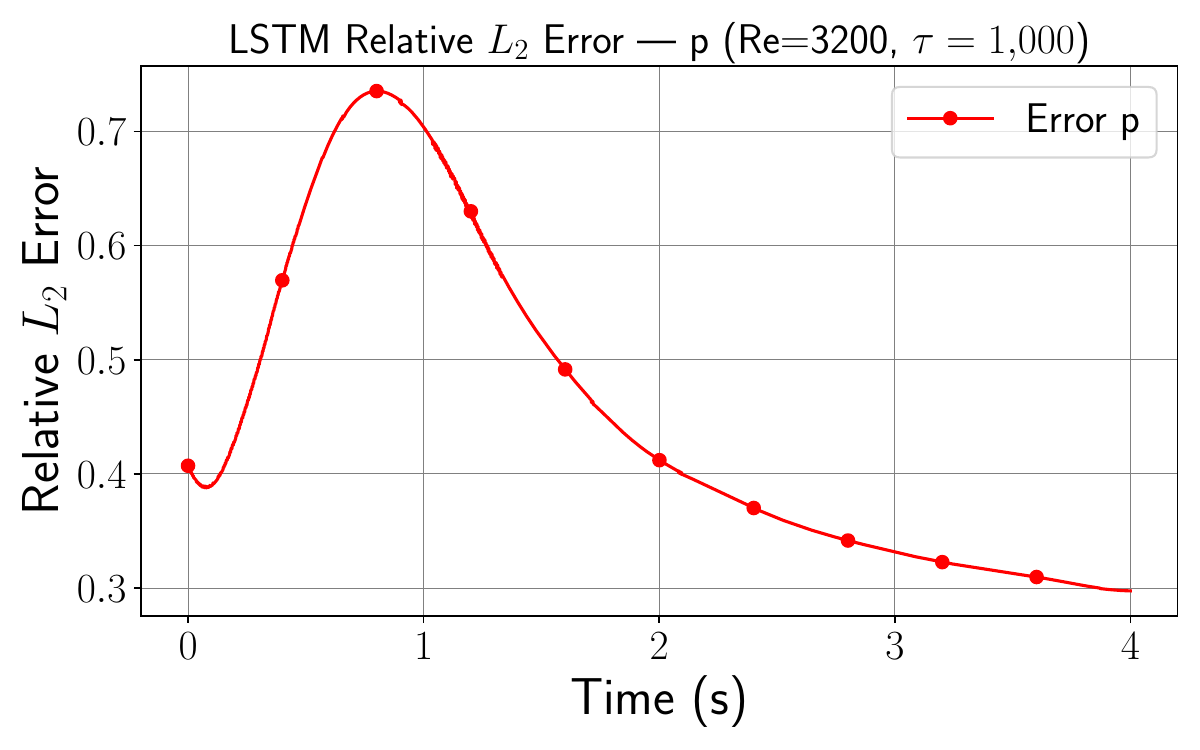}
    \caption{Pressure error $e_p(t)$}
    \label{fig:error_p_Re3200_N15}
\end{subfigure}
\caption{\small Relative $L^2$ errors between FOM and ROM solutions
over 2000 time steps for the $Re = 3200$ test case at the best overall configuration ($N = 15$ POD modes, LSTM closure, $\tau = 1{,}000$).
Left: velocity error plateaus at approximately $12.3\%$ by
$t \approx 1$\,s, remaining bounded and stable throughout.
Right: pressure error peaks at approximately $74\%$ during the
initial transient ($t \approx 0.8$\,s), then decreases steadily
to approximately $30\%$ at $t = 4$\,s, an improvement over the
$N = 10$ configuration ($40\%$). Both errors remain bounded with
no sign of divergence.}
\label{fig:errors_Re3200_N15}
\end{figure}

\subsection*{Time Extrapolation Test}

To assess the temporal generalization capability of the framework, 
a time extrapolation test was performed at $Re = 3200$: the POD 
basis, neural network closure, and pressure reference-cell parameter were 
all derived exclusively from snapshots in $t \in [0, 3]$\,s, 
and the ROM was then rolled out over the full window 
$t \in [0, 6]$\,s, with $t \in [3, 6]$\,s representing 
genuinely unseen temporal dynamics within the same Reynolds number and geometry.

\begin{figure}[H]
\centering
\includegraphics[width=0.75\linewidth]{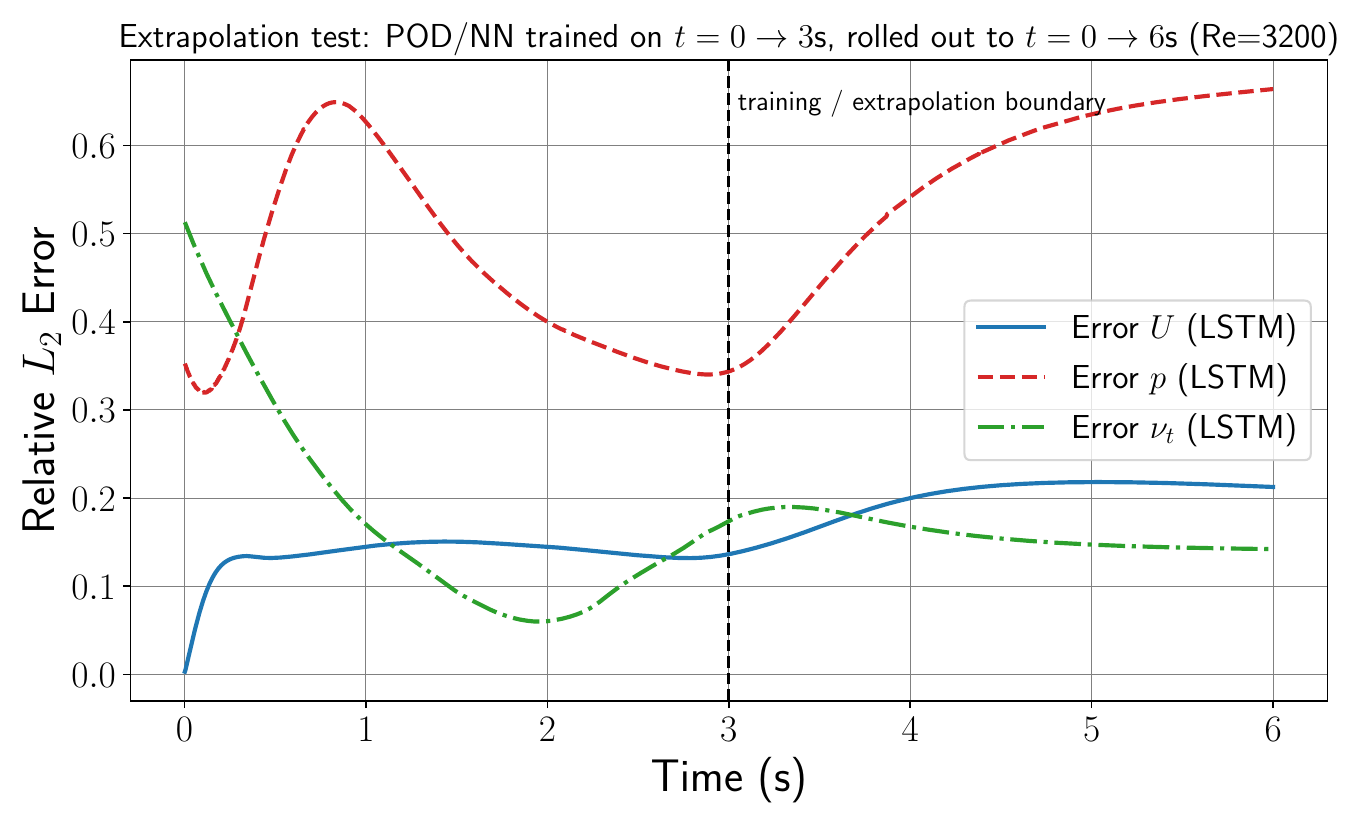}
\caption{\small Time extrapolation test at $Re = 3200$ ($N = 15$ 
POD modes, LSTM closure, $\tau = 1{,}000$). The POD basis and 
neural network were trained on $t \in [0, 3]$\,s only; the ROM 
was rolled out to $t = 6$\,s. The vertical dashed line marks the 
training/extrapolation boundary. Velocity error ($e_U$, solid blue) 
rises from a training-window mean of $13.6\%$ to an extrapolation 
mean of $19.9\%$; pressure error ($e_p$, dashed red) rises from 
$45.5\%$ to $56.3\%$, still increasing at $t = 6$\,s. By contrast, 
the turbulent-viscosity error ($e_{\nu_T}$, dash-dot green) 
decreases from a training-window mean of $17.7\%$ to $16.1\%$ in 
the extrapolation region, indicating that the LSTM closure remains 
robust beyond the training window even as the velocity and pressure 
fields degrade. The ROM remains bounded throughout without 
numerical divergence.}
\label{fig:extrapolation}
\end{figure}

Figure~\ref{fig:extrapolation} shows the relative $L^2$ 
reconstruction errors for velocity and pressure over the full 
rollout window. Within the training region ($t < 3$\,s), the velocity error 
averages $13.6\%$, the pressure error averages $45.5\%$, and the 
turbulent-viscosity error averages $17.7\%$. These training-window 
errors are slightly higher than the main $Re=3200$ results at 
$N=15$ ($12.34\%$ and $29.76\%$ for velocity and pressure 
respectively), reflecting the reduced snapshot count used for POD 
construction ($1{,}500$ vs $2{,}001$ snapshots) and the shorter 
training window ($t \in [0,3]$\,s vs $t \in [0,4]$\,s). In the extrapolation 
region ($t > 3$\,s), the velocity error rises to an average of 
$19.9\%$ and the pressure error to $56.3\%$, reaching $21.3\%$ 
and $66.4\%$ respectively at the final time step ($t = 6$\,s), 
with no sign of leveling off.

Notably, the turbulent-viscosity error 
decreases slightly in the extrapolation region (mean $16.1\%$), 
indicating that the LSTM closure remains robust beyond the training 
window even as the velocity and pressure fields degrade. The energy and enstrophy errors further confirm this picture: 
in the training region ($t < 3$\,s), the mean kinetic energy error 
is $11.1\%$ and the mean enstrophy error is $13.5\%$; in the 
extrapolation region ($t > 3$\,s), these rise to $30.9\%$ and 
$31.7\%$ respectively, roughly tripling relative to the training 
window. This amplification, more pronounced than the per-field 
velocity and pressure errors, arises because energy and enstrophy 
are quadratic aggregates over all POD modes, making them more 
sensitive to the accumulated drift in the extrapolation region. The ROM 
remains bounded and does not diverge, but the systematic error 
increase confirms that the current framework has limited temporal 
extrapolation capability beyond the training window. This highlights time extrapolation as an important direction for 
future improvement. Promising avenues to address this limitation 
include online adaptation of the closure model during rollout, 
physics-informed neural network architectures that enforce 
conservation constraints beyond the training window, and reservoir computing approaches combined with nonlinear
autoencoder-based dimensionality reduction for reduced-order forecasting
of dynamical systems
\cite{hajisharifi2025combining}.

At the best overall configuration of $N = 15$ modes, the online-phase wall-clock speedup 
reaches approximately $50\times$ over the FOM (ROM: $\sim 59$\,s 
vs.\ FOM: 2935\,s for 2000 steps), notably faster than the $N = 10$ 
configuration ($\sim 44\times$, ROM: $\sim 67$\,s). This improvement 
arises because the smaller pressure reference-cell parameter at $N = 15$ 
($\tau = 1{,}000$ vs.\ $\tau = 15{,}000$ at $N = 10$) yields a 
better-conditioned reduced pressure solve, reducing the cost of the 
online phase. The speedup at $N = 15$ is substantially larger than 
at $Re = 100$ ($\sim 8\times$), because the FOM cost grows 
substantially with the smaller time step required at higher Reynolds 
number, while the ROM online cost remains modest across both cases.
These results confirm that the proposed discretize-then-project framework 
with data-driven closure extends to a substantially more demanding 
flow regime without changing the neural-network architecture, after 
case-specific retraining and pressure reference-cell parameter retuning ($\tau$, $\mathbf{v}$).

\section{Conclusions and outlooks} \label{sec:Conclusions}

This work introduces a novel hybrid framework that combines projection-based ROM with data-driven approaches for the treatment of turbulent viscosity. Specifically, the reduced order models for the velocity and pressure variables have been constructed following a “discretize-then-project” strategy ensuring full consistency with the underlying discrete formulation of the incompressible Navier–Stokes equations. At the same time, the turbulent viscosity field has been reconstructed through machine learning techniques, enabling a flexible and efficient integration of physical and data-driven components within the reduced framework.

A central outcome of this study concerns the role of the neural network architecture adopted for the modeling of turbulent viscosity within the reduced-order setting. Based on the numerical results at $Re = 100$, the LSTM architecture shows more reliable performance than the MLP and Transformer models, consistent with its temporal memory mechanism. At $Re = 3200$, the velocity and pressure errors are comparable across architectures, while the turbulent-viscosity reconstruction is more architecture-sensitive. In particular, the LSTM provides the most stable closure behavior across the mode-count study. This behavior suggests that explicitly accounting for temporal correlations plays an important role when modeling unsteady turbulent dynamics, particularly in situations where the evolution of the closure term depends on the recent flow history.

The hybrid model was tested by comparing reduced-order results with the corresponding full-order simulations. For the configuration considered in this work, the reduced model captures the main characteristics of the flow evolution. This indicates that the combination of projection-based reduction and data-driven modeling can be used for turbulent viscosity prediction in this setting. Furthermore, for the $Re = 100$ case, the ROM online phase achieves an online-phase wall-clock speedup of $8.3\times$ over the full-order solver (55.3\,s vs.\ 458\,s for 2000 time steps), confirming the practical computational benefit of the offline/online decomposition in a repeated-query scenario.

The framework was additionally validated at $Re = 3200$, a genuinely 
unsteady regime above the onset of turbulent characteristics, where the 
turbulent viscosity field is non-negligible and time-dependent, in contrast 
to the near-steady behavior observed at $Re = 100$. A systematic sensitivity 
study identified $N = 15$ POD modes as the best overall configuration. 
Although $N = 25$ gives a slightly lower velocity error, $N = 15$ provides 
the lowest pressure error, nearly the same velocity accuracy, a lower reduced 
dimension, and more stable closure behavior. With this configuration, the 
velocity error is reduced from $13.9\%$ at $N = 10$ to $12.3\%$, and the 
pressure error is reduced from $40\%$ to $30\%$, with an online-phase 
wall-clock speedup of approximately $50\times$ at this higher Reynolds number.

The pressure 
reconstruction error is attributable to the online pressure Poisson solve 
rather than the POD basis (POD truncation error: $1.78\%$ at $N = 10$), 
as discussed in Section~\ref{sec:Re3200}. These results demonstrate that the methodology can be extended 
to a more demanding unsteady regime after case-specific retraining and pressure reference-cell parameter retuning, 
while highlighting pressure reconstruction at higher Reynolds number as an important direction for future 
refinement. A time extrapolation test confirms that the ROM remains bounded when rolled out to $t = 6$\,s after training on $t \in [0, 3]$\,s only. Velocity and pressure errors increase in the unseen region ($19.9\%$ and $56.3\%$ respectively), while the LSTM turbulent-viscosity closure remains robust (error decreasing from $17.7\%$ to $16.1\%$), motivating future work on adaptive and physics-informed closure strategies. Extension to cross-parameter generalization, training at one Reynolds 
number and rolling out at another without retraining, is identified 
as a natural direction for future work.

As discussed in Section~\ref{sec:ROM}, the framework is geometry-agnostic 
by construction, and extension to more complex geometries is planned as a 
direction of future work.

\section*{Acknowledgments}
We acknowledge the financial support under the National Recovery and Resilience Plan (NRRP), Mission 4, Component 2, Investment 1.1, Call for tender No. 1409 published on 14.9.2022 by the Italian Ministry of University and Research (MUR), funded by the European Union – NextGenerationEU– Project Title ROMEU – CUP P2022FEZS3 - Grant Assignment Decree No. 1379 adopted on 01/09/2023 by the Italian Ministry of Ministry of University and Research (MUR) and acknowledges the financial support by the European Union (ERC, DANTE, GA-101115741). Views and opinions expressed are however those of the author(s) only and do not necessarily reflect those of the European Union or the European Research Council Executive Agency. Neither the European Union nor the granting authority can be held responsible for them.

\bibliographystyle{ieeetr}
\bibliography{sample}

\section*{Appendix} \label{appendix}

The numerical simulations were performed using the ITHACA-FV 
library based on OpenFOAM. Specifically, the unsteady 
Navier-Stokes cases were set up using the solver available 
in the ITHACA-FV library.

The hyperparameters presented in Table~\ref{tab:hyperparameters} 
define the technical configuration for the non-intrusive, 
data-driven component of the hybrid framework. While velocity 
and pressure are resolved through intrusive Galerkin projection, 
these network architectures are designed to reconstruct the 
temporal evolution of the turbulent viscosity coefficients.

To ensure a robust learning process, the dataset was partitioned 
differently for each test case. For the $Re = 100$ case, 1,800 
training samples and 200 validation samples were drawn from a total 
of 2,001 high-fidelity snapshots. For the $Re = 3200$ case, the 
first 500 snapshots ($t < 1$\,s) were excluded to avoid training 
on the startup transient; the remaining snapshots were split into 
1,400 training samples (snapshots 500--1900) and 100 validation 
samples (snapshots 1900--2000). For the time extrapolation test 
(Section~\ref{sec:Re3200}), a restricted split was used: 850 
training samples (snapshots 500--1350, $t \in [1, 2.7]$\,s) and 
150 validation samples (snapshots 1350--1500, $t \in [2.7, 3]$\,s).

All models utilized 
the Adam optimizer and Mean Squared Error (MSE) loss function, 
which are standard for regressing modal coefficients in fluid 
dynamics applications. Numerical stability across the different 
architectures was maintained by applying a StandardScaler to 
normalise both input (velocity and pressure coefficients) and 
output (turbulent viscosity coefficients) sequences.

For the sequence-based models (LSTM and Transformer), a lookback 
window of 15 time steps was adopted. This choice allows the 
networks to incorporate information from recent temporal history 
when predicting the turbulent viscosity coefficients, which is 
relevant for representing the time-dependent behavior observed 
in the unsteady turbulent flow.

\begin{table}[H]
    \centering
\caption{Hyperparameters and architecture details of the 
    neural network models used in this study. Training and validation sample counts, input/output dimensions, and the
reported trainable-parameter counts refer to the $Re = 100$ case
($N = 10$ modes). For 
    $Re = 3200$ with $N = 15$ modes: 1,400 training and 
    100 validation samples; input dimension $15 \times 30$ 
    (LSTM/Transformer) or $30$ (MLP); output dimension $15$.
    (see Section~\ref{sec:Re3200}).}
\label{tab:hyperparameters}
    \begin{tabular}{lccc} 
        \toprule
        \textbf{Hyperparameters} & \textbf{LSTM} & \textbf{MLP} & \textbf{Transformer} \\
        \midrule
        \multicolumn{4}{l}{\textit{General Training}} \\
        \hspace{3mm} Epochs             & 1200 & 1000 & 1000 \\
        \hspace{3mm} Training Samples   & 1800 & 1800 & 1800 \\
        \hspace{3mm} Validation Samples & 200  & 200  & 200 \\
        \hspace{3mm} Batch Size         & 64   & 64   & 64 \\
        \hspace{3mm} Optimizer          & Adam & Adam & Adam \\
        \hspace{3mm} Loss Function      & MSE  & MSE  & MSE \\
        \hspace{3mm} Learning Rate      & $2 \times 10^{-5}$ & $3 \times 10^{-5}$ & $9 \times 10^{-5}$ \\
        \hspace{3mm} Normalization      & StandardScaler & StandardScaler & StandardScaler \\

        \midrule
        \multicolumn{4}{l}{\textit{Architecture Specifics}} \\
        \hspace{3mm} Lookback Window    & 15 & -- & 15 \\
        \hspace{3mm} Hidden Layers      & 4  & 4  & -- \\
        \hspace{3mm} Input Dimension    & $15 \times 20$ ($N=10$); $15 \times 30$ ($N=15$) 
        & $20$ ($N=10$); $30$ ($N=15$) 
        & $15 \times 20$ ($N=10$); $15 \times 30$ ($N=15$) \\
        \hspace{3mm} Output Dimension   & $10$ ($N=10$); $15$ ($N=15$) 
        & $10$ ($N=10$); $15$ ($N=15$) 
        & $10$ ($N=10$); $15$ ($N=15$) \\

        \midrule
        \multicolumn{4}{l}{\textit{Transformer Specifics}} \\
        \hspace{3mm} Attention Heads    & -- & -- & 4 \\
        \hspace{3mm} Feed Forward Dim   & -- & -- & 128 \\
        \hspace{3mm} Dropout            & -- & -- & 0.1 \\
        \hspace{3mm} Layer Norm $\epsilon$ & -- & -- & $1 \times 10^{-6}$ \\

        \midrule
        \multicolumn{4}{l}{\textit{Model Complexity}} \\
        \hspace{3mm} Trainable Parameters ($N=10$)  & $35{,}562$ & $11{,}594$ & $26{,}010$ \\
        \hspace{3mm} Final Train Loss ($Re=100$)   & $1.25 \times 10^{-5}$ & $1.94 \times 10^{-2}$ & $3.17 \times 10^{-4}$ \\
        \hspace{3mm} Final Val Loss ($Re=100$)     & $1.53 \times 10^{-5}$ & $3.69 \times 10^{-4}$ & $7.84 \times 10^{-5}$ \\
        \hspace{3mm} Final Train Loss ($Re=3200$, $N=10$)  & $8.33 \times 10^{-6}$ & $8.50 \times 10^{-3}$ & $9.64 \times 10^{-5}$ \\
\hspace{3mm} Final Val Loss ($Re=3200$, $N=10$)    & $2.23 \times 10^{-1}$ & $8.21 \times 10^{-2}$ & $1.31 \times 10^{-1}$ \\
\hspace{3mm} Final Train Loss ($Re=3200$, $N=15$)  & $2.14 \times 10^{-5}$ & $1.23 \times 10^{-2}$ & $1.36 \times 10^{-4}$ \\
\hspace{3mm} Final Val Loss ($Re=3200$, $N=15$)    & $1.44$ & $6.90 \times 10^{-1}$ & $8.64 \times 10^{-1}$ \\
        \bottomrule
    \end{tabular}
\end{table}

The internal configurations adopted for the three neural 
network models considered in this study, LSTM, MLP, and 
Transformer, are reported in Table~\ref{tab:layer_features}. 
The layer-wise dimensions in Table~\ref{tab:layer_features} correspond to the $N=10$ 
configuration; for the $Re=3200$, $N=15$ case, the same hidden-layer hyperparameters 
are retained, while the input and output dimensions increase consistently with the 
larger number of POD coefficients. 
All three architectures are employed within the non-intrusive 
component of the hybrid framework and are used to model the 
temporal evolution of the turbulent viscosity coefficients in 
the reduced-order space.

The LSTM model is designed to exploit temporal dependencies 
through its gated memory mechanism. The architecture utilizes 
a hierarchical recurrent structure: the first layer (64 units) 
is configured to return sequences, passing temporal information 
from the lookback window to the subsequent 32-unit LSTM layer, 
followed by a dense layer with ReLU activation and a linear output layer predicting the $N$ reduced turbulent-viscosity
coefficients, with $N=10$ for the architecture dimensions reported
in Table~\ref{tab:layer_features} and $N=15$ for the final $Re=3200$ configuration.

The MLP architecture consists of a sequence of fully connected 
layers without explicit temporal recurrence. The network starts 
with a dense layer of 128 units to capture global relationships 
in the input features, followed by a Dropout layer (rate 0.2) 
to limit overfitting, and terminates with a linear output layer producing the $N$ predicted
turbulent-viscosity coefficients, with $N=10$ for the architecture
dimensions reported in Table~\ref{tab:layer_features} and $N=15$ for the final $Re=3200$
configuration.

The Transformer architecture employs a self-attention mechanism 
to capture temporal dependencies across the lookback window. 
It consists of two encoder blocks, each containing a 
multi-head attention layer (4 heads, key dimension 20), 
dropout (rate 0.1), layer normalization, and a feed-forward 
network (dimension 128). The encoder output is passed through 
a dense layer (64 units, ReLU) before the final linear output layer producing the $N$
turbulent-viscosity coefficients, with $N=10$ for the architecture
dimensions reported in Table~\ref{tab:layer_features} and $N=15$ for the final $Re=3200$
configuration.

\begin{table}[H]
    \centering
    \caption{Layer-wise architecture details for the LSTM, 
    MLP, and Transformer models.}
    \label{tab:layer_features}
    \begin{tabular}{lllc} 
        \toprule
        \textbf{Model} & \textbf{Layer} & \textbf{Output Shape} 
        & \textbf{Parameters / Config} \\
        \midrule
        
        \textbf{LSTM} & 1. LSTM (64 units) & $15 \times 64$ 
                      & Return Sequences \\
                      & 2. LSTM (32 units) & $32$           
                      & No Return Sequences \\
                      & 3. Dense (32 units) & $32$          
                      & ReLU Activation \\
                      & 4. Output (Dense)   & $10$     
                      & Linear Activation \\
        \midrule
        
        \textbf{MLP}  & 1. Dense (128 units) & $128$        
                      & ReLU \\
                      & 2. Dropout           & $128$        
                      & Rate $= 0.2$ \\
                      & 3. Dense (64 units)  & $64$         
                      & ReLU \\
                      & 4. Output (Dense)    & $10$    
                      & Linear Activation \\
        \midrule

        \textbf{Transformer} 
                      & 1. MultiHeadAttention $\times 2$ & $15 \times 20$ 
                      & 4 heads, key\_dim$=20$ \\
                      & 2. Dropout $\times 4$  & $15 \times 20$        
                      & Rate $= 0.1$ \\
                      & 3. LayerNormalization $\times 4$ & $15 \times 20$         
                      & $\epsilon = 1\times10^{-6}$ \\
                      & 4. FF Dense $\times 2$ & $15 \times 128 \to 15 \times 20$    
                      & ReLU $\to$ Linear \\
                      & 5. Dense (64 units)  & $64$         
                      & ReLU \\
                      & 6. Output (Dense)    & $10$    
                      & Linear Activation \\
        
        \bottomrule
    \end{tabular}
\end{table}

\end{document}